\documentclass[vecarrow]{svmult}

\usepackage{amsmath,amssymb}
\usepackage{float}

\usepackage{makeidx}         %
\usepackage{graphicx}        %
\usepackage{multicol}        %
\usepackage[bottom]{footmisc}%

\makeindex             %
\newcommand{\circdown}[1]{\mathop\circ\limits_{\rlap{$#1$}}}
\newlength{\barlength}
\newcommand{\minusfill}{$\mathsurround=0pt\mathord- \mkern-6mu
    \cleaders\hbox{$\mkern-3mu \mathord- \mkern-3mu$}\hfill
       \mkern-6mu \mathord-$}
\newcommand{\yokobo}{\hbox to 2em{\minusfill}}
\newcommand{\equalfill}{$\mathsurround=0pt\mathord= \mkern-6mu
    \cleaders\hbox{$\mkern-3mu \mathord= \mkern-3mu$}\hfill
       \mkern-6mu \mathord=$}
\newcommand{\Longlongrightarrow}
       {\hbox to 2em{\equalfill$\mkern-3mu\Rightarrow$}}
\newcommand{\Longlongleftarrow}
       {\hbox to 2em{$\Leftarrow\mkern-3mu$\equalfill}}
\newcommand{\tsume}{\kern-.35em}
\newlength{\circlength}
\divide\circlength by 2
\advance\circlength by 0.08em

\newcommand{\tatebo}{\raisebox{\circlength}{\rotatebox{90}{\yokobo}}}
\newcommand{\Ad}{\operatorname{Ad}}
\def\adm#1{{$#1$-admissible}}
\newcommand{\rrank}{\mathbb{R}\text{-}\operatorname{rank}}

\newcommand{\gk}{$(\mathfrak{g}_\mathbb{C}, K)$}
\def   \set#1#2{\{{#1}:{#2}\}}
\newcommand{\rarrowsim}{\smash{\mathop{\,\longrightarrow\,}\limits
  ^{\lower1.5pt\hbox{$\scriptstyle\sim$}}}}
\newcommand{\diag}{\operatorname{diag}}
\def \F#1#2{\pi_{#2}^{#1}}
\def \itm#1{\newline\noindent{\rm{#1}}\enspace}
\newcommand{\End}{\operatorname{End}}
\newcommand{\Hom}{\operatorname{Hom}}
\newcommand{\Aut}{\operatorname{{Aut}}}
\newcommand{\rt}{\Delta}
\def \hwm#1#2{\pi^{{#1}}_{{#2}}}
\def\Sb #1 \endSb{_{\substack{#1}}}
\newenvironment{eq-text}{\par%
\refstepcounter{equation}
\noindent
{\upshape(\theequation)}\quad}
{\par%
\noindent \ignorespacesafterend}
\spnewtheorem{claimsubsec}{Claim}{\bfseries}{\itshape}

\spnewtheorem{thmalph}{Theorem}{\bfseries}{\itshape}

\spnewtheorem{fact}{Fact}[subsection]{\bfseries}{\itshape}
\spnewtheorem{rem}[fact]{Remark}{\itshape}{}
\spnewtheorem{lem}[fact]{Lemma}{\bfseries}{\itshape}
\spnewtheorem{prop}[fact]{Proposition}{\bfseries}{\itshape}
\spnewtheorem{defn}[fact]{Definition}{\bfseries}{}
\spnewtheorem{exam}[fact]{Example}{\itshape}{}
\spnewtheorem{fct}{Fact}{\bfseries}{\itshape}

\spnewtheorem{exa}{Example}{\itshape}{}

\spdefaulttheorem{defin}{Definition}{\bfseries}{}

\spnewtheorem*{remnonumber}{Remark}{\itshape}{}

\begin{document}

\title*{Multiplicity-free theorems of the 
Restrictions of Unitary Highest Weight Modules
        with respect to Reductive Symmetric Pairs}
\titlerunning{Multiplicity-free restrictions to symmetric pairs}
\author{Toshiyuki KOBAYASHI}%
\institute{RIMS, %
Kyoto University, Kyoto 606-8502, Japan
\texttt{toshi@kurims.kyoto-u.ac.jp}}
\maketitle

\begin{abstract}

The complex analytic methods have found
 a wide range of applications
 in the study of 
 multiplicity-free representations.
This article discusses, in particular,
 its applications
 to the question of
 restricting highest weight modules
 with respect to reductive symmetric pairs.
We present a number of
 multiplicity-free branching theorems
 that include the multiplicity-free
 property of some of known results
 such as the Clebsh--Gordan--Pieri formula for
  tensor products,
  the Plancherel theorem for
  Hermitian symmetric spaces (also for line bundle cases),
  the Hua--Kostant--Schmid $K$-type formula,
  and
  the canonical representations
   in the sense of
    Vershik--Gelfand--Graev.
Our method works in a uniform manner
 for both finite and infinite dimensional
 cases,
 for both discrete and continuous spectra,
 and for both classical and exceptional cases.


\noindent
\subclassname \ 
Primary
         22E46, %
Secondary 
   32A37, 05E15, 20G05,
          53C35. %

\end{abstract}

\noindent
\textbf{Contents}
\setcounter{minitocdepth}{1}
\dominitoc

\numberwithin{equation}{subsection}
\numberwithin{table}{subsection}

\section{Introduction and statement of main results}
\label{sec:1}

The purpose of this article is to 
 give a quite detailed account of the theory
 of multiplicity-free representations
 based on a  non-standard method
 (\textit{visible actions} on complex manifolds)
 through its application to 
 branching problems.
More precisely, 
 we address the question of restricting
 irreducible highest weight
  representations $\pi$ of reductive Lie groups $G$ 
with respect to symmetric pairs
$(G,H)$.
Then, our main goal is to
 give a simple and sufficient condition on
the triple
$(G,H,\pi)$ such that the restriction $\pi|_H$ 
is 
multiplicity-free.
We shall see that our method works
 in a uniform way for
 both infinite and finite dimensional
  representations,
  for both classical and exceptional cases,
   and for both continuous and discrete spectra.

This article is an outgrowth of the  manuscript
\cite{xkmf} which I did not publish,
but which has been circulated as a preprint. 
From then onwards, we have extended the theory,
 in particular, 
to the following three directions:
\newline
1)\enspace
the generalization of our main machinery (Theorem~\ref{thm:2.2}) to the
vector bundle case (\cite{mfbdle}),
\newline
2)\enspace
the theory of \lq{visible actions}\rq\ 
 on complex manifolds
(\cite{visiblesymm, xkgencar, xksovisible}),
\newline
3)\enspace
`multiplicity-free geometry' for coadjoint orbits
(\cite{xknasrin}).
\par
We refer the reader to our paper
 \cite{RIMS} for a precise statement of 
 the general results and an exposition
 of the related topics 
 that have recently developed.

In this article,
 we confine ourselves to the line bundle case.
On the one hand,
this is sufficiently general to produce 
 many interesting consequences,
  some of which are new
  and some others may be regarded as
  prototypes of various 
 multiplicity-free branching theorems 
(e.g.\ 
\cite{xsaid, xvdh, xkleiden, xko, xkratten,
 xnere, xokada, xstemgl, xyawa, xzhang}).
On the other hand, 
 the line bundle case is
 sufficiently simple,
 so that we can illustrate the essence of
 our main ideas without going into
 technical details. 
Thus, keeping the spirit of \cite{xkmf},
we have included here
the 
proof of our method
 (Theorem~\ref{thm:2.2}),
  its applications to multiplicity-free theorems
(Theorems~\ref{thm:A}--\ref{thm:F}),
and the explicit formulae
(Theorems~\ref{thm:gHKS}, \ref{thm:tensordeco}, and \ref{thm:upqupq}),
 except that we referred to another paper \cite{visiblesymm} for
 the proof of some algebraic lemmas
 on the triple of involutions of Lie algebras
 (Lemmas~\ref{lem:5.1}
and \ref{lem:7.5}).

\subsection{Definition of multiplicity-free representations}
\label{subsec:1.2}
Let us begin by recalling
 the concept of the multiplicity-free decomposition
 of a unitary representation.

Suppose $H$ is a Lie group of type I in the sense of von Neumann
algebras. 
Any reductive Lie group is of type I as well as any algebraic group. 
We denote by $\widehat{H}$ the unitary dual of $H$,
 that is, the set of equivalence classes of irreducible unitary representations
 of $H$.
The unitary dual $\widehat{H}$ is endowed with the Fell topology.

Suppose that $(\pi, \mathcal{H})$ is a unitary representation of $H$
 defined on a (second countable) Hilbert space $\mathcal{H}$.
By a theorem of Mautner, 
 $\pi$ is decomposed uniquely
 into irreducible unitary representations of $H$
 in terms of the direct integral of Hilbert spaces:
\begin{equation}
\label{eqn:1.2.1}
   \pi \simeq \int_{\widehat{H}} m_\pi(\mu) \mu \, d \sigma (\mu) \, ,
\end{equation}
 where $d \sigma (\mu)$ is a Borel measure on $\widehat{H}$,
 and the multiplicity function $m_\pi : \widehat{H} \to \mathbb{N}
 \cup \{\infty\}$
 is uniquely defined almost everywhere with respect to the measure $d \sigma$.

Let $\End(\mathcal{H})$ be the ring of continuous operators on $\mathcal{H}$,
 and
 $\End_H(\mathcal{H})$
the subring of $H$-intertwining operators, that is,
 the commutant of $\set{\pi(g)}{g \in H}$ in $\End(\mathcal{H})$.
\begin{defin}
\label{def:1.2}
We say that the unitary representation $(\pi, \mathcal{H})$ is
\emph{multiplicity-free} if the ring
$\End_H(\mathcal{H})$ is commutative.
\end{defin}

It is not difficult to see that this definition
 is equivalent to the following property:
$$
\text{$m_\pi (\mu) \le 1$ for almost all $\mu \in \widehat{H}$
 with respect to the measure $d \sigma (\mu)$ }
$$
by Schur's lemma for unitary representations.
In particular, it implies  that
any irreducible unitary representation $\mu$ of $H$ occurs at most
once as a subrepresentation of $\pi$.

\subsection{Multiplicities for inductions and restrictions}
\label{subsec:1.3}
With regard to
 the question of finding irreducible decompositions of unitary
representations, 
there are two fundamental settings:
one is the induced representation from  
 smaller groups
 (e.g. harmonic analysis on homogeneous spaces), 
and the other is the restriction from  larger groups
 (e.g. tensor product representations).

To be more rigorous, 
 suppose $G$ is a Lie group,
 and $H$ is a closed subgroup of $G$.  
The $G$-irreducible decomposition
 of the induced representation
  $L^2$-$\operatorname{Ind}_H^G \tau$
 ($\tau \in \widehat H$) is called
   the {\it Plancherel formula}, 
 while the $H$-irreducible decomposition 
 of the restriction $\pi |_H$ ($\pi \in \widehat G$) 
 is referred to as the {\it branching law}.  

This subsection examines multiplicities
 in the irreducible decomposition
  of the induction and the restriction
  for reductive symmetric pairs $(G,H)$
(see Subsection~\ref{subsec:3.1} for definition).

Let us start with the induced representation.
 Van den Ban \cite{xban} proved that 
 the multiplicity in the Plancherel formula for 
 $L^2$-$\operatorname{Ind}_H^G \tau$ is finite
as far as $\dim \tau < \infty$.
In particular, this is the case
 if $\tau$ is the trivial representation 
$\mathbf{1}$.
Over the past several decades, 
the induced representation 
 $L^2$-$\operatorname{Ind}_H^G \mathbf{1}$  %
 has developed its own identity
 (harmonic analysis on  reductive symmetric spaces $G/H$)
 as a rich and meaningful part of mathematics.

In contrast,
 the multiplicities of
 the branching law of the restriction
 $\pi|_H$
 ($\pi \in \widehat G$) are usually infinite.
For instance, 
 we saw in \cite{xkbdlejp} that
 this is the case if
 $(G, H)=(GL(p+q,{\mathbb {R}}), GL(p,{\mathbb{R}}) 
  \times GL(q,{\mathbb{R}}))$
where $\min(p,q) \ge 2$, 
for any tempered representation $\pi$ of $G$.
In this article,
we illuminate by Example~\ref{ex:finite infinite}
this wild behavior. 

 In light of such a wild phenomenon of branching laws 
 for reductive symmetric pairs $(G,H)$ with $H$ non-compact, 
we proposed in \cite{xkdecomp, xkdecoalg} 
to seek for
 a `nice' class of the triple
 $(G, H, \pi)$ in which a systematic
 study of the restriction $\pi|_H$
 could be launched.

Finiteness of multiplicities is
 a natural requirement for this program.
By also imposing discrete
decomposability on the restriction 
  $\pi|_H$,
 we established the general theory for
 \textit{admissible restriction}
  in \cite{xkdecomp, xkdecoalg, xkdecoass}
and found
 that there exist fairly rich triples $(G, H, \pi)$
 that enjoy this nice property.  
It is noteworthy that 
 new interesting directions
  of research
   in the framework of admissible restrictions
    have been recently developed by 
M. Duflo, D. Gross, J.-S.\ Huang, J.-S.\ Li, S.-T.\ Lee, H.-Y.\ Loke,
T. Oda, P. Pand\v{z}i\'{c}, G. Savin, B. Speh, 
J. Vargas, D. Vogan, and N. Wallach 
(see \cite{xkbeijing, xkronsdecomp} and references therein).

Multiplicity-freeness is another
 ideal situation,
in which we may expect
an especially simple and detailed study
 of the branching law of $\pi|_H$.
Thus, we aim for principles
 that lead us to  abundant
 family of multiplicity-free
 cases.
Among them, a well-known one
 is the dual pair correspondence,
which has given fruitful examples
 in infinite dimensional
theory in the following setting:
\begin{list}{}{\setlength{\leftmargin}{2\parindent}}
\item[a)]
 $G$ is the metaplectic group,
and $\pi$ is the Weil representation.
\item[b)]
$H = H_1 \cdot H_2$ forms a dual pair,
that is,
$H_1$ is the commutant of $H_2$ in $G$, and vice versa.
\end{list}
This paper uses a new principle
 that generates multiplicity-free representations.  
The general theory discussed in Section~\ref{sec:2}
brings us to uniformly bounded multiplicity theorems 
(Theorems~\ref{thm:B} and \ref{thm:D})
and 
multiplicity-free theorems (Theorems~\ref{thm:A}, \ref{thm:C},
\ref{thm:E} and \ref{thm:F})
in the following setting:
\begin{list}{}{\setlength{\leftmargin}{2\parindent}}
\item[a)]
 $\pi$
 is a unitary highest weight representation of $G$
(see Subsection~\ref{subsec:1.4}), %
\item[b)]
  $(G,H)$ is a symmetric pair 
(see Subsection~\ref{subsec:1.5}).  %
\end{list}

We note that we allow the case where continuous spectra occur
in the branching law,
and consequently,
irreducible summands are not always highest weight representations.

We remark that our bounded multiplicity theorems for the restriction 
$\pi|_H$ ($\pi$: highest weight module) may be regarded as the 
counterpart of the bounded multiplicity theorem for the induction 
$L^2$-$\operatorname{Ind}_H^G \tau$
($\tau$: finite dimensional representation) due to van den Ban.

\subsection{Unitary highest weight modules}
\label{subsec:1.4}
Let us recall the basic notion of highest weight modules.

Let $G$ be a non-compact simple Lie group,
 $\theta$ a Cartan involution of $G$,
 and $K := \set{g \in G}{\theta g = g}$.
We write
 $\mathfrak{g} = \mathfrak{k} + \mathfrak{p}$ 
 for the Cartan decomposition
 of the Lie algebra $\mathfrak{g}$ of $G$,
 corresponding to the Cartan involution $\theta$. 

We assume that $G$ is of {\it Hermitian type},
 that is,
 the Riemannian symmetric space $G/K$ carries
 the structure of a Hermitian symmetric space,
 or equivalently, 
 the center $\mathfrak{c}(\mathfrak{k})$ of $\mathfrak{k}$ is non-trivial. 
The classification of simple Lie algebras ${\mathfrak {g}}$
 of Hermitian type
 is given as follows:
$$
{\mathfrak {su}}(p,q)\,,\ 
{\mathfrak {sp}}(n,{\mathbb{R}})\,,\ 
{\mathfrak {so}}(m,2) \ (m \ne 2)\,,\ 
{\mathfrak {e}}_{6(-14)}\,,\ 
{\mathfrak {e}}_{7(-25)} \, .  
$$
Such a Lie algebras $\mathfrak{g}$ satisfies the rank condition:
\begin{equation}\label{eqn:rankGK}
\operatorname{rank} G = \operatorname{rank} K \, ,
\end{equation}
or equivalently, a Cartan subalgebra of $\mathfrak{k}$ becomes a
Cartan subalgebra of $\mathfrak{g}$. 
By a theorem of Harish-Chandra,
the rank condition \eqref{eqn:rankGK} is equivalent to the existence of
(relative) discrete series representations of $G$. 
Here, an irreducible unitary representation $(\pi,\mathcal{H})$ is
called a \textit{(relative) discrete series representation} 
of $G$ if the matrix coefficient
$g \mapsto (\pi(g)u,v)$ is square integrable on $G$
(modulo its center) for any $u,v \in \mathcal{H}$. 

If $\mathfrak{g}$ is a simple Lie algebra of Hermitian type,
then there exists a characteristic element
 $Z \in \mathfrak{c}(\mathfrak{k})$
 such that 
\begin{equation}
\label{eqn:gkppz}
      \mathfrak{g}_\mathbb{C} 
    := \mathfrak{g} \otimes \mathbb{C}
     = \mathfrak{k}_\mathbb{C} \oplus \mathfrak{p}_+ \oplus \mathfrak{p}_-
\end{equation}
 is the eigenspace decomposition of $\operatorname{ad}(Z)$
 with eigenvalues $0$, $\sqrt{-1}$ and $-\sqrt{-1}$,
 respectively.
We note that
 $\dim \mathfrak{c}(\mathfrak{k}) = 1$
 if $\mathfrak{g}$ is 
 a simple Lie algebra of Hermitian type,
and therefore
$\mathfrak{c}(\mathfrak{k}) = \mathbb{R} Z$.

Suppose $V$ is an irreducible $({\mathfrak {g}}_{\mathbb{C}}, K)$-module.  
We set
\begin{equation}
\label{eqn:Hpk}
   V^{\mathfrak{p}_+}
   := \set{v \in V}{Y v = 0 
   \ \text{ for any } Y \in \mathfrak{p}_+}\, .
\end{equation}
Since $K$ normalizes $\mathfrak{p}_+$,
$V^{\mathfrak{p}_+}$ is a $K$-submodule.
Further,
$V^{\mathfrak{p}_+}$
 is either zero or an irreducible finite
 dimensional representation of $K$.
We say
$V$ is 
 a \textit{highest weight module}
 if $V^{\mathfrak{p}_+} \neq \{0\}$.  

\begin{defin}  %
\label{def:1.4}
Suppose $\pi$ is an irreducible unitary representation
of $G$ on a Hilbert space $\mathcal{H}$.
We set
$\mathcal{H}_K:=\{v\in\mathcal{H}:
 \dim_{\mathbb{C}} \mathbb{C}\text{-span}
 \{\pi(k)v:k\in K\} < \infty \}
$.
Then, $\mathcal{H}_K$ is a dense subspace of $\mathcal{H}$,
on which the differential action $d\pi$ of the Lie algebra
$\mathfrak{g}$ (and consequently that of its complexified Lie algebra
$\mathfrak{g}_{\mathbb{C}}$) and the action of the compact subgroup $K$
is well-defined.
We say $\mathcal{H}_K$ is 
\textit{the underlying 
$(\mathfrak{g}_{\mathbb{C}},K)$-module of $(\pi,\mathcal{H})$}. 
We say $(\pi,\mathcal{H})$ is a 
\textit{unitary highest weight representation}
 of $G$ if $\mathcal{H}_K^{\mathfrak{p}_+} \ne \{0\}$.
Then, $\pi$ is 
\emph{of scalar type}
 (or \emph{of scalar minimal $K$-type})
 if $\mathcal{H}^{\mathfrak{p}_+}_K$ is one dimensional;
 $\pi$ is a 
\textit{(relative) 
holomorphic discrete series representation}
 for $G$
if the matrix coefficient
$g \mapsto (\pi(g)u,v)$
is square integrable on $G$ modulo its center for any
$u, v \in \mathcal{H}$.
Lowest weight modules and anti-holomorphic discrete series representations
 are defined similarly 
 with $\mathfrak{p}_+$ replaced by $\mathfrak{p}_-$.
\end{defin}

This definition also applies to $G$ which is not simple 
(see Subsection~\ref{subsec:8.1}).

The classification of irreducible unitary highest weight representations
 was accomplished
 by Enright--Howe--Wallach \cite{xhew}
 and H. Jakobsen \cite{xjak} independently;
 see also \cite{xej}.
There always exist infinitely many 
(relative) holomorphic discrete series representations
 of scalar type for any non-compact simple Lie group of Hermitian type.

\subsection{Involutions on Hermitian symmetric spaces}
\label{subsec:1.5}  %
Suppose $G$ is a non-compact simple Lie group of Hermitian type.
Let $\tau$ be an involutive automorphism of $G$
commuting with the Cartan involution $\theta$.
We use the same letter $\tau$ to denote
 its differential.
Then $\tau$ stabilizes $\mathfrak{k}$ and also $\mathfrak{c}(\mathfrak{k})$.
Because $\tau^2 = \operatorname{id}$
 and $\mathfrak{c}(\mathfrak{k}) = \mathbb{R} Z$,
 we have the following two possibilities:
\begin{align}
    \tau Z &= Z \, ,
\label{eqn:1.5.1}
\\
   \tau Z &= -Z \, .
\label{eqn:1.5.2}
\end{align}
Geometric meanings of these conditions
 become clear in the context of the embedding 
 $G^{\tau}/ K^{\tau} \hookrightarrow G/K$, 
where $G^\tau := \set{g \in G}{\tau g = g}$
 and $K^\tau := G^\tau \cap K$
 (see \cite{xfo, xjafbams, xjafjdg, xkobanaga}). 
 The condition \eqref{eqn:1.5.1} %
 implies: 
\newline\indent{1-a)}\enspace
 $\tau$ acts {\bf holomorphically} 
 on the Hermitian symmetric space $G/K$,
\newline\indent{1-b)}\enspace
 $G^\tau/K^\tau \hookrightarrow G/K$ defines a complex submanifold,
\newline
whereas the condition
\eqref{eqn:1.5.2} %
implies:
\newline\indent{2-a)}\enspace
 $\tau$ acts {\bf anti-holomorphically} on  $G/K$,
\newline\indent{2-b)}\enspace
 $G^\tau/K^\tau \hookrightarrow G/K$ defines a totally real submanifold.

\begin{defin}
\label{def:holo-anti}
We say the involutive automorphism $\tau$ is 
\textit{of holomorphic type} if \eqref{eqn:1.5.1} is satisfied,
and is of 
\textit{anti-holomorphic type}
if \eqref{eqn:1.5.2} is satisfied.
The same terminology will be applied also to the symmetric pair 
$(G,H)$
(or its Lie algebras $(\mathfrak{g}, \mathfrak{h})$)
corresponding to the involution $\tau$.
\end{defin}

Here, we recall that $(G,H)$ is called a 
\textit{symmetric pair} corresponding
to $\tau$ if $H$ is an open subgroup of $G^\tau$
(see Subsections~\ref{subsec:3.1} and \ref{subsec:3.2.ex}).
We note that the Lie algebra $\mathfrak{h}$ of $H$ is equal to 
$\mathfrak{g}^\tau
 := \set{X \in \mathfrak{g}}{\tau X = X}$.
The classification of  symmetric pairs
$(\mathfrak{g}, \mathfrak{g}^\tau)$
 for simple Lie algebras $\mathfrak{g}$
 was accomplished by M. Berger \cite{xber}.
 The classification of symmetric pairs
 $({\mathfrak {g}}, {\mathfrak {g}}^{\tau})$ of holomorphic type
(respectively, of anti-holomorphic type)
is regarded as a subset of
Berger's list, and will be presented
 in Table~\ref{tbl:3.3.1} (respectively, Table~\ref{tbl:3.3.2}). 

\subsection{Multiplicity-free restrictions --- infinite dimensional
case}
\label{subsec:1.6}
We are ready to state our main results.
Let $G$ be a non-compact simple Lie group of Hermitian type,
 and 
 $(G,H)$ a symmetric pair.

\begin{theopargself}
\begin{thmalph}[\textmd{(multiplicity-free restriction)}]
\label{thm:A}
If $\pi$ is an irreducible unitary highest weight representation 
 of scalar type of $G$,
 then the restriction $\pi|_H$ is multiplicity-free.
\end{thmalph}
\end{theopargself}

\par
The branching law of the restriction $\pi |_H$ may and may not contain
discrete spectra in Theorem~\ref{thm:A}.
If $(G,H)$ is of holomorphic type then the 
restriction $\pi |_H$ is discretely decomposable
(i.e.\ there is no continuous spectrum in the branching law);
see Fact~\ref{fact:3.4.1}.
Besides,
the following theorem asserts that
the multiplicities are still uniformly bounded even if we drop
 the assumption
 that $\pi$ is of scalar type.  

\begin{theopargself}
\begin{thmalph}[\textmd{(uniformly bounded multiplicities)}]
\label{thm:B}
We assume that the symmetric pair $(G,H)$ is of holomorphic type.
Let $\pi$ be an irreducible unitary highest weight representation of $G$.
\itm{1)}
The restriction $\pi|_H$ splits into a discrete Hilbert sum of
  irreducible unitary representations of $H$:
$$
  \pi|_H \simeq \sideset{}{^\oplus}\sum_{\mu \in \widehat{H}} 
m_\pi(\mu) \mu \, ,
$$ 
and the multiplicities are uniformly bounded:
$$
   C(\pi) := \sup_{\mu \in \widehat{H}}
    m_\pi(\mu) < \infty \, .
$$
\itm{2)}
 $C(\pi) = 1$ if $\pi$ is of scalar type.
\end{thmalph}
\end{theopargself}

The second statement is a direct consequence of Theorems~\ref{thm:A}
 and \ref{thm:B}~(1).
As we shall see in Section \ref{sec:6}, 
 such uniform boundedness theorem does not hold
in general if $\pi$ is not a highest weight representation
(see Examples~\ref{exam:6.3} and \ref{ex:finite infinite}).  

Here are multiplicity-free theorems for the decomposition of tensor products,
 which are parallel to Theorems~\ref{thm:A} and \ref{thm:B}:

\begin{theopargself}
\begin{thmalph}[\textmd{(multiplicity-free tensor product)}]
\label{thm:C}
Let $\pi_1$ and $\pi_2$ be irreducible
 unitary highest (or lowest)
 weight representations of scalar type.
Then the tensor product $\pi_1 \widehat\otimes \pi_2$
 is multiplicity-free as a representation of $G$.
\end{thmalph}
\end{theopargself}

Here, $\pi_1 \widehat{\otimes} \pi_2$ stands for the tensor product
representation of two unitary representations
$(\pi_1, \mathcal{H}_1)$ and $(\pi_2, \mathcal{H}_2)$ realized on the
completion $\mathcal{H}_1 \widehat{\otimes} \mathcal{H}_2$ of the
pre-Hilbert space $\mathcal{H}_1 \otimes \mathcal{H}_2$.
(We do not need to take the completion if at least one of
 $\mathcal{H}_1$ or
$\mathcal{H}_2$ is finite dimensional.)
Theorem~\ref{thm:C} asserts that multiplicities in the direct integral
of the irreducible decomposition are not greater than
one in
 both discrete and continuous spectra.
We note that continuous spectra appear in the irreducible
decomposition of the tensor product representation 
$\pi_1 \widehat{\otimes} \pi_2$  only if
\begin{align*}
&\begin{cases}
  \text{$\pi_1$ is a highest weight representation, and} \\
  \text{$\pi_2$ is a lowest weight representation,}
 \end{cases}
\end{align*}
or in reverse order.

 If $\pi_1$ and $\pi_2$ are simultaneously highest 
 weight representations
(or simultaneously lowest weight representations), 
 then the tensor product $\pi_1 \widehat\otimes \pi_2$ decomposes discretely.
Dropping the assumption of \lq\lq scalar type\rq\rq,
 we have still a uniform estimate of multiplicities:

\begin{theopargself}
\begin{thmalph}[\textmd{(uniformly bounded multiplicities)}]
\label{thm:D}
Let $\pi_1$ and $\pi_2$ be two irreducible unitary highest weight 
representations of
 $G$.
\itm{1)}
The tensor product $\pi_1 \widehat\otimes \pi_2$
 splits into a discrete Hilbert sum
 of irreducible unitary representations of $G$:
$$
  \pi_1 \widehat\otimes \pi_2 \simeq 
  \sideset{}{^\oplus}\sum_{\mu \in \widehat{G}} 
  m_{\pi_1, \pi_2}(\mu) \mu \, ,
$$ 
and the multiplicities $m_{\pi_1,\pi_2}(\mu)$
are uniformly bounded:
$$
   C(\pi_1,\pi_2) := \sup_{\mu \in \widehat{G}}
     m_{\pi_1, \pi_2}(\mu) < \infty \, .
$$
\itm{2)}
$C(\pi_1, \pi_2) = 1$ if both $\pi_1$ and $\pi_2$ are
 of scalar type.
\end{thmalph}
\end{theopargself}

\begin{remark}
\label{rem:D}
For classical groups,
we can relate the constants $C(\pi)$ and $C(\pi_1,\pi_2)$ to the 
\textit{stable constants}
 of branching coefficients
of finite dimensional representations in the sense of
 F. Sato \cite{xsatof} by using 
the see-saw dual pair correspondence due to R. Howe \cite{xhoweseesaw}.
\end{remark}

Our machinery that gives the above multiplicity-free theorems is
built on complex geometry,
and we shall explicate the general theory for the line bundle case in
Section~\ref{sec:2}.
The key idea is to transfer properties on representations 
(e.g.\ unitarity, multiplicity-freeness)
into the corresponding properties of reproducing kernels,
 which we analyze by geometric methods.

\subsection{Multiplicity-free restrictions --- finite dimensional case}
\label{subsec:1.7}
Our method yields multiplicity-free theorems not only for infinite
dimensional representations but also for finite dimensional
representations. 

This subsection presents multiplicity-free theorems that are regarded as
`finite dimensional version' of Theorems~\ref{thm:A} and \ref{thm:C}.
They give a 
 unified explanation of the multiplicity-free property of
previously known branching formulae obtained
 by combinatorial methods such as the
Littlewood--Richardson rule, 
Koike--Terada's Young diagrammatic methods,
Littelmann's path method,
minor summation formulae, etc.\ 
(see
\cite{xhowemultone, xkoiterjalg, 
xmacd, xokada, xproctor, xstem}
and references therein).
They also contain some  `new'
cases, 
 for which there are,
 to the best of our knowledge,
 no explicit branching formulae in the literature.

To state the theorems, 
let $\mathfrak{g}_{\mathbb{C}}$ be a complex simple Lie algebra,
and $\mathfrak{j}$ a Cartan subalgebra.
We fix a positive root system
$\Delta^+ (\mathfrak{g}_{\mathbb{C}}, \mathfrak{j})$,
and write $\alpha_1,\dots,\alpha_n$ for the simple roots.
Let
$\omega_1,\dots,\omega_n$ be the corresponding fundamental weights. 
We denote by 
$\pi_\lambda \equiv \F{\mathfrak{g}_{\mathbb{C}}}{\lambda}$
the irreducible finite dimensional representation of
$\mathfrak{g}_{\mathbb{C}}$ with highest weight $\lambda$.

We say $\pi_\lambda$ is of \textbf{pan type}
if $\lambda$ is a scalar multiple of some $\omega_i$
such that the nilradical of 
the maximal parabolic subalgebra corresponding to $\alpha_i$
is abelian 
(see Lemma~\ref{lem:7.3.1} for equivalent definitions).

\begin{theopargself}
\begin{thmalph}[\textmd{(multiplicity-free restriction --- finite 
dimensional case)}] 
\label{thm:E}
Let $\pi$ be an arbitrary irreducible 
finite dimensional representation of
$\mathfrak{g}_{\mathbb{C}}$ 
of pan type,
and $(\mathfrak{g}_{\mathbb{C}}, \mathfrak{h}_{\mathbb{C}})$
be any symmetric pair.
Then, the restriction $\pi|_{\mathfrak{h}_{\mathbb{C}}}$
is multiplicity-free.
\end{thmalph}
\end{theopargself}

\begin{theopargself}
\begin{thmalph}[\textmd{(multiplicity-free tensor product --- finite
dimensional case)}]
\label{thm:F}
The tensor product $\pi_1 \otimes \pi_2$ 
of any two irreducible finite dimensional
representations $\pi_1$ and $\pi_2$
of pan type is multiplicity-free.
\end{thmalph}
\end{theopargself}

Theorems~E and F are the counterpart to Theorems~\ref{thm:A} and \ref{thm:C}
for finite dimensional representations.
The main machinery of the proof is again Theorem~\ref{thm:2.2}.

Alternatively, 
one could verify Theorems~\ref{thm:E} and \ref{thm:F} by a classical
technique: 
finding an open orbit of a Borel subgroup.
For example,
Littelmann \cite{xlittelmann}
and Panyushev independently classified the pair of maximal parabolic
subalgebras $(\mathfrak{p}_1, \mathfrak{p}_2)$ such that the diagonal
action of a Borel subgroup $B$ 
of a complex simple Lie group $G_{\mathbb{C}}$
on 
$G_{\mathbb{C}} / P_1 \times G_{\mathbb{C}} / P_2$
has an open orbit.
Here, $P_1, P_2$ are the corresponding maximal parabolic subgroups of
$G_{\mathbb{C}}$. 
This gives another proof of Theorem~\ref{thm:F}.

The advantage of our method is that it enables us to understand
(or even to discover) the multiplicity-free property simultaneously, 
 for both infinite and finite dimensional
representations,
 for both continuous and discrete spectra,
and for both classical and exceptional cases
by the single principle.
This is because our main machinery (Theorem~\ref{thm:2.2}) uses only a
\textit{local} geometric assumption
(see Remark~\ref{rem:2.3.2} (2)).
Thus, we can verify it at the same time for both compact
and non-compact
 complex manifolds,
and in turn get
finite and infinite dimensional results, respectively.

Once we tell a priori that a representation is multiplicity-free, 
we may be tempted to find explicitly its irreducible decomposition.
Recently, S. Okada \cite{xokada}
 found explicit branching laws for some 
classical cases that arise in Theorems~\ref{thm:E} and \ref{thm:F} by
using minor summation formulae,
and H. Alikawa \cite{xalikawa} for
$(\mathfrak{g}, \mathfrak{h}) = (\mathfrak{e}_6, \mathfrak{f}_4)$
corresponding to Theorem~\ref{thm:E}.
We note that
the concept of  pan type representations includes
\textit{rectangular-shaped} representations of classical groups
(see \cite{xkratten, xokada}).

There are also some few cases where $\pi_1 \otimes \pi_2$
is multiplicity-free even though neither $\pi_1$ nor
$\pi_2$ is of pan type.
See the recent papers
 \cite{xkleiden} or \cite{xstemgl} for the complete list of such
pairs $(\pi_1, \pi_2)$ for 
$\mathfrak{g}_{\mathbb{C}} = \mathfrak{gl}(n, \mathbb{C})$.
The method in \cite{xkleiden} to find all such pairs is geometric and based on
the
 `vector bundle version' of Theorem~\ref{thm:2.2} 
proved in \cite{mfbdle},
whereas the method in \cite{xstemgl} is combinatorial and  based on
case-by-case argument. 

We refer the reader to our papers \cite{visiblesymm, xkgencar, xksovisible}
for some further results relevant to Theorems~\ref{thm:E} and
\ref{thm:F} along the same line of argument here.

\subsection{$SL_2$ examples}
\label{subsec:1.8}
We illustrate the above theorems by $SL_2$ examples.
\begin{exa}
\label{exam:1.8}
1)\enspace
We denote by $\pi_n$
 the holomorphic discrete series representation 
 of $G=SL(2,\mathbb{R})$
 with minimal $K$-type $\chi_n$ ($n \ge 2$),
 where we write $\chi_n$ 
 for the character of $K=SO(2)$
parametrized by $n \in \mathbb{Z}$.
Likewise $\pi_{-n}$ denotes the anti-holomorphic discrete series representation
 of $SL(2, \mathbb{R})$ with minimal $K$-type $\chi_{-n}$ ($n \ge 2$).
We note that any
 holomorphic discrete series of $SL(2,\mathbb{R})$ is of scalar
type.

We write $\pi^\varepsilon_{\sqrt{-1} \nu}$
 ($\varepsilon = \pm 1, \nu \in \mathbb{R})$ for the unitary
 principal series representations of $SL(2, \mathbb{R})$.
We have a unitary equivalence %
 $\pi^{\varepsilon}_{\sqrt{-1} \nu} \simeq 
 \pi^{\varepsilon}_{- \sqrt{-1} \nu}$.
We write $\chi_\zeta$ 
 for the unitary character of $SO_0(1,1) \simeq \mathbb{R}$
parametrized by $\zeta \in \mathbb{R}$.

Let $m \ge n \ge 2$.
Then,
 the following branching formulae %
 hold.
All of them are multiplicity-free, 
as is `predicted' by Theorems~\ref{thm:A} and
\ref{thm:C}:
\begin{subequations}
  \renewcommand{\theequation}{\theparentequation) (\alph{equation}}%
\begin{align}
   \pi_n|_{SO_0(1,1)}
   &\simeq \int_{-\infty}^\infty \chi_{\zeta} \, d \zeta \, ,  
\label{eqn:1.8.1a}
\\
   \pi_n|_{SO(2)}
   &\simeq \sideset{}{^\oplus}\sum_{k\in \mathbb{N}} \chi_{n + 2 k} \, ,  
\label{eqn:1.8.1b}
\\
   \pi_m \widehat\otimes \pi_{-n}
   &\simeq
    \int_{0}^\infty \pi^{(-1)^{m-n}}_{\sqrt{-1} \nu}  d \nu \quad
    \oplus \sum \Sb k \in \mathbb{N} \\ 0 \le 2k \le m-n-2 \endSb 
     \pi_{m-n-2k} \, ,
\label{eqn:1.8.1c}
\\
   \pi_m \widehat\otimes \pi_n
   &\simeq \sideset{}{^\oplus}\sum_{k\in \mathbb{N}} \pi_{m+n + 2 k} \, .
\label{eqn:1.8.1d}
\end{align}
\end{subequations}
The key assumption of our main machinery (Theorem~\ref{thm:2.2}) that
leads us to Theorems~\ref{thm:A} and \ref{thm:C} is illustrated
by the following geometric results in this $SL_2$ case:
\newline{i)}\enspace
Given any element $z$ in the Poincar\'e disk $D$,
 there exists $\varphi \in \mathbb{R}$ such that
 $e^{\sqrt{-1}\varphi} z = \overline{z}$.
In fact, one can take $\varphi = -2 \arg z$.
This is the geometry that explains the multiplicity-free property of
 \eqref{eqn:1.8.1b}. %
\newline{ii)}\enspace
Given any two elements $z, w \in D$,
 there exists a linear fractional transform $T$ on $D$ such that 
 $T (z) =  \overline{z}$ and $T (w) =  \overline{w}$. 
 This is the geometry for  \eqref{eqn:1.8.1d}. %

These are examples of the geometric view point that we pursued in
\cite{visiblesymm} for symmetric pairs.
\itm{2)}\enspace
Here is a \lq\lq{finite dimensional version}\rq\rq\
 of the above example.
Let $\pi_n$ be the irreducible $n+1$-dimensional
 representation of $SU(2)$.  
Then we have the following branching formulae:
For $m, n \in \mathbb{N}$, 
\addtocounter{equation}{-1}
\begin{subequations}
  \renewcommand{\theequation}{\theparentequation) (\alph{equation}}%
\abovedisplayskip0pt
\begin{align}
\addtocounter{equation}{4}
\hphantom{\pi_n|_{SO_0(1,1)}}
\llap{$\pi_n|_{SO(2)}$}
       &\simeq 
\rlap{$\chi_{n} \oplus \chi_{n-2} \oplus \cdots \oplus \chi_{-n} \, ,$}
{\hphantom{\int_{0}^\infty \pi^{(-1)^{m-n}}_{\sqrt{-1} \nu}  d \nu \quad
    \oplus \sum \Sb k \in \mathbb{N} \\ 0 \le 2k \le m-n-2 \endSb 
     \pi_{m-n-2k} \, ,}}
{\vphantom{\bigg|}}
\label{eqn:1.8.1e}
\\
         \pi_m \otimes \pi_n
       &\simeq 
       \pi_{n+m} \oplus \pi_{n+m-2} \oplus \cdots \oplus \pi_{|n-m|}
\, .
\label{eqn:1.8.1f}
\end{align}
\end{subequations}
The formula \eqref{eqn:1.8.1e} corresponds to the character formula,
whereas \eqref{eqn:1.8.1f} is known as the Clebsch--Gordan formula.
The multiplicity-free property of these formulae is
the simplest example of Theorems~\ref{thm:E} and \ref{thm:F}.

\subsection{Analysis on multiplicity-free representations}
\label{subsec:analysismf}

Multiplicity-free property arouses our interest in developing
beautiful analysis on such representations,
as we discussed in Subsection~\ref{subsec:1.7} for finite dimensional
cases. 
This subsection picks up some recent topics about detailed analysis on
multiplicity-free representations for infinite dimensional cases.

Let $G$ be a connected,
simple non-compact Lie group of Hermitian type.
We begin with branching laws without continuous spectra,
and then discuss branching laws with continuous spectra.

\noindent
{1)}\enspace %
(Discretely decomposable case)\enspace
Let $(G,H)$ be a symmetric pair of holomorphic type.
Then, any unitary highest weight representation $\pi$ of $G$
decomposes discretely when restricted to $H$ (Fact~\ref{fact:3.4.1}).
\newline\indent{1-a)}
Suppose now that $\pi$ is a holomorphic discrete series
representation. 
L.-K. Hua \cite{xhua}, B. Kostant, W. Schmid \cite{xschmidherm} 
and K. Johnson \cite{xjohnson}
found an explicit formula of the restriction $\pi|_K$
 ($K$-type formula).
This turns out to be multiplicity-free.
Alternatively, the
 special case of Theorem~\ref{thm:B}~(2) by setting $H=K$
 gives a new proof of this multiplicity-free property.
\newline\indent{1-b)}
Furthermore,
 we consider a generalization of
 the Hua--Kostant--Schmid formula
from compact $H$ to noncompact $H$,
for which Theorem~\ref{thm:B} (2) still ensures that the
generalization will be multiplicity-free.
This generalized formula is stated in Theorem \ref{thm:gHKS},
which was originally given in
 \cite[Theorem~C]{xkmfjp}.
In Section~\ref{sec:8},
we give a full account of its proof.
W.~Bertram and J. Hilgert \cite{xbehi}
obtained some special cases independently,
and Ben Sa\"{\i}d \cite{xsaid} studied a quantative estimate of this multiplicity-free
$H$-type formula
(see also \cite{xyawa,xzhangtensor} for some \textit{singular} cases).
\newline\indent{1-c)}
The branching formulae of the restriction of 
\textit{singular} highest
weight representations $\pi$ are also interesting.
For instance,
the restriction of
 the Segal--Shale--Weil representation $\varpi$ of $Mp(n, \mathbb{R})$
 with respect to $U(p, n-p)$ (more precisely, its double covering)
 decomposes discretely into a multiplicity-free sum of 
the so called {\it{ladder representations}}
 of $U(p, n-p)$  (e.g.\ \cite[Introduction]{xkv}).
This multiplicity-free property 
 is a special case of Howe's correspondence
 because $(U(p, n-p), U(1))$ forms a dual pair in $Mp(n, \mathbb{R})$,
 and also is a special case of Theorem~\ref{thm:A}
 because $(\mathfrak{sp}(n, \mathbb{R}), {\mathfrak{u}}(p, n-p))$
 forms a symmetric pair.
Explicit branching laws for most of classical cases corresponding 
 to Theorems~\ref{thm:B}~(2) and \ref{thm:D}~(2) 
(see Theorems \ref{thm:gHKS}, \ref{thm:tensordeco}, \ref{thm:upqupq})
 can be obtained by using the 
 \lq\lq{see-saw
 dual pair}\rq\rq,
 which we hope to report in another paper.  
\itm{2)}
(Branching laws with continuous spectra)\enspace
Suppose $\pi_1$ is a highest weight module
 and $\pi_2$ is a lowest weight module,
 and both being of scalar type.
\newline\indent{2-a)}
If both $\pi_1$ and $\pi_2$ are discrete series
 representations in addition,
 then the tensor product $\pi_1 \widehat \otimes \pi_2$ is unitarily
 equivalent to the regular representation on $L^2(G/K, \chi)$,
 the Hilbert space of $L^2$-sections
 of the $G$-equivalent line bundle $G\times_K \mathbb{C}_\chi \to G/K$
 associated to some unitary character $\chi$ of $K$
 (R. Howe \cite{xhoweseesaw}, J. Repka \cite{xrep}).
In particular,
  Theorem~\ref{thm:C}
 gives a new proof of the multiplicity-free property
 of the Plancherel formula
 for $L^2(G/K,\chi)$.
Yet another proof
 of the multiplicity-free property
  of $L^2(G/K, \chi)$
  was given in
\cite[Theorem~21]{RIMS}
by still applying Theorem~\ref{thm:2.2} to the \textit{crown domain} 
(equivalently, the Akhiezer--Gindikin domain)
of the Riemannian symmetric space $G/K$.
The explicit decomposition of $L^2(G/K, \chi)$
 was found by J. Heckman \cite{xhecsch}
and N. Shimeno \cite{xshimeno}
that generalizes the work of Harish-Chandra, S. Helgason, 
  and S. Gindikin--F. Karpelevich
 for the trivial bundle case.  

In contrast to Riemannian symmetric spaces, 
it is known that  \lq\lq multiplicity-free property\rq\rq\
 in the Plancherel formula
 fails for  (non-Riemannian)  symmetric spaces 
$G/H$ in general 
 (see \cite{xbsann, xde} for the description of the multiplicity
of the most continuous series representations for $G/H$ 
 in terms of  Weyl groups).
\newline\indent{2-b)}
Similarly to the case 2-a),
 the restriction
  $\pi|_H$ for
 a symmetric pair $(G,H)$ of non-holomorphic type 
 is multiplicity-free and is 
 decomposed into only continuous spectra 
 if $\pi$ is a holomorphic discrete series of scalar type.
This case was studied by G. \'Olafsson--B. \O rsted (\cite{xoo}).
\newline\indent{2-c)}
Theorem~\ref{thm:C} applied 
to non-discrete series representations
 $\pi_1$ and $\pi_2$
 (i.e. tensor products of \textit{singular}
  unitary highest weight representations)
 provides new settings of multiplicity-free branching laws.  
They might be interesting
 from the view point of representation theory because they construct
 \lq\lq small\rq\rq\ representations as discrete summands. 
(We note that irreducible unitary representations of
reductive Lie groups have not been classified even in the spherical case.
See \cite{xbarbasch} for the split case.) 
They might be interesting also from the view point of spectral theory and
harmonic analysis
 which is relevant to the
 {\it{canonical representation}} in the sense of Vershik--Gelfand--Graev.
Once we know the branching law is a priori multiplicity-free,
 it is promising to obtain its explicit formula.
Some special cases have been worked on in this direction so far,
 for $G= SL(2,\mathbb{R})$ by V. F. Molchanov \cite{xmol};  
 for $G= SU(2,2)$ by B. \O rsted and G. Zhang \cite{xoz};  
 for $G= SU(n,1)$ by G. van Dijk and S. Hille \cite{xvdh};
for $G = SU(p,q)$ by Y. Neretin and G. Ol'shanski\u\i\ 
\cite{xnere, xneol}.
See also G. van Dijk--M. Pevzner \cite{xDijkPev}, 
M. Pevzner \cite{xpevzner}
and G. Zhang \cite{xzhang}. 
Their results show that a different family of 
irreducible unitary representations
(sometimes, spherical complementary series representations)
 can occur in the same branching laws and each multiplicity is not
greater than one.
\end{exa}

\subsection{Organization of this article}
\label{subsec:1.9}
This paper is organized as follows:
In Section~\ref{sec:2}, %
 we  give a proof of an abstract multiplicity-free theorem
(Theorem~\ref{thm:2.2}) in the line bundle setting.
This is an extension of
 a theorem of Faraut--Thomas \cite{xft},
 whose idea may go back to Gelfand's proof 
 \cite{xgelsph}
 of the commutativity of the Hecke algebra $L^1(K \backslash G/K)$.
Theorem~\ref{thm:2.2} is a main method in this article to find various
multiplicity-free theorems.
In Section~\ref{sec:3},
we use Theorem~\ref{thm:2.2} to give a proof of Theorem~\ref{thm:A}.
The key idea is the reduction of the geometric condition
\eqref{eqn:2.2.3}
(\textit{strongly visible action} in the sense of \cite{RIMS})
to the existence problem of
 a \lq\lq nice\rq\rq\ involutive automorphism $\sigma$ of $G$
 satisfying a certain rank condition.
Section~\ref{sec:4} considers the multiplicity-free theorem for the
tensor product representations of two irreducible highest (or lowest)
weight modules and gives a proof of Theorem~\ref{thm:C}.
Sections~\ref{sec:5} and \ref{sec:6} examine our assumptions in our
multiplicity-free theorems (Theorems~\ref{thm:A} and \ref{thm:C}). 
That is, we drop the assumption of `scalar type' in
Section~\ref{sec:5} and prove that multiplicities are still uniformly
bounded (Theorems~\ref{thm:B} and \ref{thm:D}).
We note that multiplicities can be greater than one in this generality. 
In Section~\ref{sec:6},
we leave unchanged the assumption that $(G,H)$ is a symmetric pair,
and relax the assumption that $\pi$ is a highest weight module.
We
  illustrate by examples 
a wild behavior of multiplicities without this assumption.
In Section~\ref{sec:7},
 analogous results of Theorems~\ref{thm:A} and \ref{thm:C} are proved
 for finite dimensional representations of compact groups.
In Section~\ref{sec:8},
we present explicit branching laws that are assured a priori to be
multiplicity-free by Theorems~\ref{thm:A} and \ref{thm:C}.
Theorem~\ref{thm:tensordeco}
 generalizes the Hua--Kostant--Schmid formula.
In Section~\ref{sec:9} (Appendix)
 we present some basic results on homogeneous line bundles
 for the convenience of the reader,
 which give a sufficient condition for the assumption~\eqref{eqn:2.2.2} in 
Theorem~\ref{thm:2.2}.

\section{Main machinery from complex geometry}
\label{sec:2}

J. Faraut and E. Thomas \cite{xft},
 in the case of trivial twisting parameter,
 gave a sufficient condition for the commutativity
 of $\End_H(\mathcal{H})$ by using the theory of reproducing kernels,
 which we extend to the general,
 twisted case in this preliminary section.
The proof parallels to theirs,
 except that we need just find an additional condition
 \eqref{eqn:2.2.2} when we formalize Theorem~\ref{thm:2.2} in the line
 bundle setting.

\subsection{Basic operations on holomorphic line bundles}
\label{subsec:2.1}
Let $\mathcal{L} \to D$ be a holomorphic line bundle
 over a complex manifold $D$.
We denote by 
$\mathcal{O}(\mathcal{L}) \equiv \mathcal{O}(D,\mathcal{L})$ 
 the space of holomorphic sections of $\mathcal{L} \to D$.
Then $\mathcal{O}(\mathcal{L})$ carries a Fr\'echet topology
 by the uniform convergence on compact sets.
If a Lie group $H$ acts holomorphically and equivariantly
 on the holomorphic line bundle $\mathcal{L} \to D$, 
 then $H$ defines a (continuous) representation on $\mathcal{O}(\mathcal{L})$
 by the pull-back of sections.

Let $\{U_\alpha\}$ be trivializing neighborhoods of $D$,
 and $g_{\alpha \beta} \in \mathcal{O}^\times(U_\alpha \cap U_\beta)$
 the transition functions of the holomorphic line bundle $\mathcal{L} \to D$.
Then an anti-holomorphic line bundle
 $\overline{\mathcal{L}} \to D$ is a complex line bundle with
 the transition functions $\overline{g_{\alpha \beta}}$.
We denote by $\overline{\mathcal{O}}(\overline{\mathcal{L}})$
 the space of anti-holomorphic sections for $\overline{\mathcal{L}} \to D$.

Suppose $\sigma$ is an anti-holomorphic diffeomorphism of $D$.
Then the pull-back $\sigma^* \mathcal{L} \to D$ is an anti-holomorphic line
 bundle over $D$.
In turn,
 $\overline{\sigma^* \mathcal{L}} \to D$ is a holomorphic line bundle over $D$
 (see Appendix for more details).

\subsection{Abstract multiplicity-free theorem}
\label{subsec:2.2}
Here is the main machinery to prove various multiplicity-free theorems
of branching laws including
Theorems~\ref{thm:A} and \ref{thm:C}
(infinite dimensional representations) and Theorems~\ref{thm:E} and
\ref{thm:F} (finite dimensional representations). 

\begin{theorem}
\label{thm:2.2}
Let $(\pi,\mathcal{H})$ be a unitary representation of a Lie group $H$.
Assume that there exist
 an $H$-equivariant holomorphic line bundle $\mathcal{L} \to D$ and
 an  anti-holomorphic involutive diffeomorphism $\sigma$ of $D$
 with the following three conditions:

\begin{eq-text}
\label{eqn:2.2.1}
 There is an injective
 (continuous) $H$-intertwining map $\mathcal{H} \to \mathcal{O}(\mathcal{L})$.
\end{eq-text}
\begin{eq-text}
\label{eqn:2.2.2}
  There exists an isomorphism of $H$-equivariant holomorphic line bundles
   $\Psi: \mathcal{L} \overset{\sim}{\to} \overline{\sigma^* \mathcal{L}}$.
\end{eq-text}
\begin{eq-text}
\label{eqn:2.2.3}
  Given $x \in D$,
 there exists $g \in H$ such that $\sigma (x) = g \cdot x$.
\end{eq-text}

Then,
 the ring\/
$\End_H(\mathcal{H})$ of continuous $H$-intertwining operators on 
$\mathcal{H}$
is commutative.
Consequently, $(\pi, \mathcal{H})$ is multiplicity-free 
(see Definition~\ref{def:1.2}).
\end{theorem}

\subsection{Remarks on Theorem~\ref{thm:2.2}}

This subsection gives brief comments on Theorem~\ref{thm:2.2}.
First, we consider a special case,
and also a generalization.

\begin{rem}[specialization and generalization]%
1)\enspace
Suppose $\mathcal{L} \to D$ is the trivial line bundle.
Then, the condition \eqref{eqn:2.2.2} %
is automatically satisfied.
In this case,
 Theorem~\ref{thm:2.2} %
was proved in \cite{xft}.  
\newline
2)\enspace
An extension of  Theorem~\ref{thm:2.2} to the equivariant vector
bundle $\mathcal{V} \to D$ is the main subject of \cite{mfbdle},
where a more general multiplicity-free
 theorem is obtained under an additional condition
that the isotropy representation of
$H_x = \{ h \in H: h \cdot x = x \}$
on the fiber $\mathcal{V}_x$ is multiplicity-free for generic 
$x \in D$. 
Obviously, the $H_x$-action on $\mathcal{V}_x$ is multiplicity-free for
the case $\dim \mathcal{V}_x = 1$,
namely,
 for the line bundle case.
\end{rem}

Next, we examine the conditions \eqref{eqn:2.2.2} and
\eqref{eqn:2.2.3}. 

\begin{rem}
\label{rem:2.3.2}
1)\enspace
In many cases, the condition \eqref{eqn:2.2.2} is naturally
satisfied. 
We shall explicate how to construct the bundle isomorphism
$\Psi$
 in Lemma~\ref{lem:9.6} for a Hermitian symmetric space $D$.
\newline
2)\enspace
As the proof below shows,
Theorem~\ref{thm:2.2} still holds if we replace $D$ by an
$H$-invariant open subset $D'$.
Thus, the condition \eqref{eqn:2.2.3} is \textbf{local}.
The concept of `\textbf{visible action}' (see \cite{xkleiden, mfbdle,
xkgencar}) arises from the condition \eqref{eqn:2.2.3}
on the base space $D$.
\newline
3)\enspace
The condition \eqref{eqn:2.2.3} %
is automatically satisfied if $H$ acts transitively
 on $D$.
But we are interested in a more general setting where each $H$-orbit 
 has a positive codimension in $D$.
We find in Lemma~\ref{lem:3.2}
a sufficient condition for \eqref{eqn:2.2.3} %
in terms of rank condition 
 for a symmetric space $D$.
\end{rem}

\subsection{Reproducing kernel}
\label{subsec:2.4}
This subsection gives a quick summary for the reproducing kernel of a
Hilbert space $\mathcal{H}$ realized in the space
$\mathcal{O}(\mathcal{L})$
 of holomorphic
sections for a holomorphic line bundle $\mathcal{L}$
(see \cite{mfbdle} for a generalization to the vector bundle case). 
Since the reproducing kernel $K_{\mathcal{H}}$ contains all the
information on the Hilbert space $\mathcal{H}$,
our strategy is to make use of
$K_{\mathcal{H}}$ in order to prove Theorem~\ref{thm:2.2}.

Suppose that there is an injective and continuous map for a Hilbert
space $\mathcal{H}$ into the Fr\'{e}chet space
$\mathcal{O}(\mathcal{L})$. 
Then, 
 the point evaluation map
$$
   \mathcal{O}(\mathcal{L}) \supset \mathcal{H} \to \mathcal{L}_z \simeq \mathbb{C}\, , \quad f \mapsto f(z)
$$ 
 is continuous with respect to the Hilbert topology on $\mathcal{H}$.

Let $\{\varphi_\nu\}$ be an orthonormal basis of $\mathcal{H}$. 
We define
$$
   K_{\mathcal{H}}(x,y) \equiv  K(x,y)
 := \sum_{\nu} \varphi_\nu(x) \overline{\varphi_\nu(y)}
          \in 
        \mathcal{O}(\mathcal{L}) \widehat{\otimes} 
     \overline{\mathcal{O}}(\overline{\mathcal{L}}) \, .
$$
Then, $K(x,y)$ is well-defined as a holomorphic section of
$\mathcal{L} \to D$ 
for the first variable,
and as an anti-holomorphic section of
$\overline{\mathcal{L}} \to D$
for the second variable.
The definition is 
independent of the choice of
 an orthonormal basis $\{\varphi_\nu\}$.
$K(x,y)$ is called the reproducing kernel of
$\mathcal{H}$.

\begin{lemma}%
\label{lem:2.4}
{\rm 1)}\enspace
For each $y \in D$,
 $K(\cdot, y) \in \mathcal{H} \otimes \overline{\mathcal{L}_y} \ (\simeq \mathcal{H})$
 and 
 $(f(\cdot), K(\cdot,y))_{\mathcal{H}} = f(y)$ for any $f \in \mathcal{H}$.
\itm{2)}
Let $K_i(x,y)$ 
 be the reproducing kernels of
 Hilbert spaces $\mathcal{H}_i \ \subset \mathcal{O}(\mathcal{L})$
 with inner products $(\ , \ )_{\mathcal{H}_i}$,
 respectively, for $i = 1, 2$.
If $K_1 \equiv K_2$,
 then $\mathcal{H}_1 = \mathcal{H}_2$ and $(\ , \ )_{\mathcal{H}_1} = (\ , \ )_{\mathcal{H}_2}$.
\itm{3)}
If $K_1(x,x) = K_2(x,x)$ for any $x \in D$,
 then $K_1 \equiv K_2$.
\end{lemma}

\begin{proof}
(1) and (2) are standard.
We review only the way how to recover $\mathcal{H}$ together
with its inner product
 from a given reproducing kernel.
For each $y \in D$,
 we fix an isomorphism $\mathcal{L}_y \simeq \mathbb{C}$.
Through this isomorphism,
 we can regard
 $K(\cdot, y) \in \mathcal{H} \otimes \overline{\mathcal{L}_{y}}$
 as an element of $\mathcal{H}$.
The Hilbert space $\mathcal{H}$ is the completion of
 the $\mathbb{C}$-span of $\set{K(\cdot, y)}{y \in D}$ 
 with pre-Hilbert structure
\begin{equation}
  (K(\cdot, y_1), K(\cdot, y_2))_{\mathcal{H}} := K(y_2, y_1) \in 
   \mathcal{L}_{y_2} \otimes \overline{\mathcal{L}_{y_1}} 
   \ (\simeq \mathbb{C}) \, .
\label{eqn:2.4.1}
\end{equation}
This procedure is
 independent of the choice of 
 the isomorphism ${\mathcal{L}}_y \simeq \mathbb{C}$.
Hence, the Hilbert space $\mathcal{H}$ together with its inner product
 is recovered.
\itm{3)}
We denote by $\overline{D}$ the complex manifold
 endowed with the conjugate complex structure on $D$.
Then,
 $\overline{\mathcal{L}} \to \overline{D}$ is a holomorphic line bundle,
 and $K(\cdot,\cdot) \equiv K_{\mathcal{H}}(\cdot,\cdot)$ 
 is a holomorphic section of the holomorphic line bundle
 $\mathcal{L} \boxtimes \overline{\mathcal{L}} \to D \times \overline{D}$. 
As the diagonal embedding
 $\iota: D \to D \times \overline{D}, z \mapsto (z, z)$
 is totally real,
 $(K_1-K_2)|_{\iota(D)} \equiv 0$ implies 
 $K_1-K_2\equiv 0$ by the unicity theorem of holomorphic functions.
\qed
\end{proof}

\subsection{Construction of $J$}
\label{subsec:2.5}
Suppose we are in the setting of Theorem~\ref{thm:2.2}. %
We define an anti-linear map
$$
J : \mathcal{O}(\mathcal{L}) \to \mathcal{O}(\mathcal{L})\, ,
 \quad f\mapsto J f 
$$
by $J f(z) := \overline{f(\sigma(z))}$ ($z \in D$).
$Jf$ is regarded as an element of $\mathcal{O}(\mathcal{L})$
through the isomorphism 
$
 \Psi_* : \mathcal{O}(\mathcal{L}) \simeq 
 \mathcal{O}(\overline{\sigma^* \mathcal{L}})
$
(see \eqref{eqn:2.2.2}).
\begin{lemma}
\label{lem:2.5}
In the setting of Theorem~\ref{thm:2.2}, %
 we identify $\mathcal{H}$ with a subspace of $\mathcal{O}(\mathcal{L})$.
Then, the anti-linear map 
 $J$ is an isometry from $\mathcal{H}$ onto $\mathcal{H}$.
\end{lemma}
\begin{proof}
We put $\widetilde{\mathcal{H}} := J(\mathcal{H})$,
 equipped with the inner product
\begin{equation}
  (J f_1, J f_2)_{\widetilde{\mathcal{H}}} := (f_2, f_1)_{\mathcal{H}}
 \quad \text{for } f_1, f_2 \in \mathcal{H} \, .
\label{eqn:2.5.1}
\end{equation}
If $\{\varphi_\nu\}$ is an orthonormal basis of $\mathcal{H}$,
 then $\widetilde{\mathcal{H}}$ is a Hilbert space with
 orthonormal basis $\{J \varphi_\nu\}$.
Hence,
 the reproducing kernel of $\widetilde{\mathcal{H}}$ is given by
 $K_{\widetilde{\mathcal{H}}}(x,y) = K_{\mathcal{H}}(\sigma(y), \sigma(x))$
 because
\begin{equation}
 K_{\widetilde{\mathcal{H}}}(x,y) 
 =\sum_\nu J \varphi_\nu(x) \overline{J \varphi_\nu(y)}
 =\sum_\nu  \overline{\varphi_\nu(\sigma(x))} \ 
 \overline{\overline{\varphi_\nu(\sigma(y))}}
= K_{\mathcal{H}}(\sigma(y), \sigma(x)) \, .
\label{eqn:2.5.2}
\end{equation}
We fix $x \in D$ and
 take $g \in H$ such that $\sigma (x) = g \cdot x$
 (see \eqref{eqn:2.2.3}). %
Substituting $x$ for $y$ in \eqref{eqn:2.5.2}, %
 we have
$$
 K_{\widetilde{\mathcal{H}}}(x,x) 
= K_{\mathcal{H}}(\sigma (x), \sigma (x))
= K_{\mathcal{H}}(g \cdot x, g \cdot x)
= K_{\mathcal{H}}(x,  x) \, .
$$
Here,
 the last equality holds because
$\{ \varphi_\nu (g \cdot {}) \}$ is also an orthonormal basis
 of $\mathcal{H}$ as $(\pi, \mathcal{H})$
 is a unitary representation of $H$.
Then,
 by Lemma~\ref{lem:2.4}, %
 the Hilbert space $\widetilde{\mathcal{H}}$ coincides with $\mathcal{H}$
  and 
\begin{equation}
  (J f_1, J f_2)_{\mathcal{H}} = (f_2, f_1)_{\mathcal{H}} 
 \quad \text{for } f_1, f_2 \in \mathcal{H}\, .
\label{eqn:2.5.3}
\end{equation}
This is what we wanted to prove.
\qed
\end{proof}

\subsection{Proof of $A^* = JAJ^{-1}$}
\begin{theopargself}
\begin{lemma}[\textmd{(see \cite{xft})}]%
\label{lem:2.6}
Suppose $A \in \End_H(\mathcal{H})$.  %
Then the adjoint operator  $A^*$  of $A$ is given by
\begin{equation}
    A^* = J A J^{-1} \, .
\label{eqn:JAJ}
\end{equation}
\end{lemma}
\end{theopargself}

\begin{proof}
We divide the proof into three steps.
\newline
$\underline {\text{Step 1}}$ (positive self-adjoint case):\enspace
Assume $A \in \End_H(\mathcal{H})$ is a positive self-adjoint operator.
Let $\mathcal{H}_A$ be the Hilbert completion
 of $\mathcal{H}$ by the pre-Hilbert structure
\begin{equation}
  (f_1, f_2)_{\mathcal{H}_A} := (A f_1, f_2)_{\mathcal{H}} 
 \quad \text{for } f_1, f_2 \in \mathcal{H} \, .  
\label{eqn:2.6.2}
\end{equation}
If $f_1, f_2 \in \mathcal{H}$ and $g \in H$,
 then
\begin{align*}
   (\pi(g) f_1, \pi(g) f_2)_{{\mathcal{H}}_A}
  &= (A \pi(g) f_1, \pi(g) f_2)_{\mathcal{H}}
\\
  &= (\pi(g) A f_1, \pi(g) f_2)_{\mathcal{H}}
  = (A f_1,  f_2)_{\mathcal{H}}
  = (f_1,  f_2)_{{\mathcal{H}}_A} \, .
\end{align*}
Therefore,
 $(\pi, \mathcal{H})$ extends to a unitary representation on ${\mathcal{H}}_A$.
Applying (2.5.3) to both $\mathcal{H}_A$ and $\mathcal{H}$,
 we have 
\begin{multline*}
    (A f_1, f_2)_{\mathcal{H}}
 =  (f_1, f_2)_{\mathcal{H}_A}
 =  (J f_2, J f_1)_{\mathcal{H}_A}
 =  (A J f_2, J f_1)_{\mathcal{H}}
\\
 =  (J f_2, A^* J f_1)_{\mathcal{H}}
 =  (J f_2, J J^{-1} A^* J f_1)_{\mathcal{H}}
 =  (J^{-1} A^* J f_1, f_2)_{\mathcal{H}} \, .
\end{multline*}
Hence, $A = J^{-1} A^* J$,
 and \eqref{eqn:JAJ} %
follows.

\noindent
$\underline{\text{Step 2}}$ (self-adjoint case):\enspace
Assume $A \in \End_H(\mathcal{H})$ is a self-adjoint operator.
Let $A = \int \lambda d E_\lambda$ be the spectral decomposition of $A$.
Then every projection operator $E_\lambda \in \End(\mathcal{H})$ also
 commutes with $\pi(g)$ for all $g \in H$,
 namely, 
 $E_\lambda \in \End_H(\mathcal{H})$.
We define
$$
 A_+ := \int_{\lambda \ge 0} \lambda d E_\lambda \, ,
 \qquad
 A_- := \int_{\lambda< 0} \lambda d E_\lambda \, .
$$
Then $A = A_+ + A_-$.
Let $I$ be the identity operator on $\mathcal{H}$.
As a positive self-adjoint operator $A_+ + I$
 is an element of $\End_H(\mathcal{H})$,
 we have 
$
    (A_+ + I)^* = J (A_+ + I) J^{-1}
$
 by Step 1,
whence 
$
    A_+^* = J A_+ J^{-1}.
$
Applying Step 1 again to $- A_-$,
 we have
$
    A_-^* = J A_- J^{-1}.
$
Thus,
$$
  A^* = A_+^* + A_-^* = J A_+ J^{-1} + J A_- J^{-1}
  = J (A_+ + A_-) J^{-1}
 = J A J^{-1} \, .
$$

\noindent
$\underline{\text{Step 3}}$ (general case):\enspace
Suppose $A \in \End_H(\mathcal{H})$. 
Then $A^*$ also commutes with $\pi(g)$ ($g \in H$) because $\pi$ is unitary.
We put $B := \frac{1}{2}(A+ A^*)$ and $C := \frac{\sqrt{-1}}{2}(A^* - A)$.
Then,
 both $B$ and $C$ are self-adjoint operators commuting with $\pi(g)$
 ($g \in H$).
It follows from Step 2 that
$
    B^* = J B J^{-1}
$
and
$
    C^* = J C J^{-1} .
$
As $J$ is an anti-linear map,
 we have
$$
   (\sqrt{-1}\, C)^* 
 = - \sqrt{-1}\, C^* 
 = - \sqrt{-1}\, J C J^{-1}
 = J (\sqrt{-1}\, C) J^{-1} \, .
$$
Hence,
 $A = B + \sqrt{-1}\, C$ also satisfies
$
    A^* = J A J^{-1}.
$
\qed
\end{proof}

\subsection{Proof of Theorem~\ref{thm:2.2}}%
\label{subsec:2.7}
We are now ready to complete the proof of Theorem~\ref{thm:2.2}.
Let $A$, $B \in \End_H(\mathcal{H})$. %
By Lemma~\ref{lem:2.6}, %
 we have
\begin{align*}
&  AB = J^{-1}(AB)^* J = (J^{-1} B^* J) (J^{-1} A^* J) = BA \, .
\end{align*}
Therefore,
 $\End_H(\mathcal{H})$  is commutative.
\qed

\section{Proof of Theorem~\protect\ref{thm:A}}
\label{sec:3}
This section gives a proof of Theorem~\ref{thm:A}
by using Theorem~\ref{thm:2.2}.
The core of the proof is to reduce the geometric condition
\eqref{eqn:2.2.3} to an algebraic condition
(the existence of a certain involution of the Lie algebra).
This reduction is stated in Lemma~\ref{lem:3.2}.
The reader who is familiar with symmetric pairs can skip
Subsections~\ref{subsec:3.1}, \ref{subsec:3.2.ex}, \ref{subsec:3.3} and
\ref{subsec:3.5}.

\subsection{Reductive symmetric pairs}
\label{subsec:3.1}
Let $G$ be a Lie group.
Suppose that $\tau$ is an involutive automorphism of $G$.
We write
$$G^\tau := \{ g \in G: \tau g = g \}$$
for the fixed point subgroup of $\tau$,
and denote by $G_0^\tau$ its connected component containing the unit
element. 
The pair $(G,H)$ 
(or the pair $(\mathfrak{g}, \mathfrak{h})$ of their Lie algebras) 
is called a \textit{symmetric pair}
 if the subgroup $H$ is an open subgroup of $G^\tau$,
that is, if $H$ satisfies
$$
G_0^\tau \subset H \subset G^\tau .
$$
It is called a \textit{reductive symmetric pair}
 if $G$ is a reductive Lie
group; 
 a \textit{semisimple symmetric pair}
 if $G$ is a semisimple Lie group.   
Obviously, a semisimple symmetric pair is a reductive symmetric pair.

We shall use the same letter $\tau$ to denote the differential of $\tau$.
We set
$$
\mathfrak{g}^{\pm \tau} := \set{Y \in \mathfrak{g}}{\tau Y = \pm Y}\,.
$$
Then,
 it follows from $\tau^2 =\operatorname{id}$
 that we have a direct sum decomposition 
$$
\mathfrak{g} = \mathfrak{g}^\tau \oplus \mathfrak{g}^{-\tau}.
$$

Suppose now
 that $G$ is a semisimple Lie group.  
It is known that there exists
 a Cartan involution $\theta$ of $G$ commuting with $\tau$.
Take such $\theta$, and we write
 $K := G^\theta =\set{g \in G}{\theta g = g}$.
Then, $K$ is compact if $G$ is a linear Lie group. 
The direct sum decomposition 
$$
\mathfrak{g} = \mathfrak{k}\oplus \mathfrak{p} 
\equiv \mathfrak{g}^\theta \oplus \mathfrak{g}^{-\theta}
$$
is called a Cartan decomposition. 
Later, 
we shall allow $G$ to be non-linear,
 in particular,
 $K$ is not necessarily compact.
The \textit{real rank} of $\mathfrak{g}$,
 denoted by $\rrank \mathfrak{g}$,
 is defined to be the dimension of a maximal abelian subspace of
 $\mathfrak{g}^{-\theta}$.

As $(\tau\theta)^2 = \operatorname{id}$,
the pair $(\mathfrak{g}, \mathfrak{g}^{\tau\theta})$
also forms a symmetric pair.
The Lie group
$$
G^{\tau\theta} = \{ g \in G : (\tau\theta)(g) = g \}
$$
is a reductive Lie group with Cartan involution
$\theta|_{G^{\tau\theta}}$,
and its Lie algebra $\mathfrak{g}^{\tau\theta}$ is reductive with
Cartan decomposition
\begin{equation}
\label{eqn:gtaut}
\mathfrak{g}^{\tau\theta} = \mathfrak{g}^{\tau\theta,\theta}
   \oplus \mathfrak{g}^{\tau\theta,-\theta}
 = \mathfrak{g}^{\tau,\theta} \oplus \mathfrak{g}^{-\tau,-\theta} \, .
\end{equation}
Here, we have used the notation
$\mathfrak{g}^{-\tau,-\theta}$ and alike,
defined as follows:
$$
\mathfrak{g}^{-\tau,-\theta} 
  := \{ Y \in \mathfrak{g} :
     (-\tau) Y = (-\theta) Y = Y \} \, .
$$
Then, the dimension of a maximal abelian subspace $\mathfrak{a}$ of
$\mathfrak{g}^{-\tau,-\theta}$ is equal to the real rank of
$\mathfrak{g}^{\tau\theta}$, 
which is referred to as the \textit{split rank} 
of the semisimple symmetric space $G/H$.
We shall write $\rrank G/H$ or $\rrank \mathfrak{g}/\mathfrak{g}^\tau$
for this dimension.
Thus, 
\begin{equation}
    \rrank \mathfrak{g}^{\theta\tau} = \rrank \mathfrak{g}/\mathfrak{g}^\tau
 \, .
\label{eqn:3.1.1}
\end{equation}
In particular, we have
 $\rrank \mathfrak{g} = \rrank \mathfrak{g}/\mathfrak{k}$ 
if we take $\tau$ to be $\theta$.

The Killing form on the Lie algebra $\mathfrak{g}$ is non-degenerate
on $\mathfrak{g}$, 
and is also non-degenerate when restricted to ${\mathfrak {h}}$.
Then, it induces an $\Ad(H)$-invariant non-degenerate bilinear form on
$\mathfrak{g}/\mathfrak{h}$,
and therefore a $G$-invariant pseudo-Riemannian structure on the
homogeneous space $G/H$,
so that $G/H$ becomes a symmetric space with respect to the 
Levi--Civita connection
 and is called a \textit{semisimple symmetric space}.  
In this context, 
the subspace ${\mathfrak {a}}$
 has the following geometric meaning: 
Let $A := \exp(\mathfrak{a})$,
the connected abelian subgroup of $G$ with Lie algebra $\mathfrak{a}$.
Then, the orbit $A \cdot o$ through
$o := eH \in G/H$
becomes a flat,
totally geodesic submanifold in $G/H$.
Furthermore, we have a (generalized) Cartan decomposition:

\begin{fct}[\textmd{see \cite[Section 2]{xfjd}}]
\label{fact:genCar}
$G = KAH$.
\end{fct}

{\renewcommand{\proofname}{Sketch of Proof}
\begin{proof}
The direct sum decomposition of the Lie algebra
$$
\mathfrak{g} = \mathfrak{k} \oplus \mathfrak{g}^{-\tau,-\theta}
     \oplus \mathfrak{g}^{\tau,-\theta}
$$
lifts to a diffeomorphism:
$$
\mathfrak{g}^{-\tau,-\theta} + \mathfrak{g}^{\tau,-\theta}
  \overset{\sim}{\to}
  K\backslash G\, ,
\quad
  (X,Y) \mapsto K e^X e^Y.
$$
Since $\exp(\mathfrak{g}^{\tau,-\theta}) \subset H$,
the decomposition $G = KAH$ follows if we show
\begin{equation}
\label{eqn:HKa}
\Ad(H \cap K) \mathfrak{a} 
= \mathfrak{g}^{-\tau,-\theta}.
\end{equation}
The equation \eqref{eqn:HKa} is well-known
 as the key ingredient of the original Cartan decomposition
$G^{\tau\theta} = K^\tau A K^\tau$
in light of \eqref{eqn:gtaut}.
\qed
\end{proof}
}

Furthermore, suppose that 
 $\sigma$ is an involutive automorphism of $G$ such that
 $\sigma$, $\tau$ and $\theta$ commute with one another.
We set
$$
 G^{\sigma, \tau} := G^\sigma \cap G^\tau 
 =
 \set{g \in G}{\sigma g = \tau g = g} \, .
$$
Then
  $(G^\sigma, G^{\sigma, \tau})$ forms a reductive symmetric pair, 
 because $\sigma$ and $\tau$ commute.
The commutativity of $\sigma$ and $\theta$ implies that
 the automorphism $\sigma: G \to G$ 
 stabilizes $K$ and induces a diffeomorphism of $G/K$,
 for which we use the same letter $\sigma$.

\subsection{Examples of symmetric pairs}
\label{subsec:3.2.ex}
This subsection presents some basic examples of semisimple
(and therefore, reductive) symmetric
pairs. 

\begin{exam}[group manifold]
\label{ex:gpmfd}
Let $G'$ be a semisimple Lie group,
and $G := G' \times G'$.
We define an involutive automorphism $\tau$ of $G$ by
$\tau(x,y) := (y,x)$.
Then, $G^\tau = \{ (g,g): g \in G' \}$
is the diagonal subgroup,
denoted by $\diag(G')$,
which is isomorphic to $G'$.
Thus, $(G' \times G', \diag(G'))$
forms a semisimple symmetric pair.

We set
\begin{align*}
I_{p,q} &:=
  \begin{pmatrix}
   \;
     \begin{matrix}
            1 \\ &\ \rotatebox{-5}{$\ddots$} \\  &&\ 1
     \end{matrix}
      \kern-3em\smash{\raisebox{4.9ex}{\rotatebox{-45}{$\overbrace{~~~~~~~~~~~~~}$}} }
             \kern-1.5em\smash{\raisebox{3.2ex}{$p$}}
     &\mbox{\Huge$0$}
  \\
    \mbox{\Huge$0$}
    &\begin{matrix}
           \kern-.1em -1 \\ &\rotatebox{-5}{$\ddots$} \\  &&\kern-.1em -1
     \end{matrix}
      \kern-2.8em\smash{\raisebox{4.8ex}{\rotatebox{-45}{$\overbrace{~~~~~~~~~~~~~}$}} }
             \kern-1.5em\smash{\raisebox{3.2ex}{$q$}}
   \kern1em
  \end{pmatrix}
\\
J &:=
   \begin{pmatrix}
     \parbox[c]{3em}{\hfil\Huge$0$\hfil}
     & \parbox[c]{3em}{\hfil\Large$I_n$\hfil}
     \\[3	ex]
     \parbox[c]{3em}{\hfil\Large$-I_n$\hfil}
     &\parbox[c]{3em}{\hfil\Huge$0$\hfil}
   \end{pmatrix}
\end{align*}
\end{exam}

\begin{exam}
\label{ex:SLSU}
Let $G = SL(n,\mathbb{C})$,
and fix $p,q$ such that $p+q=n$.
Then,
$$
\tau(g) := I_{p,q} \, g^* I_{p,q}
\quad (g \in G)
$$
defines an involutive automorphism of $G$,
and $G^\tau = SU(p,q)$
(the indefinite unitary group).
Thus, $(SL(n,\mathbb{C}), SU(p,q))$
forms a semisimple symmetric pair.
\end{exam}

\begin{exam}
\label{ex:SLSL}
Let $G = SL(n,\mathbb{C})$,
and $\sigma(g) := \overline{g}$.
Then $\sigma$ is an involutive automorphism of $G$,
and $G^\sigma = SL(n,\mathbb{R})$.
We note that $\sigma$ commutes with the involution $\tau$ in the
previous example,
and
\begin{align*}
G^{\sigma,\tau}& = \{ g \in SL(n,\mathbb{C}):
     \overline{g} = g = I_{p,q} \, \overline{{}^t g} \, I_{p,q} \}
\\
   & = SO(p,q) \, .
\end{align*}
Thus, $(SL(n,\mathbb{C}), SL(n,\mathbb{R}))$,
$(SU(p,q), SO(p,q))$, $(SL(n,\mathbb{R}), SO(p,q))$
are examples of semisimple symmetric pairs.
\end{exam}

\begin{exam}
\label{ex:SLSp}
Let $G := SL(2n,\mathbb{R})$,
and $\tau(g) := J \, {}^t g^{-1} J^{-1}$.
Then,
$G^\tau = Sp(n,\mathbb{R})$
(the real symplectic group).
Thus, $(SL(2n,\mathbb{R}), Sp(n,\mathbb{R}))$ forms a semisimple
symmetric pair.
\end{exam}

\subsection{Reduction of visibility to real rank condition}
\label{subsec:3.2}
The following lemma gives a sufficient condition for 
\eqref{eqn:2.2.3}.
Then, it plays a key role when we apply
Theorem~\ref{thm:2.2}
 to the branching problem for the restriction from $G$ to $G^\tau$
(with the notation of Theorem~\ref{thm:2.2}, 
$D = G/K$ and $H = G_0^\tau$).
This lemma is also used in reducing
`visibility' of an action to an algebraic condition
(\cite[Lemma~2.2]{visiblesymm}). 

\begin{lemma}
\label{lem:3.2}
Let $\sigma$ and $\tau$ be involutive automorphisms of $G$.
We assume that the pair $(\sigma, \tau)$ satisfies
 the following two conditions:

\begin{eq-text}
  $\sigma$, $\tau$ and $\theta$ commute with one another.
\label{eqn:3.2.1}
\end{eq-text}
\begin{eq-text}
  $\rrank \mathfrak{g}^{\tau\theta}
 = \rrank \mathfrak{g}^{\sigma, \tau\theta}$.
\label{eqn:3.2.2}
\end{eq-text}

Then for any $x \in G/K$, there exists $g \in G^\tau_0$ such that
$
      \sigma (x) = g \cdot x
$.
\end{lemma}

\begin{proof}
It follows from the condition \eqref{eqn:3.2.1} %
 that
 $\theta|_{G^\sigma}$ is a Cartan involution of
 a reductive Lie group $G^\sigma$
 and 
 that  $\tau|_{G^\sigma}$ is an involutive automorphism of $G^\sigma$
 commuting with $\theta|_{G^\sigma}$.
Take a maximal abelian subspace $\mathfrak{a}$ in
$$
   \mathfrak{g}^{-\theta, \sigma, \tau\theta}
   :=
   \set{Y \in \mathfrak{g}}{(-\theta) Y = \sigma Y = \tau\theta Y = Y} \,
   .
$$
From definition, 
we have 
 $\dim \mathfrak{a} = \rrank \mathfrak{g}^{\sigma, \tau\theta}$,
which in turn equals
$\rrank \mathfrak{g}^{\tau\theta}$ by the condition
\eqref{eqn:3.2.2}.
This means that $\mathfrak{a}$ is also a maximal abelian subspace in
$$
   \mathfrak{g}^{-\theta, \tau\theta}
   =
   \set{Y \in \mathfrak{g}}{(-\theta) Y = \tau\theta Y = Y} \, .
$$
Let $A = \exp(\mathfrak{a})$.
Then it follows from Fact~\ref{fact:genCar} that
we have a generalized Cartan decomposition
\begin{equation}
   G = G^\tau_0 A K \, .
\label{eqn:3.2.3}
\end{equation}
Let $o := e K \in G/K$.
Fix $x \in G/K$.
Then,
 according to the decomposition \eqref{eqn:3.2.3}, %
 we find $h \in G^\tau_0$ and $a \in A$ such that
$$
        x = h a \cdot o \, .
$$
We set $g := \sigma(h)\; h^{-1}$.
We claim $g \in G^\tau_0$.
In fact,
 by using $\sigma \tau = \tau \sigma$ and $\tau h = h$,
 we have
$$
 \tau (g) = \tau \sigma(h)\; \tau(h^{-1}) 
          = \sigma \tau (h)\; \tau(h)^{-1}
          = \sigma (h)\; h^{-1} = g \, .
$$
Hence, 
 $g \in G^\tau$.
Moreover,
 since the image of the continuous map
$$
     G^\tau_0 \to G\, , 
  \quad h \mapsto \sigma(h)\; h^{-1}
$$
 is connected,
 we have $g \in G^\tau_0$.

On the other hand,
 we have $\sigma(a) = a$ because $\mathfrak{a} \subset
   \mathfrak{g}^{-\theta, \sigma, -\tau} \subset
 \mathfrak{g}^\sigma$.
Therefore we have
$$
    \sigma (x) = \sigma(h)\; \sigma(a) \cdot o
               = \sigma(h)\; h^{-1} h a \cdot o
               = g \cdot x \, ,
$$
 proving the lemma.
\qed
\end{proof}

\subsection{Hermitian Symmetric Space $G/K$}
\label{subsec:3.3}
Throughout the rest of this section,
 we assume that $G$ is a simple, non-compact,
 Lie group
of Hermitian type.
We retain the notation of 
Subsection~\ref{subsec:1.4}. 

Let $G_\mathbb{C}$ be a connected complex Lie group with Lie algebra
 $\mathfrak{g}_\mathbb{C}$,
 and $Q^-$ the maximal parabolic subgroup of $G_\mathbb{C}$ 
 with Lie algebra $\mathfrak{k}_\mathbb{C} + \mathfrak{p}_-$.
Then we have an open
  embedding $G/K \hookrightarrow G_\mathbb{C}/Q^-$
because 
$\mathfrak{g}_{\mathbb{C}} = \mathfrak{g} +
    (\mathfrak{k}_{\mathbb{C}} + \mathfrak{p}_-)$.
Thus, a $G$-invariant complex structure on $G/K$ is induced
from $G_{\mathbb{C}} / Q^-$.
(We remark that the embedding $G/K \hookrightarrow G_\mathbb{C}/Q^-$
 is well-defined,  
 even though $G$ is not necessarily a subgroup of $G_{\mathbb{C}}$.)

Suppose $\tau$ is an involutive automorphism of $G$
 commuting with $\theta$.
We recall from Subsection~\ref{subsec:1.5} that
we have either
\begin{alignat*}{2}
    \tau Z &= Z 
  &\quad&\mbox{(holomorphic type),}
\tag{\ref{eqn:1.5.1}}
\\
\intertext{or}
   \tau Z &= -Z 
  &\quad&\mbox{(anti-holomorphic type).}
\tag{\ref{eqn:1.5.2}}
\end{alignat*}
Here is the classification of semisimple symmetric pairs
 $(\mathfrak{g}, \mathfrak{g}^\tau)$ 
 with $\mathfrak{g}$ simple
 such that the pair $(\mathfrak{g}, \mathfrak{g}^\tau)$
 satisfies the condition \eqref{eqn:1.5.1}
 (respectively, \eqref{eqn:1.5.2}).
Table~\ref{tbl:3.3.2} %
is equivalent to the classification of
 totally real symmetric spaces $G^\tau/K^\tau$ of the Hermitian
 symmetric space $G/K$
 (see \cite{xfo, xjafbams, xjafjdg, xkobanaga}).

\begin{table}[H]
\caption{}
\label{tbl:3.3.1}
$$
\vbox{
\offinterlineskip
\def\tablerule{\noalign{\hrule}}
\def\ct{&\cr\tablerule}
\halign{\strut#&\vrule#&
            \;\;\hfil#\hfil\hfil\;\;&\vrule#&
            \;\;\hfil#\hfil\hfil\;\;&\vrule#\cr\tablerule
&&\multispan3\hfil $(\mathfrak{g}, \mathfrak{g}^\tau)$ is of holomorphic 
type
 \hfil \ct
&& $\mathfrak{g}$ && $\mathfrak{g}^\tau$   \ct
&&  $\mathfrak{su}(p,q)$     && $\hphantom{mmi}\mathfrak{s}(\mathfrak{u}(i,j) + \mathfrak{u}(p-i,q-j))$ \ct
&&  $\mathfrak{su}(n,n)$     && $\hphantom{mmi}\mathfrak{so}^*(2n)$ \ct
&&  $\mathfrak{su}(n,n)$     && $\hphantom{mmm}\mathfrak{sp}(n, \mathbb{R})$ \ct
&&  $\mathfrak{so}^*(2n)$    && $\hphantom{mm}\mathfrak{so}^*(2p) + \mathfrak{so}^*(2n-2p)$ \ct
&&  $\mathfrak{so}^*(2n)$    && $\hphantom{mmm}\mathfrak{u}(p,n-p)$ \ct
&&  $\mathfrak{so}(2,n)$     && $\hphantom{m,}\mathfrak{so}(2,p)+\mathfrak{so}(n-p)$ \ct
&&  $\mathfrak{so}(2,2n)$    && $\hphantom{mm}\mathfrak{u}(1,n)$ \ct
&&  $\mathfrak{sp}(n,\mathbb{R})$&& $\hphantom{mmm}\mathfrak{u}(p,n-p)$ \ct
&&  $\mathfrak{sp}(n,\mathbb{R})$&& $\hphantom{mi}\mathfrak{sp}(p,\mathbb{R})+\mathfrak{sp}(n-p,\mathbb{R})$\ct
&&  $\mathfrak{e}_{6(-14)}$  && $\hphantom{,}\mathfrak{so}(10)+\mathfrak{so}(2)$\ct
&&  $\mathfrak{e}_{6(-14)}$  && $\mathfrak{so}^*(10)+\mathfrak{so}(2)$\ct
&&  $\mathfrak{e}_{6(-14)}$  && $\mathfrak{so}(8,2)+\mathfrak{so}(2)$\ct
&&  $\mathfrak{e}_{6(-14)}$  && $\hphantom{i}\mathfrak{su}(5,1)+\mathfrak{sl}(2, \mathbb{R})$\ct
&&  $\mathfrak{e}_{6(-14)}$  && $\mathfrak{su}(4,2)+\mathfrak{su}(2)$\ct
&&  $\mathfrak{e}_{7(-25)}$  && $\hphantom{i}\mathfrak{e}_{6(-78)} +\mathfrak{so}(2)$\ct
&&  $\mathfrak{e}_{7(-25)}$  && $\hphantom{i}\mathfrak{e}_{6(-14)} +\mathfrak{so}(2)$\ct
&&  $\mathfrak{e}_{7(-25)}$  && $\mathfrak{so}(10,2) +\mathfrak{sl}(2,\mathbb{R})$\ct
&&  $\mathfrak{e}_{7(-25)}$  && $\mathfrak{so}^*(12) +\mathfrak{su}(2)$\ct
&&  $\mathfrak{e}_{7(-25)}$  && $\hphantom{mi}\mathfrak{su}(6,2)$\ct
\hfil\cr}}
$$
\end{table}

\begin{table}[H]
\caption{}
\label{tbl:3.3.2}
$$
\vbox{
\offinterlineskip
\def\tablerule{\noalign{\hrule}}
\def\ct{&\cr\tablerule}
\halign{\strut#&\vrule#&
            \;\;\hfil#\hfil\hfil\;\;&\vrule#&
            \;\;\hfil#\hfil\hfil\;\;&\vrule#\cr\tablerule
&&\multispan3\hfil $(\mathfrak{g}, \mathfrak{g}^\tau)$ is of
anti-holomorphic type
\hfil \ct
&& $\mathfrak{g}$ && \quad $\mathfrak{g}^\tau$   \ct
&&  $\mathfrak{su}(p,q)$    && $\hphantom{mm}\mathfrak{so}(p,q)$ \ct
&&  $\mathfrak{su}(n,n)$     && $\mathfrak{sl}(n,\mathbb{C}) + \mathbb{R}$ \ct
&&  $\mathfrak{su}(2p,2q)$     && \hphantom{mm}$\mathfrak{sp}(p,q)$ \ct
&&  $\mathfrak{so}^*(2n)$    && $\hphantom{mm}\mathfrak{so}(n,\mathbb{C})$ \ct
&&  $\mathfrak{so}^*(4n)$    && $\mathfrak{su}^*(2n) + \mathbb{R}$ \ct
&&  $\mathfrak{so}(2,n)$     && $\hphantom{mm}\,\,\,\mathfrak{so}(1,p)+\mathfrak{so}(1,n-p)$ \ct
&&  $\mathfrak{sp}(n,\mathbb{R})$&& $\hphantom{mm}\,\,\mathfrak{gl}(n,\mathbb{R})$ \ct
&&  $\mathfrak{sp}(2n,\mathbb{R})$&& $\hphantom{mm}\mathfrak{sp}(n,\mathbb{C})$\ct
&&  $\mathfrak{e}_{6(-14)}$  && $\hphantom{mm}\mathfrak{f}_{4(-20)}$\ct
&&  $\mathfrak{e}_{6(-14)}$  && $\hphantom{mm}\mathfrak{sp}(2,2)$\ct
&&  $\mathfrak{e}_{7(-25)}$  && $\phantom{mi}\;\mathfrak{e}_{6(-26)} +\mathfrak{so}(1,1)$\ct
&&  $\mathfrak{e}_{7(-25)}$  && $\hphantom{mmmi}\mathfrak{su}^*(8)$ \ct
\hfil\cr}}
$$
\end{table}

\subsection{Holomorphic realization of highest weight representations}
\label{subsec:3.5}
It is well-known that an irreducible  highest weight representation
$\pi$ of $G$ can be realized as a subrepresentation of the space
 of global holomorphic sections of an equivariant
 holomorphic vector bundle over the Hermitian symmetric space
$G/K$.
We supply a proof here for the convenience of the reader
in a way that we shall use later.

\begin{lemma}
\label{lem:3.5}
Let $(\pi, \mathcal{H})$ 
 be an irreducible unitary highest weight module.
We write $\chi$ for the representation of $K$ on 
$U := \mathcal{H}_K^{\mathfrak{p}_+}$
(see Definition~\ref{def:1.4}).
Let
 $\mathcal{L} := G \times_K U \to G/K$
 be the $G$-equivariant holomorphic vector bundle associated to $\chi$.
Then, there is a natural injective continuous $G$-homomorphism
$\mathcal{H} \to \mathcal{O}(\mathcal{L})$.
\end{lemma}

\begin{proof}
Let $(\ , \ )_{\mathcal{H}}$ be a $G$-invariant inner product on $\mathcal{H}$.
We write $(\ , \ )_U$ for the induced inner product on $U$.
Then, $K$ acts unitarily on $\mathcal{H}$,
and in particular on $U$.
We consider the map
$$
   G \times \mathcal{H} \times U \to \mathbb{C}\, ,
    \ \
   (g,v,u) \mapsto (\pi(g)^{-1} v, u)_{\mathcal{H}} = 
   (v, \pi(g) u)_{\mathcal{H}} \, .
$$
For each fixed $g \in G$ and $v \in \mathcal{H}$,
 the map
 $U \to \mathbb{C}, \ u \mapsto (\pi(g)^{-1} v, u)_{\mathcal{H}}$
 is an anti-linear functional on $U$.
Then there exists a unique element
 $F_v(g) \in U$ by the Riesz representation theorem
 for the finite dimensional Hilbert space $U$ such that
$$
    (F_v(g), u)_U = (\pi(g)^{-1} v, u)_{\mathcal{H}}
    \ \
    \text{ for any } u \in U\, .
$$
Then it is readily  seen that
 $F_v(g k) = \chi(k)^{-1} F_v(g)$
 and
  $F_{\pi(g') v}(g) = F_v({g'}^{-1} g)$
 for any $g, g' \in G$, $k \in K$ and $v \in \mathcal{H}$.
As $u$ is a smooth vector in $\mathcal{H}$,
 $(F_v(g), u)_U = (v, \pi (g) u)_{\mathcal{H}}$
 is a $C^\infty$-function on $G$.
Then
 $F_v(g)$ is a $C^\infty$-function on $G$ with value in $U$
 for each fixed $v \in \mathcal{H}$.
Thus,
 we have a non-zero $G$-intertwining operator given by
$$
    F: {\mathcal{H}} \to C^\infty(G \times_K U)\, ,
    \quad
    v \mapsto F_v \, .
$$
As $U$ is annihilated by $\mathfrak{p}_+$,
 $F_v$ is a holomorphic section of the holomorphic vector bundle
 $G \times_K U \to G/K$,
 that is,
 $F_v \in  \mathcal{O}(G \times_K U)$.
Then, the non-zero map $F: \mathcal{H} \to  \mathcal{O}(G \times_K U)$ is injective
 because $\mathcal{H}$ is irreducible.
Furthermore, $F$ is continuous by the closed graph theorem.
Hence, Lemma~\ref{lem:5.1} is proved.
\qed
\end{proof}

\subsection{Reduction to real rank condition}
\label{subsec:3.6}
The next Lemma is a stepping-stone to Theorem~\ref{thm:A}.
It becomes also a key lemma to the theorem that
 the  action of a  subgroup
$H$ on the bounded symmetric domain $G/K$ is `strongly visible' 
for any symmetric pair $(G,H)$
(see \cite{visiblesymm}).

\begin{lemma}
\label{lem:5.1}
Suppose $\mathfrak{g}$ is a real simple Lie algebra of Hermitian type.
Let $\tau$ be an involutive automorphism of\/ $\mathfrak{g}$, 
 commuting with a fixed Cartan involution $\theta$.
Then there exists an involutive automorphism $\sigma$ of\/ $\mathfrak{g}$
 satisfying the following three conditions:
\begin{eq-text}\label{eqn:3.6.one}
  $\sigma$, $\tau$ and $\theta$ commute with one another.
\end{eq-text}
\begin{eq-text}\label{eqn:3.6.2}
  $\rrank \mathfrak{g}^{\tau \theta} 
 = \rrank \mathfrak{g}^{\sigma, \tau \theta}$.
\end{eq-text}
\begin{eq-text}\label{eqn:3.6.3}
  $\sigma Z = -Z$.
\end{eq-text}
\end{lemma}

\begin{proof}
We shall give a proof in the special case
$\tau = \theta$ in Subsection~\ref{subsec:4.1}.
For the general case,
see   \cite[Lemma~3.1]{visiblesymm}
or \cite[Lemma 5.1]{xkmf}.
\qed
\end{proof}

\subsection{Proof of Theorem~\ref{thm:A}}
\label{subsec:proofthmA}

Now, we are ready to complete the proof of Theorem~\ref{thm:A}.

Without loss of generality,
we may and do assume that $G$ is simply connected.
Let $(\pi,\mathcal{H})$ be an irreducible unitary highest weight
representation of scalar type.
We define a holomorphic line bundle by
$\mathcal{L} := G \times_K \mathcal{H}_K^{\mathfrak{p}_+}$ 
over the Hermitian symmetric space
$D := G/K$.
Then, it follows from Lemma~\ref{lem:3.5} that
there is an injective continuous $G$-intertwining map
$\mathcal{H} \to \mathcal{O}(\mathcal{L})$.

Suppose $(G,H)$ is a symmetric pair.
We first note that for an involutive automorphism $\tau$ of $G$,
there exists $g \in G$ such that
$\tau^g \theta = \theta \tau^g$
if we set
$$
\tau^g(x) := g \tau(g^{-1} x g)g^{-1}
$$
for $x \in G$.
Then, $G^{\tau_g} = gHg^{-1}$
is $\theta$-stable.
Since the multiplicity-free property of the restriction $\pi|_H$ is
unchanged if we replace $H$ by 
$gHg^{-1}$,
we may and do assume that $\theta H = H$,
in other words, $\theta\tau = \tau\theta$.

Now, by  applying Lemma~\ref{lem:5.1},
we can take $\sigma$ satisfying
\eqref{eqn:3.6.one}, \eqref{eqn:3.6.2} and \eqref{eqn:3.6.3}.
We use the same letter $\sigma$ to denote its lift to $G$.
It follows from \eqref{eqn:3.6.3} that
 the induced involutive diffeomorphism $\sigma : G/K \to G/K$
 is anti-holomorphic (see Subsection~\ref{subsec:1.5}).
In light of the conditions \eqref{eqn:3.6.one} and \eqref{eqn:3.6.2}, 
we can apply Lemma~\ref{lem:3.2} to see that for any $x \in D$ there
exists $g \in H$ such that $\sigma(x) = g \cdot x$.

Moreover, by using 
 Lemma~\ref{lem:9.6} in the Appendix,
 we have an isomorphism
 $\overline{\sigma^* \mathcal{L}} \simeq \mathcal{L}$
 as $G$-equivariant holomorphic line bundles over $G/K$.
Therefore,
 all the assumptions of Theorem~\ref{thm:2.2} are satisfied.
Thus, we conclude that the restriction
$\pi|_H$ is multiplicity-free by Theorem~\ref{thm:2.2}.
\qed
\section{Proof of Theorem~\protect\ref{thm:C}}
\label{sec:4}
In this section we give a proof of Theorem~\ref{thm:C}.

Throughout this section,
 we may and do assume that $G$ is simply connected
 so that any automorphism of $\mathfrak{g}$ lifts to $G$.
We divide the proof of Theorem~\ref{thm:C} 
into the following  cases:
\newline
Case I.\enspace
Both $\pi_1$ and $\pi_2$ are highest weight modules.
\newline
Case I$'$.\enspace
Both $\pi_1$ and $\pi_2$ are lowest weight modules.
\newline
Case II.\enspace
 $\pi_1$ is a highest weight module, and $\pi_2$ is a lowest 
weight module.
\newline
Case II$'$.\enspace
$\pi_1$ is a lowest weight module, and $\pi_2$ is a highest weight
module.

\subsection{Reduction to real rank condition}
\label{subsec:4.1}
The following lemma is a special case of Lemma~\ref{lem:5.1} with
$\tau = \theta$.
We shall see that Theorem~\ref{thm:C} in Case~I 
(likewise, Case~I$'$)
reduces to this algebraic result.

\begin{lem} %
\label{lem:4.1.1}
Suppose $\mathfrak{g}$ is a real simple Lie algebra of
  Hermitian type.
Let $\theta$ be a Cartan involution.
Then there exists an involutive automorphism $\sigma$ of $\mathfrak{g}$
 satisfying the following three conditions:
\begin{eq-text}\label{eqn:4.1.1}
  $\sigma$ and $\theta$ commute.
\end{eq-text}
\begin{eq-text}\label{eqn:4.1.2}
  $\rrank \mathfrak{g}= \rrank \mathfrak{g}^\sigma$.
\end{eq-text}
\begin{eq-text}\label{eqn:4.1.3}
  $\sigma Z = -Z$.
\end{eq-text}
\end{lem}

\begin{proof}
We give a proof of the Lemma based on the classification of 
simple Lie algebras
$\mathfrak{g}$ of Hermitian type.

We recall that
 for any involutive automorphism $\sigma$ of $G$,
 there exists $g \in G$
 such that $\sigma^g \theta = \theta \sigma^g$.
Thus,
\eqref{eqn:4.1.1} is always satisfied
 after replacing $\sigma$
 by some $\sigma^g$.  
The remaining conditions \eqref{eqn:4.1.2} and \eqref{eqn:4.1.3} (cf.\ 
Table~\ref{tbl:3.3.2}) are satisfied
 if we choose $\sigma \in \operatorname{Aut}(G)$ in the following 
Table~\ref{tbl:4.1.2}
 for each simple non-compact Lie group $G$ of Hermitian type:
\addtocounter{table}{1}
\begin{table}[H]
\caption{}
\label{tbl:4.1.2}
$$
\vbox{
\offinterlineskip
\def\tablerule{\noalign{\hrule}}
\halign{\strut#&\vrule#&
            \;\;\hfil#\hfil\hfil\;\;&\vrule#&
            \;\;\hfil#\hfil\hfil\;\;&\vrule#&
                \hfil#\hfil\hfil\,&\vrule#&
            \;\;\hfil#\hfil\hfil\;\;&\vrule#\cr\tablerule
&&\multispan7\hfil $(\mathfrak{g}, \mathfrak{g}^\sigma)$ satisfying \eqref{eqn:4.1.2}
 and \eqref{eqn:4.1.3}
 \hfil &\cr\tablerule
&& ${\mathfrak{g}}$ && $\mathfrak{g}^\sigma$
 &&&& $\rrank \mathfrak{g} = \rrank \mathfrak{g}^\sigma$  &\cr\tablerule
&& ${\mathfrak{su}}(p,q)$ && ${\mathfrak{so}}(p,q)$ &&&& $\min(p,q)$ &\cr\tablerule
&& ${\mathfrak{so}}^*(2n)$ && ${\mathfrak{so}}(n,\mathbb{C})$ &&&& $[\frac{1}{2} n]$ &\cr\tablerule
&& ${\mathfrak{sp}}(n,\mathbb{R})$ && ${\mathfrak{gl}}(n,\mathbb{R})$ &&&& $n$ &\cr\tablerule
&& ${\mathfrak{so}}(2,n)$ && ${\mathfrak{so}}(1,n-1)+{\mathfrak{so}}(1,1)$ &&&& $\min(2,n)$ &\cr\tablerule
&& ${\mathfrak{e}}_{6 (-14)}$  && ${\mathfrak{sp}}(2,2)$&&&& $2$ &\cr\tablerule
&& ${\mathfrak{e}}_{7(-25)}$   && $\mathfrak{su}^*(8)$ &&&& $3$&\cr\tablerule
\hfil\cr}}
$$
\end{table}
\noindent
Here,
 we have proved Lemma.
\qed
\end{proof}

\addtocounter{fact}{1}
\begin{rem}%
\label{rem:4.1.3}
The choice of $\sigma$ 
in Lemma~\ref{lem:4.1.1}
 is not unique.
For example,
 we may choose
 $\mathfrak{g}^\sigma \simeq {\mathfrak{e}}_{6(-26)} \oplus \mathbb{R}$
 instead of the above choice
$\mathfrak{g}^\sigma \simeq \mathfrak{su}^*(8)$
 for $\mathfrak{g} = {\mathfrak{e}}_{7(-25)}$.
\end{rem}

\subsection{Proof of Theorem~\ref{thm:C} in Case I}
\label{subsec:3.7}

Let $G$ be a non-compact simply-connected,
simple Lie group such that $G/K$ is a Hermitian symmetric space.

Let $(\pi_1,\mathcal{H}_1)$ and $(\pi_2,\mathcal{H}_2)$ 
 be two irreducible unitary highest weight representations of scalar type.
By Lemma~\ref{lem:3.5}, 
we can realize $(\pi_i,\mathcal{H}_i)$ in 
the space $\mathcal{O}(\mathcal{L}_i)$
of holomorphic sections of 
the holomorphic line bundle
$\mathcal{L}_i := G \times_K (\mathcal{H}_i)_K^{\mathfrak{p}_+}$
$(i = 1,2)$ over the Hermitian symmetric space $G/K$.
We now define a holomorphic line bundle
$\mathcal{L} := \mathcal{L}_1 \boxtimes \mathcal{L}_2$
over $D := G/K \times G/K$ as the outer tensor product of
$\mathcal{L}_1$ and $\mathcal{L}_2$.
Then, we have naturally an injective continuous
$(G \times G)$-intertwining map
$\mathcal{H}_1 \widehat{\otimes} \mathcal{H}_2 \to
 \mathcal{O}(\mathcal{L})$.

Let us take an involution $\sigma'$ of $\mathfrak{g}$ as in
Lemma~\ref{lem:4.1.1} 
(but we use the letter $\sigma'$ instead of $\sigma$),
and lift it to $G$.
We set $\sigma := \sigma' \times \sigma'$.
Then it follows from \eqref{eqn:4.1.3} that $\sigma'$ acts
anti-holomorphically on $G/K$,
and so does $\sigma$ on $D$.
Furthermore, we have isomorphisms of holomorphic line bundles
$\overline{(\sigma')^* \mathcal{L}_i} \simeq \mathcal{L}_i$
$(i = 1,2)$
by Lemma~\ref{lem:9.6} and thus
$\overline{\sigma^* \mathcal{L}} \simeq \mathcal{L}$.

We now introduce another involutive automorphism $\tau$ of $G \times G$
  by $\tau(g_1, g_2) := (g_2, g_1)$.
Then $(G \times G)^\tau = \diag(G) :=\set{(g,g)}{g \in G}$.
We shall use the same letter $\theta$ to denote the Cartan involution 
$\theta \times \theta$ on $G \times G$
(and $\theta \oplus \theta$ on $\mathfrak{g} \oplus \mathfrak{g}$).
Then,
we observe the following isomorphisms:
\begin{alignat*}{2}
&(\mathfrak{g} \oplus \mathfrak{g})^{\tau\theta}
 = \{ (X, \theta X): X \in \mathfrak{g} \}
&&\simeq \mathfrak{g} \, ,
\\
&(\mathfrak{g} \oplus \mathfrak{g})^{\sigma,\tau\theta}
 = \{ (X, \theta X) : X \in \mathfrak{g}^{\sigma'} \}
&&\simeq \mathfrak{g}^{\sigma'} .
\end{alignat*}%
Thus,
 the condition \eqref{eqn:4.1.2}
implies 
$$
\rrank (\mathfrak{g} \oplus 
\mathfrak{g})^{\tau\theta}
 =\rrank (\mathfrak{g} \oplus \mathfrak{g})^{\sigma, \tau\theta}
\, .
$$
Therefore,
 given $(x_1,x_2) \in D \simeq (G \times G)/(K \times K)$,
 there exists $(g,g) \in (G \times G)^\tau$ satisfying
 $(g \cdot x_1, g \cdot x_2) = (\sigma'(x_1), \sigma'(x_2))$
 $(= \sigma(x_1,x_2))$
 by Lemma~\ref{lem:3.2}. %

Let us apply Theorem~\ref{thm:2.2}
 to the setting $(\mathcal{L} \to D,
  \mathcal{H}_1\widehat\otimes\mathcal{H}_2, \diag(G),
 \sigma)$.
Now that all the assumptions of Theorem~\ref{thm:2.2} are satisfied,
  we conclude that
 the tensor product $\pi_1 \widehat\otimes \pi_2$ is
 multiplicity-free as a $G$-module,
that is,
Theorem~\ref{thm:C} holds in the case I. 
\qed
\subsection{Proof of Theorem~\ref{thm:C} in Case II}
\label{subsec:4.2}
Let us give a proof of Theorem~\ref{thm:C} in the case II.
We use the same $\tau$ as in Subsection~\ref{subsec:3.7},
that is, $\tau(g_1,g_2) := (g_2,g_1)$
and define a new involution $\sigma$ by
$\sigma := \tau \theta$,
that is,
$\sigma(g_1,g_2) = (\theta g_2,\theta g_1)$
for $g_1,g_2 \in G$.
Obviously,
 $\sigma$,
 $\tau$ and the Cartan involution $\theta$ of $G \times G$
 all commute.

We write $M$ for the Hermitian symmetric space $G/K$,
 and $\overline M$ for the conjugate complex manifold.
Then $\sigma$ acts anti-holomorphically on 
$D := M \times \overline{M}$
because so does $\tau$ and because $\theta$ acts holomorphically.

By the obvious identity
$(\mathfrak{g} \oplus \mathfrak{g})^{\tau\theta}
 = (\mathfrak{g} \oplus \mathfrak{g})^{\sigma,\tau\theta}$,
we have 
$\rrank (\mathfrak{g} \oplus \mathfrak{g})^{\tau\theta}
 = \rrank (\mathfrak{g} \oplus \mathfrak{g})^{\sigma,\tau\theta}$
$(= \rrank \mathfrak{g})$.
Therefore, it follows from Lemma~\ref{lem:3.2} that
 for any $(x_1, x_2) \in D$
 there exists $(g,g) \in (G \times G)^\tau$ such that
 $\sigma(x_1, x_2) = (g,g) \cdot (x_1, x_2)$.

Suppose $\pi_1$ (respectively, $\pi_2$) is a unitary highest
(respectively, lowest) weight representation of scalar type.
We set 
$\mathcal{L}_1 := G \times_K (\mathcal{H}_1)_K^{\mathfrak{p}_+}$
and
$\mathcal{L}_2 := G \times_K (\mathcal{H}_2)_K^{\mathfrak{p}_-}$.
Then,
$\mathcal{L}_1 \to M$ and
$\mathcal{L}_2 \to \overline{M}$ are both holomorphic line bundles,
and we can realize $\pi_1$ in
$\mathcal{O}(M,\mathcal{L}_1)$,
and $\pi_2$ in
$\mathcal{O}(\overline{M},\mathcal{L}_2)$,
respectively.
Therefore,
 the
 outer tensor product $\pi_1 \boxtimes \pi_2$
 is realized in a subspace of holomorphic sections
 of the holomorphic line bundle 
$\mathcal{L} := \mathcal{L}_1 \boxtimes \mathcal{L}_2$
over $D = M \times \overline{M}$.

Now, we apply Theorem~\ref{thm:2.2} to
$(\mathcal{L} \to D, \mathcal{H}_1 \widehat{\otimes} \mathcal{H}_2,
 \diag(G), \sigma)$. 
The condition \eqref{eqn:2.2.2} holds by Lemma~\ref{lem:9.6}.
Hence, all the assumptions of Theorem~\ref{thm:2.2} are satisfied,
and therefore,
 Theorem~\ref{thm:C} holds in the case II.
\qed

Hence, Theorem~\ref{thm:C} has been proved.

\section{Uniformly bounded multiplicities --- 
Proof of Theorems~\protect\ref{thm:B} and \protect\ref{thm:D}}
\label{sec:5}

This section gives the proof of Theorems~\ref{thm:B} and \ref{thm:D}.
Since the proof of Theorem~\ref{thm:B} parallels to that of
Theorem~\ref{thm:D}, 
we  deal mostly with Theorem~\ref{thm:D} here.
Without loss of generality,
we assume $G$ is a non-compact simple Lie group of Hermitian type.

\subsection{General theory of restriction}
\label{subsec:3.4}
A unitary representation $(\pi,\mathcal{H})$ of a group $L$ is 
\textit{discretely decomposable} if $\pi$ is unitarily equivalent to
the discrete Hilbert sum of irreducible unitary representations of
$L$: 
$$
\pi \simeq 
   \sideset{}{^\oplus}\sum_{\mu \in \widehat{L}}
   m_\pi (\mu) \mu \, .
$$
Furthermore,
we say $\pi$ is $L$-\textit{admissible}
(\cite{xkdecomp}) if all the multiplicities
$m_\pi (\mu)$ are finite.
In this definition,
we do not require $m_\pi(\mu)$ to be uniformly bounded with respect to
$\mu$. 

Suppose $L'$ is a subgroup of $L$. 
Then, the restriction of $\pi$ to $L'$ is regarded as a unitary
representation of $L'$.
If $\pi$ is $L'$-admissible,
then $\pi$ is $L$-admissible
(\cite[Theorem~1.2]{xkdecomp}).  

We start with recalling from \cite{xkdisc}
 a discrete decomposability theorem of branching laws
 in the following settings:

\begin{fct}%
\label{fact:3.4.1}
{\rm 1)}
Suppose $\tau$ is of holomorphic type 
(see Definition~\ref{def:holo-anti})
 and set $H := G_0^\tau$.
If  $\pi$ is an irreducible unitary highest weight representation of $G$,
 then $\pi$ is \adm{(H \cap K)}.
In particular, $\pi$ is $H$-admissible. 
The restriction $\pi |_H$ splits into a discrete Hilbert sum of 
irreducible unitary highest weight representations
 of $H$:
\begin{equation}
    \pi|_{H} \simeq {\sum_{\mu \in \widehat{H}}}^\oplus
 m_\pi(\mu) \mu
   \qquad
   \text{(discrete Hilbert sum)},
\label{eqn:3.4.1}
\end{equation}
 where the multiplicity $m_{\pi}(\mu)$ is finite for every $\mu$.  
\newline
{\rm 2)}
Let $\pi_1, \pi_2$ be two irreducible unitary highest weight 
representations of $G$.
Then the tensor product $\pi_1 \widehat\otimes \pi_2$
 is \adm {K} under the diagonal action.
Furthermore, $\pi_1 \widehat{\otimes} \pi_2$
 splits into a discrete Hilbert sum of irreducible unitary
highest weight representations of
  $G$, each occurring with finite multiplicity.
Furthermore,
 if at least one of $\pi_1$ or $\pi_2$ is a holomorphic discrete series
 representation for $G$,
 then any irreducible summand is 
 a holomorphic discrete series representation.
\end{fct}%
\begin{proof}
See \cite[Theorem 7.4]{xkdisc} 
for the proof.  
The main
 idea of the proof is taking normal derivatives of holomorphic
 sections, 
which 
goes back to S. Martens \cite{xmartens}.  
The same idea was also employed in a number of papers including Lipsman
 (\cite[Theorem~4.2]{xlipad}) 
 and Jakobsen--Vergne (\cite[Corollary~2.3]{xjv}).  
\qed
\end{proof}

\begin{remark}
\label{rem:Fact51}
Fact~\ref{fact:3.4.1} (1) holds more generally for a closed subgroup
$H$ satisfying the following two conditions:

1)\enspace
$H$ is $\theta$-stable.

2)\enspace
The Lie algebra $\mathfrak{h}$ of $H$ contains $Z$.

\noindent
Here,
we recall that $Z$ is the generator of the center of $\mathfrak{k}$.
The proof is essentially the same as that of Fact~\ref{fact:3.4.1} (1).
\end{remark}

Theorem~\ref{thm:B} (2) follows from Theorem~\ref{thm:A} and
Fact~\ref{fact:3.4.1} (1).
Likewise, Theorem~\ref{thm:D} (2) follows from Theorem~\ref{thm:C}
and Fact~\ref{fact:3.4.1} (2).
What remains to show for Theorems~\ref{thm:B} and \ref{thm:D} is 
the uniform boundedness of multiplicities.

\subsection{Remarks on Fact~\ref{fact:3.4.1}}
\label{fact3.4.1}

Some remarks on Fact~\ref{fact:3.4.1} %
are in order.
\begin{rem}
\label{rem:3.4.2}
A Cartan involution $\theta$ is clearly of holomorphic type because 
$\theta Z = Z$.
If $\theta = \tau$ then $H = K$ and any irreducible summand $\mu$ is
finite dimensional.
In this case, the finiteness of $m_\pi (\mu)$ in
 Fact~\ref{fact:3.4.1} %
 (1) is a special case
 of Harish-Chandra's admissibility theorem
 (this holds for
 any irreducible unitary representation $\pi$ of $G$).  
\end{rem}

\begin{rem}
Fact~\ref{fact:3.4.1} %
 asserts in particular
 that there is no continuous spectrum
 in the irreducible decomposition formula.
The crucial assumption for this is that $(G,H)$ is of holomorphic
type. 
In contrast, the restriction $\pi|_H$ is not
discretely decomposable if $(G,H)$ is of anti-holomorphic type and if
$\pi$ is a holomorphic discrete series representation of $G$
 (\cite[Theorem~5.3]{xkdecomp}).
In this setting,
 R. Howe, J. Repka, G. \'{O}lafsson, B. \O rsted, 
van Dijk, S. Hille, M. Pevzner, V. Molchanov,
Y. Neretin, G. Zhang and others
 studied irreducible decompositions 
 of the restriction $\pi|_{H}$
 by means of the $L^2$-harmonic analysis on Riemannian symmetric spaces
$H/H \cap K$
 (\cite{xvd,xvdh,xDijkPev,xhoweseesaw,xmol, xnere, xoo, xoz, xrep}).
The key idea in Howe and Repka
\cite{xhoweseesaw, xrep} is that a holomorphic
function on $G/K$ is uniquely determined by its restriction to the
totally real submanifold $H/H \cap K$
(essentially, the unicity theorem of holomorphic functions),
and that any function on $H/H\cap K$ can be approximated
(in a sense)
by holomorphic functions on $G/K$ 
(essentially, the Weierstrass polynomial approximation theorem).
\end{rem}

\begin{rem}
A finite
 multiplicity theorem of the branching law \eqref{eqn:3.4.1} %
with respect to semisimple symmetric pairs $(G,H)$ 
holds for more general $\pi$ 
(i.e.\  $\pi$ is not  a highest weight module),
 under the assumption that $\pi$ is discretely decomposable
 as an $(\mathfrak{h}_{\mathbb{C}}, H \cap K)$-module  
 (see \cite[Corollary~4.3]{xkdecoass}, \cite{xkbeijing}).
However, 
 the multiplicity of the branching law can be infinite 
if the restriction is not discretely decomposable
 (see Example~\ref{ex:finite infinite}).
\end{rem}

\begin{rem}
 Theorems~\ref{thm:B} and \ref{thm:D} assert that multiplicities
 $m_\pi (\mu)$ in Fact~\ref{fact:3.4.1} are {\bf uniformly bounded}
 when we vary $\mu$.
This is a distinguished feature for the restriction of 
highest weight representations $\pi$.
A similar statement may fail 
if $\pi$ is not a highest weight module 
 (see Example~\ref{exam:6.3}).
\end{rem}

\subsection{Reduction to the scalar type case}
\label{subsec:4.3}
In order to deduce  Theorem~\ref{thm:D}~(1) from Theorem~\ref{thm:D}~(2),
 we use the idea of `coherent family' of representations of reductive
 Lie groups (for example, see \cite{xvg}).
For this,
we prepare the following  Lemma \ref{lem:4.3} and Proposition \ref{lem:4.4}.

\begin{lemma}%
\label{lem:4.3}
Suppose that $(\pi, \mathcal{H})$ is an irreducible
  unitary highest weight representation of $G$.
Then there exist an irreducible unitary highest weight 
representation $\pi'$
 of scalar type
 and a finite dimensional representation $F$ of $G$
 such that the underlying $(\mathfrak{g}_{\mathbb{C}},K)$-module
 $\pi_K$ occurs as a subquotient of the tensor product $\pi_K' \otimes F$.
\end{lemma}

\begin{proof}
Without loss of generality,
we may and do assume that $G$ is simply connected.
Since $G$ is a simple Lie group of Hermitian type,
the center $\mathfrak{c}(\mathfrak{k})$ of $\mathfrak{k}$ is one
dimensional. 
We take its generator $Z$ as in Subsection~\ref{subsec:1.5},
and write $C$ for the connected subgroup with Lie algebra
$\mathfrak{c}(\mathfrak{k})$. 
Then, $K$ is isomorphic to the direct product group of $C$ and a
semisimple group $K'$.

As $(\pi, \mathcal{H})$ is an irreducible unitary highest weight
representation of $G$,
 $\mathcal{H}_K^{\mathfrak{p}_+}$
 is an irreducible (finite dimensional) unitary representation of $K$.
The $K$-module $\mathcal{H}_K^{\mathfrak{p}_+}$
 has an expression $\sigma\otimes \chi_0$,
 where $\sigma \in \widehat{K}$ such that $\sigma|_C$ is trivial
 and $\chi_0$ is a unitary character of $K$.

Let $\chi'$ be a unitary character of $K$ such that
 $\chi'$ is trivial on the center $Z_G$ of $G$
 (namely, $\chi'$ is well-defined as a representation of
 $\operatorname{Ad}_G(K) \simeq K/Z_G$).
For later purposes,
 we take $\chi'$ such that $-\sqrt{-1}\, d \chi'(Z) \gg 0$.
There exists an irreducible finite dimensional representation $F$
 of $G$
 such that $F^{\mathfrak{p}_+} \simeq \sigma \otimes \chi'$
 as $K$-modules
 because $\sigma \otimes \chi'$ is well-defined
 as an algebraic representation of $\Ad_G(K)$.  

We set
 $\chi := \chi_0 \otimes (\chi')^*$ of $K$.  
Because $- \sqrt{-1}\, d \chi (Z) \ll 0$, 
the irreducible highest weight $(\mathfrak{g}_{\mathbb{C}},K)$-module 
$V'$ such that
$(V')^{\mathfrak{p}_+} \simeq \chi$
is unitarizable.
Let $(\pi', \mathcal{H}')$ denote the irreducible unitary
 representation of $G$ whose underlying 
$(\mathfrak{g}_{\mathbb{C}}, K)$-module 
$\mathcal{H}_K'$ is isomorphic to
$V'$.
Since $\mathcal{H}'_K$ is an irreducible
 $(\mathfrak{g}_{\mathbb{C}},K)$-module, 
 $\mathcal{H}_K' \otimes F$ is a $(\mathfrak{g}_{\mathbb{C}},K)$-module
 of finite length.
Furthermore,
as $\mathcal{H}'_K$ is a highest weight module,
so are
 all subquotient modules of $\mathcal{H}'_K \otimes F$.
Then, $\mathcal{H}_K$ arises as a subquotient of $\mathcal{H}_K' \otimes F$
because the $K$-module $\mathcal{H}_K^{\mathfrak{p}_+}$
occurs as a subrepresentation of 
$(\mathcal{H}_K'\otimes F)^{\mathfrak{p}_+}$
 in view of
$$
     \mathcal{H}_K^{\mathfrak{p}_+} 
     \simeq 
     \sigma \otimes \chi_0
     \simeq
     \chi \otimes (\sigma \otimes \chi')
     \simeq
     (\mathcal{H}_K')^{\mathfrak{p}_+} \otimes F^{\mathfrak{p}_+}
     \subset (\mathcal{H}_K' \otimes F)^{\mathfrak{p}_+} \, .  
$$
Hence,
 we have shown Lemma~\ref{lem:4.3}. %
\qed
\end{proof}

\subsection{Uniform estimate of multiplicities for tensor products}
\label{subsec:4.4}
Let $(\pi,X)$ be a \gk-module
 of finite length.  
This means that $\pi$ admits a chain of submodules
\begin{equation}
\label{eqn:YiX}
0 = Y_0 \subset Y_1 \subset \cdots \subset Y_N = X
\end{equation}
such that $Y_i/Y_{i-1}$ is irreducible for $i=1,\dots,N$.
The number $N$ is independent of the choice of the chain \eqref{eqn:YiX},
and we will write
$$
m(\pi) := N \, .
$$
That is, 
 $m(\pi)$ is
 the number of irreducible \gk-modules 
 (counted with multiplicity)
 occurring as subquotients in $\pi$.  
Here is a uniform estimate of $m(\pi)$ under the operation of tensor products:

\begin{prop}%
\label{lem:4.4}
Let $F$ be a finite dimensional representation of
 a real reductive connected Lie group $G$.
Then there exists a constant $C \equiv C(F)$ such that
$$
   m(\pi \otimes F) \le C 
$$
 for any irreducible \gk-module $\pi$. 
\end{prop}

Before entering the proof, we fix some terminologies:

\begin{defn}
\label{def:Grothen}
We write $\mathcal{F}(\mathfrak{g}_{\mathbb{C}},K)$ 
for the category of $(\mathfrak{g}_{\mathbb{C}},K)$-modules of finite
length. 
The \textit{Grothendieck group} 
$\mathcal{V}(\mathfrak{g}_{\mathbb{C}},K)$ of
$\mathcal{F}(\mathfrak{g}_{\mathbb{C}},K)$ is the abelian group
generated by $(\mathfrak{g}_{\mathbb{C}},K)$-modules of finite length,
modulo the equivalence relations
$$
X \sim Y + Z
$$
whenever there is a short exact sequence
$$
0 \to Y \to X \to Z \to 0
$$
of $(\mathfrak{g}_{\mathbb{C}},K)$-modules.
Then
$$
m: \mathcal{F}(\mathfrak{g}_{\mathbb{C}},K) \to \mathbb{N}
$$
induces a group homomorphism of abelian groups:
$$
m: \mathcal{V}(\mathfrak{g}_{\mathbb{C}},K) \to \mathbb{Z} \, .
$$
\end{defn}

The Grothendieck group 
$\mathcal{V}(\mathfrak{g}_{\mathbb{C}},K)$ is isomorphic to the free
abelian group having irreducible
$(\mathfrak{g}_{\mathbb{C}},K)$-modules as its set of finite generators. 

Suppose $(\pi,X)$ is a $(\mathfrak{g}_{\mathbb{C}},K)$-module of
finite length. 
Then, in the Grothendieck group
$\mathcal{V}(\mathfrak{g}_{\mathbb{C}},K)$, 
we have the relation
\begin{equation}
\label{eqn:XmpiY}
X = \bigoplus_Y m_\pi(Y) Y \, ,
\end{equation}
where the sum is taken over irreducible
$(\mathfrak{g}_{\mathbb{C}},K)$-modules. 
Then we have
\begin{equation}
\label{eqn:mpiY}
m(\pi) = \sum_Y m_\pi(Y) \, .
\end{equation}
Suppose $(\pi',X')$ is also a $(\mathfrak{g}_{\mathbb{C}},K)$-modules
of finite length. 
We set
\begin{align}
[\pi:\pi'] :={}
&\dim\Hom_{(\mathfrak{g}_{\mathbb{C}},K)}
(\bigoplus_Y m_\pi(Y)Y, \bigoplus_Y m_{\pi'}(Y)Y)
\label{eqn:numberGro}
\\
={}
&\sum_Y m_\pi(Y) m_{\pi'}(Y) \, .
\label{eqn:numberGro2}
\end{align}
The definition \eqref{eqn:numberGro} makes sense in a more general
setting where one of $X$or $X'$ is not of finite length.
To be more precise, we recall from \cite[Definition 1.1]{xkdecoass}:
\begin{defn}\label{def:infdeco}
Let $\mathcal{A}(\mathfrak{g}_{\mathbb{C}},K)$ be the category of
$(\mathfrak{g}_{\mathbb{C}},K)$-modules $(\pi,X)$ having the following
properties: 
\itm{1)}
($K$-admissibility)
$\dim\Hom_K(\tau,\pi) < \infty$ for any $\tau\in\widehat{K}$.
\itm{2)}
(discretely decomposability, see \cite[Definition 1.1]{xkdecoass})
$X$ admits an increasing filtration
$$
0 = Y_0 \subset Y_1 \subset Y_2 \subset \cdots
$$
of $\mathfrak{g}_{\mathbb{C}}$-modules such that $Y_i/Y_{i-1}$ is of
finite length and that $X = \bigcup_{i=1}^\infty Y_i$.
\end{defn}

We refer the reader to \cite{xkdecoass} for algebraic results on
discretely decomposable 
$(\mathfrak{g}_{\mathbb{C}},K)$-modules such as:
\begin{lem}\label{lemma:infdeco}
Suppose $X \in \mathcal{A}(\mathfrak{g}_{\mathbb{C}},K)$.
\itm{1)}
Any submodule or quotient of $X$ is an object of
$\mathcal{A}(\mathfrak{g}_{\mathbb{C}},K)$. 
\itm{2)}
The tensor product $X\otimes F$ is also an object of 
$\mathcal{A}(\mathfrak{g}_{\mathbb{C}},K)$
for any finite dimensional $(\mathfrak{g}_{\mathbb{C}},K)$-module.
\end{lem}

For $X \in \mathcal{A}(\mathfrak{g}_{\mathbb{C}},K)$,
we can take the filtration $\{Y_i\}$ such that $Y_i/Y_{i-1}$ is
irreducible as a $(\mathfrak{g}_{\mathbb{C}},K)$-module for any $i$.
Then, for any irreducible $(\mathfrak{g}_{\mathbb{C}},K)$-module,
$$
\#\{i: \text{$Y_i/Y_{i-1}$ is isomorphic to $Y$}\}
$$
is finite and independent of the filtration, 
which we will denote by $m_\pi(Y)$.

\begin{defn}\label{def:mfagk}
Suppose $X \in \mathcal{A}(\mathfrak{g}_{\mathbb{C}},K)$.
We say the
 $(\mathfrak{g}_{\mathbb{C}},K)$-module $X$ is \textit{multiplicity-free} if 
$$
m_\pi(Y) \le 1 
\quad\text{for any irreducible $(\mathfrak{g}_{\mathbb{C}},K)$-module
$Y$}.
$$
\end{defn}

This concept coincides with Definition \ref{def:1.2} if $X$ is the
underlying $(\mathfrak{g}_{\mathbb{C}},K)$-module of a unitary
representation of $G$.
The point of Definition \ref{def:mfagk}
 is that we allow the case where $X$ is not
unitarizable.

Generalizing \eqref{eqn:numberGro2},
we set
$$
[\pi:\pi'] := \sum_Y m_\pi(Y) m_{\pi'} (Y)
$$
for $\pi,\pi' \in \mathcal{A}(\mathfrak{g}_{\mathbb{C}},K)$. 
Here are immediate results from the definition:
\begin{lem}\label{lem:pi}
Let $\pi,\pi' \in \mathcal{A}(\mathfrak{g}_{\mathbb{C}},K)$. 
\smallskip
\itm{1)}
$[\pi:\pi'] < \infty$ if at least one of $\pi$ and 
$\pi'$ belongs to $\mathcal{F}(\mathfrak{g}_{\mathbb{C}},K)$.
\smallskip
\itm{2)}
$\dim\Hom_{(\mathfrak{g}_{\mathbb{C}},K)} (\pi,\pi')
 \le [\pi: \pi']$.
\smallskip
\itm{3)}
$[\pi: \pi'] = [\pi': \pi]$.
\smallskip
\itm{4)}
$m_\pi(Y) = [\pi: Y]$\quad if $Y$ is an irreducible
$(\mathfrak{g}_{\mathbb{C}},K)$-module.
\smallskip
\itm{5)}
$[\pi:\pi'] \le m(\pi)$\quad if $\pi'$ is multiplicity-free.
\end{lem}

Now, we return to Proposition \ref{lem:4.4}.

{
\renewcommand{\proofname}{Proof of Proposition~\ref{lem:4.4}} %
\begin{proof}
We divide the proof into three steps:
\newline
$\underline{\text{Step 1}}$
 ($\pi$ is a finite dimensional representation): \enspace
We shall prove 
\begin{equation}
     m(\pi \otimes F) \le \dim F
\label{eqn:4.4.1}
\end{equation}
 for any finite dimensional representation $\pi$ of $G$.

Let $\mathfrak{b} = \mathfrak{t} + \mathfrak{u}$ be a Borel subalgebra of $\mathfrak{g}_\mathbb{C}$
 with $\mathfrak{u}$ nilradical.
We denote by $H^j(\mathfrak{u}, V)$ the $j$th cohomology group
 of the Lie algebra $\mathfrak{u}$
 with coefficients in a $\mathfrak{u}$-module $V$.
Since the Lie algebra $\mathfrak{b}$ is solvable,
 we can choose a $\mathfrak{b}$-stable filtration
$$
    F = F_k \supset F_{k-1} \supset \dots \supset F_0 = \{0\} 
$$
 such that $\dim F_i/F_{i-1} = 1$.

Let us show  by induction on $i$
 that 
\begin{equation}
   \dim H^0(\mathfrak{u}, \pi \otimes F_i) \le i \, .
\label{eqn:4.4.2}
\end{equation}
This will imply $m(\pi\otimes F) =
 \dim H^0(\mathfrak{u},\pi \otimes F)\le k = \dim F$.

The inequality \eqref{eqn:4.4.2} %
is trivial if $i=0$.
Suppose \eqref{eqn:4.4.2} %
holds for $i-1$.
The short exact sequence of $\mathfrak{b}$-modules
$$
   0 \to \pi \otimes F_{i-1} \to \pi \otimes F_i \to \pi \otimes (F_i/F_{i-1})
       \to 0
$$
 gives rise to a long exact sequence
\begin{multline*}
     0 \to 
     H^0(\mathfrak{u}, \pi \otimes F_{i-1})
 \to
     H^0(\mathfrak{u}, \pi \otimes F_i)
 \to H^0(\mathfrak{u}, \pi \otimes (F_i/F_{i-1}))
\\
 \to
     H^1(\mathfrak{u}, \pi \otimes F_{i-1})
  \to \dots
\end{multline*}
 of $\mathfrak{t}$-modules.
In particular,
 we have
\begin{equation}
  \dim   H^0(\mathfrak{u}, \pi \otimes F_i)
   \le  \dim H^0(\mathfrak{u}, \pi \otimes F_{i-1})
 +
 \dim H^0(\mathfrak{u}, \pi \otimes (F_i/F_{i-1})) \, .
\label{eqn:4.4.3}
\end{equation}
Because $F_i/F_{i-1}$ is trivial as a $\mathfrak{u}$-module,
 we have
\begin{equation}
     H^0(\mathfrak{u}, \pi \otimes (F_i/F_{i-1}))
   =
     H^0(\mathfrak{u}, \pi) \otimes (F_i/F_{i-1}) \, .
\label{eqn:4.4.4}
\end{equation}
By definition $H^0(\mathfrak{u}, \pi)$ is 
 the space of highest weight vectors, and therefore
 the dimension of the right-hand side of \eqref{eqn:4.4.4} %
is one.  
Now,
 the inductive assumption combined with \eqref{eqn:4.4.3} %
implies
$\dim   H^0(\mathfrak{u}, \pi \otimes F_i) \le i$,
 as desired.
\newline
$\underline{\text{Step 2}}$
 ($\pi$ is a principal series representation): \enspace
In this step, we consider the case where $\pi$ is a principal series
 representation.
We note that $\pi$ may be reducible here.

Let $P = L N$ be a Levi decomposition of a minimal parabolic subgroup $P$
 of $G$,
 $W$ an irreducible (finite dimensional) representation of $L$,
 and $\operatorname{Ind}_P^G (W)$
 the underlying \gk-module of a principal series representation
 induced from the representation $W \boxtimes 
\mathbf{1} %
$ of $P = L N$
 (without $\rho$-shift).
Then, the socle filtration is unchanged so far as the
 parameter lies in the equisingular set,
and thus, there are only finitely many possibilities of the socle
 filtration of $\operatorname{Ind}_P^G (W)$ for irreducible
 representations $W$ of $L$.
We denote by $m(G)$ the maximum of $m(\operatorname{Ind}_P^G (W))$
for irreducible representations $W$ of $L$.

Let $F$ be a finite dimensional representation of $G$.
Then
 we have an isomorphism of \gk-modules
$$
\operatorname{Ind}_P^G (W) \otimes F \simeq
 \operatorname{Ind}_P^G (W \otimes F) \, ,
$$
 where $F$ is regarded as a $P$-module on the right-hand side.
We take a $P$-stable filtration 
$$
 W_n := W \otimes F \supset W_{n-1} \supset \dots
 \supset W_0 = \{0\}
$$
 such that each $W_{i}/W_{i-1}$ is irreducible as a $P$-module.
We notice that $n \le \dim F$ by applying Step~1 to the $L$-module $F|_L$.
As $\operatorname{Ind}_P^G (W \otimes F)$
 is isomorphic to $\bigoplus_{i=1}^n \operatorname{Ind}_P^G (W_{i}/W_{i-1})$
 in the Grothendieck group $\mathcal{V}(\mathfrak{g}_\mathbb{C}, K)$, 
 we have shown that 
$$
     m(\operatorname{Ind}_P^G (W) \otimes F)
     \le n \; m(G)
     \le (\dim F) \; m(G)
$$
 for any irreducible
 finite dimensional representation $W$ of $L$.
\newline
$\underline{\text{Step 3}}$  (general case): \enspace
By Casselman's subrepresentation theorem (see \cite[Chapter 3]{xwal}),
any irreducible $(\mathfrak{g}_{\mathbb{C}},K)$-module $\pi$ is
realized as a subrepresentation of some induced representation 
$\operatorname{Ind}_P^G (W)$.
Then
$$
m(\pi \otimes F) \le m(\pi \otimes \operatorname{Ind}_P^G (W)) \le C
$$
by step 2.
Thus, Proposition \ref{lem:4.4} is proved.
\qed
\end{proof}
}

\subsection{Proof of Theorem~\ref{thm:D}}
\label{subsec:4.5}
Now let us complete the proof of Theorem~\ref{thm:D}.  

Let $\pi = \pi_1 \boxtimes \pi_2$ be an irreducible
 unitary highest weight representation
 of $G':=G \times G$.  
It follows from Lemma~\ref{lem:4.3} %
 that there exist an irreducible unitary highest weight representation 
 $\pi'= \pi_1' \boxtimes \pi_2'$
 of scalar type
 and a finite dimensional representation $F$ of $G'$
 such that $\pi_K$ occurs as a subquotient of $\pi_K' \otimes F$.  

By using the notation \eqref{eqn:numberGro},
we set
$[V_1:V_2] := [(V_1)_K:(V_2)_K]$
 for $G$-modules $V_1$ and $V_2$ of finite length.  
Then, for $\mu\in\widehat{G}$, we have
\begin{align}
m_{\pi_1,\pi_2}(\mu) 
  &= \dim\Hom_G (\mu, \pi|_{\diag(G)})
\nonumber
\\
  &\leq  \left[\mu: \pi |_{\diag(G)}\right]
\nonumber
\\
    &\leq \left[\mu:(\pi' \otimes F)|_{\diag(G)}\right]
\nonumber
\\
  &= \left[\mu \otimes (F^* |_{\diag (G)}): \pi'|_{\diag(G)}\right]
\nonumber
\\
    &\leq m(\mu \otimes (F^* |_{\diag (G)}))
\label{eqn:4.5.1}%
\\
    &\leq C(F^*) \, .  
\nonumber
\end{align}
Here the inequality \eqref{eqn:4.5.1} %
follows from Lemma \ref{lem:pi} (5) because
 $\pi'|_{\diag(G)}\simeq \pi_1' \widehat \otimes \pi_2'$
is multiplicity-free (see Theorem~\ref{thm:D}~(2)).  
In the last inequality, $C(F^*)$
 is the constant in Proposition~\ref{lem:4.4}. %
This completes the proof of Theorem~\ref{thm:D} (1).  
\qed

\begin{remark}
The argument in Subsections~\ref{subsec:8.3} and 
\ref{subsec:pf tensordeco}
gives a different and more straightforward proof of
Theorem~\ref{thm:D}. 
\end{remark}

\section{Counter examples}
\label{sec:6}

In this section,
 we analyze the assumptions in Theorems~\ref{thm:A}
 and \ref{thm:B} by counterexamples, that is, 
 how the conclusions fail if we relax the assumptions on the
 representation $\pi$.

Let $(G, H)$ be a reductive symmetric pair corresponding to 
 an involutive automorphism $\tau$ of $G$,  
and $\pi$ an irreducible unitary representation of $G$.  
We shall see that the multiplicity of an irreducible summand occurring
 in the restriction $\pi|_H$ can be:
\itm{1)}
\textbf{greater than one}
if $\pi$ is not of scalar type (but we still assume that $\pi$ is a
 highest weight module);
\itm{2)}
\textbf{finite but not uniformly bounded}
if $\pi$ is not a highest weight module
(but we still assume that $\pi|_H$ decomposes discretely);
\itm{3)}
\textbf{infinite}
if $\pi|_H$ contains continuous spectra.

Although our concern in this paper is mainly with
 a non-compact subgroup $H$,
 we can construct such  examples for (1) and (2) even 
for $H = K$ 
(a maximal compact subgroup modulo the center of $G$).

 Case (1) will be discussed in Subsection~\ref{subsec:6.2}, %
  (2) in Subsection~\ref{subsec:6.3}, %
 and (3) in Subsection~\ref{subsec:6.4}, %
respectively.
To construct an example for (3), we use those for (1) and (2).

\subsection{Failure of multiplicity-free property}
\label{subsec:6.2}
Let $G = Sp(2, \mathbb{R})$.
Then, the maximal compact subgroup $K$ is isomorphic to $U(2)$.
We take a compact Cartan subalgebra $\mathfrak{t}$.
Let $\{f_1, f_2\}$ be the standard basis
of $\sqrt{-1}\, \mathfrak{t}^*$
 such that
 $\Delta(\mathfrak{g}, \mathfrak{t}) = 
   \{\pm f_1 \pm f_2, \pm 2 f_1, \pm 2 f_2\}$,
 and
we fix a positive system 
$\Delta^+(\mathfrak{k}, \mathfrak{t}) := \{f_1 - f_2\}$.
In what follows,
 we shall use the notation $(\lambda_1, \lambda_2)$
 to denote the character $\lambda_1 f_1 + \lambda_2 f_2$ of $\mathfrak{t}$.

Given $(p,q) \in \mathbb{Z}^2$ with $p \ge q$,
 we denote by $\F{U(2)}{(p,q)}$ the irreducible representation of 
$U(2)$
 with highest weight $(p,q) = p f_1 + q f_2$.
Then $\dim \F{U(2)}{(p,q)} = p-q+1$.

The set of holomorphic discrete series representations of $G$
 is parametrized by $\lambda := (\lambda_1, \lambda_2) \in \mathbb{N}^2$
 with $\lambda_1 > \lambda_2 > 0$.
We set
$\mu \equiv (\mu_1,\mu_2):=(\lambda_1+1,\lambda_2+2)$
and denote by
$\hwm{G}{\mu} \equiv \hwm{Sp(2,\mathbb{R})}{(\mu_1,\mu_2)}$
 the holomorphic discrete series representation of $G$
 characterized by
\begin{alignat*}{2}
 &\ \text{$Z(\mathfrak{g})$-infinitesimal character }
 = (\lambda_1, \lambda_2)
 &&\quad\text{(Harish-Chandra parameter)},
\\
 &\ \text{minimal $K$-type }
 = \hwm{U(2)}{(\mu_1, \mu_2)}
 &&\quad\text{(Blattner parameter)}.
\end{alignat*}
We note that $\hwm{G}{\mu}$ is of scalar type
 if and only if $\mu_1 = \mu_2$.

We know from Theorem~\ref{thm:B} that multiplicities of $K$-type
$\tau$ occurring in $\hwm{G}{\mu}$ are uniformly bounded for fixed
$\mu = (\mu_1,\mu_2)$.
Here is the formula:

\begin{exa}[upper bound of $K$-multiplicities
 of holomorphic discrete series]%
\label{exam:6.1}
\begin{equation}
   \sup_{\tau \in \widehat{K}} \dim {\Hom}_K(\tau, \hwm{G}{\mu}|_K)
   = \left[\frac{\mu_1 - \mu_2+2}{2}\right] \, .
\label{eqn:6.2.1}
\end{equation}
The right side of \eqref{eqn:6.2.1} %
$=1$
 if and only if
 either of the following two cases holds:
\begin{subequations}
  \renewcommand{\theequation}{\theparentequation) (\alph{equation}}%
\begin{alignat}{2}
 \mu_1 &= \mu_2
 &&\text{(i.e. $\hwm{G}{\mu}$ is of scalar type),}
\label{eqn:6.2.2a}%
\\
 \mu_1 &= \mu_2 + 1 \quad
 &&\text{(i.e. $\hwm{G}{\mu}$ is of two dimensional minimal $K$-type).}  
\label{eqn:6.2.2b}%
\end{alignat}
\end{subequations}
Thus,
 the branching law of the restriction $\hwm{G}{\mu} |_K$ 
is   multiplicity-free  
 if and only if $\mu_1 = \mu_2$
or $\mu_1 = \mu_2 + 1$.
The multiplicity-free property for
$\mu_1 = \mu_2$
(i.e.\ for $\hwm{G}{\mu}$ of scalar type)
follows from Theorem~\ref{thm:A}.
The multiplicity-free property for
$\mu_1 = \mu_2 + 1$
is outside of the scope of this paper, but
can be explained in the general framework of the 
 `vector bundle version' of Theorem~\ref{thm:2.2} 
(see \cite[Theorem~2]{RIMS}, \cite{mfbdle}).
\end{exa}

\begin{proof}
It follows from the Blattner formula
 for a holomorphic discrete series representation 
(\cite{xjohnson}, \cite{xschmidherm}) that
 the $K$-type formula of $\hwm{G}{\mu}$ is given by
\begin{align}
  \hwm{G}{\mu}|_K
 &\simeq
 \F{U(2)}{(\mu_1,\mu_2)}
 \otimes
 S(\mathbb{C}^3)
\nonumber
\\
 &=\F{U(2)}{(\mu_1,\mu_2)}
 \otimes
 \bigoplus \Sb a \ge b \ge 0 \\ (a,b) \in \mathbb{N}^2 \endSb
          \F{U(2)}{(2 a, 2 b)} \, ,
\label{eqn:holoU2}
\end{align}
 where $K=U(2)$ acts on $\mathbb{C}^3 \simeq S^2(\mathbb{C}^2)$ 
  as the symmetric tensor of the natural representation.
We write $n_\mu(p,q)$ for the multiplicity of the $K$-type
$\hwm{U(2)}{(p,q)}$ occurring in 
$\hwm{G}{\mu} \equiv \hwm{Sp(2,\mathbb{R})}{(\mu_1,\mu_2)}$,
that is,
$$
n_\mu(p,q) := \dim \Hom_K
   (\hwm{K}{(p,q)}, \hwm{G}{\mu} |_K) \, .
$$
Then, applying the Clebsch--Gordan formula \eqref{eqn:1.8.1f}
to \eqref{eqn:holoU2}, 
we obtain
$$
n_\mu(p,q) = \# \{ (a,b) \in \mathbb{N}^2: 
   \text{$(a,b)$ satisfies $a \ge b \ge 0$,
         \eqref{eqn:ab1} and \eqref{eqn:ab2}} \},
$$
where
\begin{align}
& \,
p+q = \mu_1 + \mu_2 + 2a + 2b \, ,
\label{eqn:ab1}
\\
&\max (2a+\mu_2, 2b+\mu_1) \le p \le 2a+\mu_1 \, .
\label{eqn:ab2}
\end{align}
In particular,
for fixed $(\mu_1, \mu_2)$ and $(p,q)$,
the integer $b$ is determined by $a$ from \eqref{eqn:ab1},
whereas the integer $a$ satisfies the inequalities
$p-\mu_1 \le 2a \le p-\mu_2$.
Therefore,
$$
n_\mu(p,q) \le \left[ \frac{(p-\mu_2)-(p-\mu_1)}{2} \right] +1
   = \left[ \frac{\mu_1-\mu_2+2}{2} \right] \, .
\qquad\qed
$$
\end{proof}

\subsection{Failure of uniform boundedness}
\label{subsec:6.3}
We continue the setting of Subsection~\ref{subsec:6.2}.
Let $B$ be a Borel subgroup of
$G_{\mathbb{C}} \simeq Sp(2,\mathbb{C})$.
Then, there exist $4$ closed orbits of $K_\mathbb{C} \simeq GL(2,\mathbb{C})$
 on the full flag variety $G_{\mathbb{C}}/B$.
(By the Matsuki duality,
there exist $4$ open orbits of $G = Sp(2,\mathbb{R})$ on
$G_{\mathbb{C}}/B$. 
This observation will be used in the proof of
Example~\ref{ex:finite infinite}.) 
By the Beilinson--Bernstein correspondence,
we see that
 there are $4$ series of discrete series representations of $G$.
Among them,
 two are holomorphic and anti-holomorphic discrete series representations,
that is,
$\hwm{G}{\mu}$
and
$(\hwm{G}{\mu})^*$
(the contragredient representation)
with notation as in Subsection~\ref{subsec:6.2}.
The other two series are non-holomorphic discrete series
representations. 
Let us parametrize them.
For
 $\lambda := (\lambda_1, \lambda_2) \in \mathbb{Z}^2$
 ($\lambda_1 > -\lambda_2 >0$),
we write $W_\lambda$ for the discrete series representation of $G$
 characterized by
\begin{alignat*}{2}
 &\ \text{$Z(\mathfrak{g})$-infinitesimal character }
  = (\lambda_1, \lambda_2)
 &&\quad\text{(Harish-Chandra parameter)},
\\
 &\ \text{minimal $K$-type }
 = \F{U(2)}{(\lambda_1 + 1, \lambda_2)}
 &&\quad\text{(Blattner parameter)}.
\end{alignat*}
Then, non-holomorphic discrete series representations are either
$W_\lambda$ or its contragredient representation $W_\lambda^*$ for
some $\lambda \in \mathbb{Z}^2$ with
$\lambda_1 > -\lambda_2 > 0$.
We define a $\theta$-stable Borel subalgebra
 $\mathfrak{q} = \mathfrak{t}_{\mathbb{C}} + \mathfrak{u}$
of 
$\mathfrak{g}_{\mathbb{C}} = 
  \mathfrak{k}_{\mathbb{C}} + \mathfrak{p}_{\mathbb{C}}$
 such that 
$$
  \Delta(\mathfrak{u} \cap \mathfrak{p}_{\mathbb{C}}, \mathfrak{t})
 := \{2 f_1, f_1  + f_2, -2 f_2\}\, ,
 \quad
  \Delta(\mathfrak{u} \cap \mathfrak{k}_{\mathbb{C}}, \mathfrak{t})
   := \{f_1 - f_2\} \, .
$$
Then,
 the Harish-Chandra module $(W_\lambda)_K$
 is isomorphic to the cohomological parabolic induction
 $\mathcal{R}_{\mathfrak{q}}^1 (\mathbb{C}_{(\lambda_1, \lambda_2)})$
 of degree $1$ as $(\mathfrak{g}_{\mathbb{C}},K)$-modules with
 the notation 
and the normalization as
in \cite{xvr}.
We set $\mu_1 := \lambda_1+1$ and $\mu_2 := \lambda_2$. 

\begin{exa}[multiplicity of $K$-type of non-holomorphic
 discrete series $W_\lambda$]
\label{exam:6.3}
We write $m_\lambda(p,q)$ for the multiplicity of the $K$-type
 $\F{U(2)}{(p,q)}$ occurring in $W_\lambda$,
 that is,
$$
  m_\lambda(p,q) := \dim \Hom_K(\F{U(2)}{(p,q)}, W_\lambda|_K) \, .
$$
Then,
 $m_\lambda(p,q) \neq 0$
 only if $(p,q) \in \mathbb{Z}^2$
 satisfies
\begin{equation}
     p \ge \mu_1\, ,
     \
     p-q \ge \mu_1 - \mu_2
\text{ and }
     p-q \in 2 \mathbb{Z} + \mu_1 + \mu_2 \, .
\label{eqn:6.3.1}
\end{equation}
Then,
\begin{equation}
m_\lambda(p,q)
   = 1 + \min(\left[\frac{p-\mu_1}{2}\right],
             \frac{p-q-\mu_1 +\mu_2}{2}) \, .
\label{eqn:6.3.2}
\end{equation}
In particular, for each fixed $\lambda$,
 the $K$-multiplicity in $W_\lambda$ is not uniformly bounded, namely,
$$
     \sup_{\tau \in \widehat{K}} \dim \Hom_K(\tau, W_\lambda|_K)
   =\sup_{(p,q) \text{ satisfies \eqref{eqn:6.3.1} %
}} m_{\lambda}(p,q)
   = \infty \, .  
$$
\end{exa}

\begin{proof}
For $p,q\in\mathbb{Z}$,
we write $\mathbb{C}_{(p,q)}$ for the one dimensional representation
of $\mathfrak{t}_{\mathbb{C}}$ corresponding to the weight
$pf_1 + qf_2 \in \mathfrak{t}_{\mathbb{C}}^*$.
According to the $\mathfrak{t}_{\mathbb{C}}$-module isomorphism:
$$
\mathfrak{u} \cap \mathfrak{p}_{\mathbb{C}}
 \simeq \mathbb{C}_{(2,0)} \oplus \mathbb{C}_{(1,1)} \oplus
        \mathbb{C}_{(0,-2)} \, ,
$$
the symmetric algebra 
 $S(\mathfrak{u} \cap \mathfrak{p}_{\mathbb{C}})$ is decomposed into 
 irreducible representations of $\mathfrak{t}_{\mathbb{C}}$ as
\begin{align}
   S(\mathfrak{u} \cap \mathfrak{p}_{\mathbb{C}}) 
   &\simeq \bigoplus_{a,b,c \in \mathbb{N}} 
    S^a(\mathbb{C}_{(2,0)})
    \otimes S^b(\mathbb{C}_{(1,1)})
    \otimes S^c(\mathbb{C}_{(0,-2)})
\nonumber
\\
   &\simeq \bigoplus \Sb a, b, c \in \mathbb{N} \endSb 
                 \mathbb{C}_{(2a + b, b-2c)} \, .
\label{eqn:Supabc}
\end{align}
We denote by
 $H^j(\mathfrak{u} \cap \mathfrak{k}_\mathbb{C}, \pi)$
 the $j$th cohomology group of the Lie algebra $\mathfrak{u} \cap \mathfrak{k}_\mathbb{C}$
 with coefficients in the
 $\mathfrak{u} \cap \mathfrak{k}_\mathbb{C}$-module $\pi$.
If $\pi$ is a $\mathfrak{k}_{\mathbb{C}}$-module,
 then $H^j(\mathfrak{u} \cap \mathfrak{k}_\mathbb{C}, \pi)$ 
becomes naturally a $\mathfrak{t}_{\mathbb{C}}$-module.
Then, Kostant's version of the Borel--Weil--Bott theorem (e.g.\ 
\cite[Chapter 3]{xvg})
shows that
\begin{equation}
\label{eqn:ukcohom}
   H^j(\mathfrak{u} \cap \mathfrak{k}_\mathbb{C}, \F{U(2)}{(p,q)})
  = \begin{cases} \mathbb{C}_{(p,q)} & (j=0)\, , \\
           \mathbb{C}_{(q-1,p+1)} & (j=1)\, , \\
           \{0\}                  & (j\neq 0,1)\, .
    \end{cases}
\end{equation}
By using the Blattner formula due to Hecht--Schmid
 (e.g.\ \cite[Theorem~6.3.12]{xvg}),
 the $K$-type formula of $W_\lambda$ is given by
\begin{align*}
m_\lambda(p,q)&=
 \dim \Hom_K(\F{U(2)}{(p,q)}, W_\lambda|_K)
\\
&=
 \sum_{j=0}^1 (-1)^j \dim
 \Hom_{\mathfrak{t}_{\mathbb{C}}}
 (H^j(\mathfrak{u} \cap \mathfrak{k}_\mathbb{C}, \F{U(2)}{(p,q)}),
 S(\mathfrak{u} \cap \mathfrak{p}_\mathbb{C}) \otimes
 \mathbb{C}_{(\mu_1, \mu_2)}) \, .
\\
\intertext{Now, comparing \eqref{eqn:Supabc} with
 the above formula \eqref{eqn:ukcohom} 
 as $\mathfrak{t}_{\mathbb{C}}$-modules, we see}
m_\lambda(p,q)
&=\#\set{(a,b,c) \in \mathbb{N}^3}{p = 2a + b + \mu_1, q = b - 2c + \mu_2}
\\
 &\phantom{={}}
-\#\set{(a,b,c) \in \mathbb{N}^3}{q-1 = 2a + b + \mu_1, p+1 = b - 2c + \mu_2}
\\
&=\#\set{(a,b,c) \in \mathbb{N}^3}{p = 2a + b + \mu_1, q = b - 2c + \mu_2}
\\
&= 1 + \min(\left[\frac{p-\mu_1}{2}\right],
             \frac{p-q-\mu_1 +\mu_2}{2}) \, .
\end{align*}
Thus, the formula \eqref{eqn:6.3.2} %
has been verified.
\qed
\end{proof}

\subsection{Failure of finiteness of multiplicities}
\label{subsec:6.4}%

Multiplicities of the branching laws can be infinite in general even
for reductive symmetric pairs $(G,H)$.
In this subsection, we review from \cite[Example~5.5]{xkaspm} 
a curious example
of the branching law,
in which
 the multiplicity of a discrete summand is non-zero and finite
 and that of
another discrete summand is infinite.
Such a phenomenon happens only when continuous spectra
appear.

\begin{exa}[infinite and finite multiplicities]
\label{ex:finite infinite}
Let $(G_{\mathbb{C}},G)$ be a reductive symmetric pair
$(Sp(2,\mathbb{C}), Sp(2,\mathbb{R}))$.
We note that $(G_{\mathbb{C}},G)$ is locally isomorphic to the
symmetric pair $(SO(5,\mathbb{C}), SO(3,2))$.
We take a Cartan subgroup $H = TA$ of $G_{\mathbb{C}}$.
We note that 
$T \simeq \mathbb{T}^2$ and $A \simeq \mathbb{R}^2$, 
and identify $\widehat{T}$ with $\mathbb{Z}^2$.

Let $\varpi \equiv \varpi_{(a,b)}^{Sp(2,\mathbb{C})}$
be the unitary principal series representation of $G_{\mathbb{C}}$
induced unitarily from the character $\chi$ of a Borel subgroup $B$
containing $H=TA$ such that
$$
\chi |_H \simeq \mathbb{C}_{(a,b)} \boxtimes \mathbf{1} \, .
$$
We assume $a,b \ge 0$ and set
$$
c(\mu_1,\mu_2;a,b) :=
   \# \{ (s,t,u) \in \mathbb{N}^3:
   a = \mu_1+2s+t, \,
   b=\mu_2+t+2u \} \, .
$$
Then, the discrete part of the branching law of the restriction
$\varpi_{(a,b)}^{Sp(2,\mathbb{C})} |_{Sp(2,\mathbb{R})}$
is given by the following spectra:
\begin{equation}
\bigoplus_{\mu_1 \ge \mu_2 \ge 3}
   c(\mu_1,\mu_2;a,b)
   ( \hwm{Sp(2,\mathbb{R})}{(\mu_1,\mu_2)} \oplus
   \left( \hwm{Sp(2,\mathbb{R})}{(\mu_1,\mu_2)} \right)^* )
 \oplus
\sideset{}{^\oplus}\sum_{\lambda_1>-\lambda_2>0}
   \infty (W_\lambda \oplus W^*_\lambda)
\, ,
\label{eqn:discSp}
\end{equation}
with the notation as in Examples~\ref{exam:6.1} and
\ref{exam:6.3}.

The first term of \eqref{eqn:discSp} 
is a finite sum because there are at most finitely many
$(\mu_1,\mu_2)$ such that $c(\mu_1,\mu_2;a,b) \ne 0$
for each fixed $(a,b)$.
For instance, 
 the first term of \eqref{eqn:discSp}
amounts to
$$
\bigoplus_{\substack{3\le\mu_1\le a \\ \mu_1\equiv a\bmod 2}}
   \hwm{Sp(2,\mathbb{R})}{(\mu_1,3)} \oplus
\bigoplus_{\substack{3\le\mu_1\le a \\ \mu_1\equiv a\bmod 2}}
   \left(\hwm{Sp(2,\mathbb{R})}{(\mu_1,3)}\right)^*
\quad\text{(multiplicity-free)}
$$
if $b=3$.

The second term of \eqref{eqn:discSp} is nothing other than the
direct sum of all non-holomorphic discrete series
   representations of $G=Sp(2,\mathbb{R})$
with infinite multiplicities for any $a$ and $b$.
\end{exa}

{\renewcommand{\proofname}{Sketch of Proof}
\begin{proof}
There exist $4$ open $G$-orbits on $G_{\mathbb{C}}/B$,
for which the isotropy subgroups are all isomorphic to
$T \simeq \mathbb{T}^2$. 
By the Mackey theory,
the restriction $\varpi^{G_{\mathbb{C}}}_{(a,b)}|_G$ 
is unitarily equivalent to the direct sum of
the regular representations realized on $L^2$-sections of
$G$-equivariant line bundles 
$G \times_T \mathbb{C}_{(\pm a, \pm b)} \to G/T$.
That is,
$$
\varpi^{G_{\mathbb{C}}}_{(a,b)}|_G 
\simeq \bigoplus_{\varepsilon_1,\varepsilon_2=\pm1}
   L^2 (G/T, \mathbb{C}_{(\varepsilon_1 a,\varepsilon_2 b)} ) \, .
$$
Therefore, an irreducible unitary representation $\sigma$ of $G$
occurs as a discrete spectrum in
$\varpi_{(a,b)}^{G_{\mathbb{C}}} |_G$
if and only if $\sigma$ occurs as a discrete summand in
$L^2(G/T,\mathbb{C}_{(\varepsilon_1 a, \varepsilon_2 b)})$
for some $\varepsilon_1, \varepsilon_2 = \pm 1$.
Further, 
the multiplicity is given by
$$
\dim\Hom_G (\sigma, \varpi_{(a,b)}^{G_{\mathbb{C}}} |_G)
   = \sum_{\varepsilon_1,\varepsilon_2 = \pm1}
     \dim\Hom_{\mathbb{T}^2}
     ( \mathbb{C}_{(\varepsilon_1 a,\varepsilon_2 b)}, 
             \sigma |_{\mathbb{T}^2} )
$$
by the Frobenius reciprocity theorem.

Since $T$ is compact,
$\sigma$ must be a discrete series representation of
$G = Sp(2,\mathbb{R})$ if $\sigma$ occurs in
$L^2(G/T, \mathbb{C}_{(\varepsilon_1 a, \varepsilon_2 b)})$
as a discrete summand.
We divide the computation of multiplicities into the following two
cases: 

Case I.
$\sigma$ is a 
 holomorphic series representation
or its contragredient representation.
Let
$\sigma = \hwm{Sp(2,\mathbb{R})}{\mu}$. 
Combining \eqref{eqn:holoU2} with the weight formulae
\begin{align*}
S (\mathbb{C}^3) |_{\mathbb{T}^2}
& \simeq
  \bigoplus_{s,t,u\in\mathbb{N}}
  S^s (\mathbb{C}_{(2,0)}) \otimes
  S^t (\mathbb{C}_{(1,1)}) \otimes
  S^u (\mathbb{C}_{(0,2)})
 \simeq
  \bigoplus_{s,t,u\in\mathbb{N}}
  \mathbb{C}_{(2s+t, t+2u)} \, ,
\\
\hwm{U(2)}{(\mu_1,\mu_2)} |_{\mathbb{T}^2}
& \simeq \bigoplus_{\substack{p+q=\mu_1+\mu_2 \\
                              \mu_2\le p \le\mu_1}}
  \mathbb{C}_{(p,q)} \, ,
\end{align*}
we have
$$
\dim\Hom_{\mathbb{T}^2}
   (\mathbb{C}_{(a,b)}, \hwm{Sp(2,\mathbb{R})}{\mu} |_{\mathbb{T}^2})
   = c(\mu_1,\mu_2;a,b)\, .
$$

Case II.
$\sigma$ is a non-holomorphic discrete series representation.
Let $\sigma = W_\lambda$.
It follows from the $K$-type formula
\eqref{eqn:6.3.2} of $W_\lambda$ that we have
$$
\dim\Hom_{\mathbb{T}^2}
   (\mathbb{C}_{(a,b)}, W_{\lambda} |_{\mathbb{T}^2} )
   = \sum_{p\ge q} m_\lambda(p,q)\dim\Hom_{\mathbb{T}^2}
       (\mathbb{C}_{(a,b)},\hwm{U(2)}{(p,q)})
   = \infty \, .
$$
Likewise for $\sigma = W_\lambda^*$
(the contragredient representation).

Hence, the discrete part of the branching law is given by
\eqref{eqn:discSp}.
\qed
\end{proof}
}

\section{Finite Dimensional Cases
--- Proof of Theorems~\protect\ref{thm:E} and \protect\ref{thm:F}}
\label{sec:7}
\subsection{Infinite v.s.\ finite dimensional representations}
\label{subsec:7.1}
Our method applied to infinite dimensional representations in
 Sections~\ref{sec:3} and \ref{sec:4} also applies to
 {\bf finite} dimensional representations,
 leading us to  multiplicity-free theorems,
as stated in Theorems~\ref{thm:E} and \ref{thm:F} in
 Section~\ref{sec:1}, 
 for
the restriction with respect to symmetric pairs.

The comparison with multiplicity-free theorems in the 
infinite dimensional case
 is illustrated by the following correspondence:
\begin{alignat*}{3}
   &\text{a non-compact simple group $G$}
   &&\quad \leftrightarrow
   &&\quad \text{a compact simple group $G_U$}
\\
   &\text{a unitary highest weight module}
   &&\quad \leftrightarrow
   &&\quad \text{a finite dimensional module}
\\
   &\text{{scalar type (Definition~\ref{def:1.4})}}
   &&\quad \leftrightarrow
   &&\quad \text{\lq\lq pan type\rq\rq\ 
     (Definition~\ref{def:7.3.2})} %
\\
   &\text{Theorems~\ref{thm:A} and \ref{thm:B}}
   &&\quad \leftrightarrow
   &&\quad \text{Theorems~\ref{thm:E} and \ref{thm:F}}.
\end{alignat*}
The main goal of this section is to give a proof of
Theorems~\ref{thm:E} and \ref{thm:F}
by using Theorem~\ref{thm:2.2}.
Geometrically, our proof is built
 on the fact 
that the $H_U$ action on the Hermitian
symmetric space is strongly visible if $(G_U,H_U)$ is a symmetric
pair (see \cite{visiblesymm}).

\subsection{Representations associated to maximal parabolic subalgebras}
\label{subsec:7.2}
Let $\mathfrak{g}_\mathbb{C}$ be a complex simple Lie algebra.
We take a Cartan subalgebra $\mathfrak{j}$ of $\mathfrak{g}_\mathbb{C}$,
 and fix a positive system $\Delta^+(\mathfrak{g}_\mathbb{C}, \mathfrak{j})$.
We denote by $\{\alpha_1, \dots, \alpha_n\}$ the set of simple roots,
 and by $\{\omega_1, \dots, \omega_n\} \ (\subset \mathfrak{j}^*)$
 the set of the fundamental weights.

We denote by $\F{\mathfrak{g}_\mathbb{C}}{\lambda}$ 
  irreducible finite dimensional representation of 
$\mathfrak{g}_\mathbb{C}$
 with highest weight $\lambda = \sum_{i=1}^n m_i \omega_i$
 for $m_1, \dots, m_n \in \mathbb{N}$.
It is also regarded as a holomorphic representation of $G_\mathbb{C}$,
 a simply connected complex Lie group with Lie algebra $\mathfrak{g}_\mathbb{C}$.

We fix a simple root $\alpha_i$, 
 and define a maximal parabolic subalgebra
$$
   {\mathfrak{p}}^{-}_{i\mathbb{C}} := 
  {\mathfrak{l}}_{i\mathbb{C}} + {\mathfrak{n}}^{-}_{i\mathbb{C}}
$$
such that the nilradical
 ${\mathfrak{n}}^{-}_{i\mathbb{C}}$ and  the Levi part 
${\mathfrak{l}}_{i\mathbb{C}}\ (\ \supset \mathfrak{j})$
 are given by
\begin{align*}
   \Delta({\mathfrak{l}}_{i\mathbb{C}}, \mathfrak{j})
 &= \mathbb{Z}\text{-span of } 
 \{\alpha_1, \dots, \overset{\wedge}{\alpha_i}, \dots, \alpha_n\}
 \cap \Delta(\mathfrak{g}_\mathbb{C}, \mathfrak{j}) \, ,
\\
   \Delta({\mathfrak{n}}^{-}_{i\mathbb{C}}, \mathfrak{j})
 &= \Delta^-(\mathfrak{g}_\mathbb{C}, \mathfrak{j}) \setminus 
   \Delta({\mathfrak{l}}_{i\mathbb{C}}, \mathfrak{j}) \, .
\end{align*}
We shall see that irreducible finite dimensional representations
realized on generalized flag varieties $G_{\mathbb{C}}/P_{\mathbb{C}}$
is multiplicity-free with respect to any symmetric pairs
if $P_{\mathbb{C}}$ has an abelian unipotent radical.

We write ${P}^{-}_{i\mathbb{C}} = {L}_{i\mathbb{C}} {N}^{-}_{i\mathbb{C}}$
 for the corresponding maximal parabolic subgroup of $G_\mathbb{C}$.

Let
$\Hom(\mathfrak{p}^{-}_{i\mathbb{C}}, \mathbb{C})$
be the set of Lie algebra homomorphisms over $\mathbb{C}$.
Since any such homomorphism vanishes on the derived ideal
$[\mathfrak{p}^{-}_{i\mathbb{C}}, \mathfrak{p}^{-}_{i\mathbb{C}}]$,
 $\Hom({\mathfrak{p}}^{-}_{i\mathbb{C}}, \mathbb{C})$
 is naturally identified with
$$
\Hom({\mathfrak{p}}^{-}_{i\mathbb{C}}/[{\mathfrak{p}}^{-}_{i\mathbb{C}}, 
{\mathfrak{p}}^{-}_{i\mathbb{C}}], \mathbb{C}) \simeq \mathbb{C} \omega_i
\, .
$$

Next, let 
$\Hom(P^{-}_{i\mathbb{C}}, \mathbb{C}^\times)$
be the set of complex Lie group homomorphisms.
Then, we can regard  
 $\Hom({P}^{-}_{i\mathbb{C}}, \mathbb{C}^\times) \subset
\Hom({\mathfrak{p}}^{-}_{i\mathbb{C}}, \mathbb{C})$. 
As its subset,
 $\Hom({P}^{-}_{i\mathbb{C}}, \mathbb{C}^\times)$ 
 is identified with 
 $\mathbb{Z} \omega_i$
 since $G_\mathbb{C}$ is simply connected.

For $k \in \mathbb{Z}$,
 we write $\mathbb{C}_{k \omega_i}$ for the corresponding character of
 ${P}^{-}_{i\mathbb{C}}$,
 and denote by
\begin{equation}
 \mathcal{L}_{k \omega_i}
 := G_{\mathbb{C}} \times_{{P}^{-}_{i\mathbb{C}}} \mathbb{C}_{k \omega_i} 
 \to G_\mathbb{C}/{P}^{-}_{i\mathbb{C}}
\label{eqn:7.2.1}
\end{equation}
  the associated holomorphic line bundle.
We naturally have a representation of $G_\mathbb{C}$
 on the space of holomorphic sections
 $\mathcal{O}\left(\mathcal{L}_{k \omega_i}\right)$.
Then, by the Borel--Weil theory, 
$\mathcal{O}(\mathcal{L}_{k\omega_i})$ is non-zero and irreducible if 
$k \ge 0$ and
we have
 an isomorphism of representations of $G_\mathbb{C}$
 (also of $\mathfrak{g}_\mathbb{C}$):
\begin{equation}
\F{\mathfrak{g}_\mathbb{C}}{k \omega_i}
 \simeq \mathcal{O}\left(\mathcal{L}_{k \omega_i}\right).
\label{eqn:7.2.2}
\end{equation}

\subsection{Parabolic subalgebra with abelian nilradical}
\label{subsec:7.3}
A parabolic subalgebra with abelian nilradical
 is automatically a maximal parabolic subalgebra.
Conversely, the nilradical of a maximal parabolic subalgebra is not
 necessarily abelian.
We recall 
 from Richardson--R\"ohrle--Steinberg \cite{xrrs}
the following equivalent characterization
of such parabolic algebras:
\begin{lem}
\label{lem:7.3.1}
Retain the setting of Subsection~\ref{subsec:7.2}. %
Then, the following four conditions on the pair
 $(\mathfrak{g}_\mathbb{C}, \alpha_i)$ are equivalent:
\newline{\rm{i)}}\enspace
The nilradical ${\mathfrak{n}}^{-}_{i\mathbb{C}}$ is abelian.
\newline{\rm{ii)}}\enspace
$(\mathfrak{g}_\mathbb{C}, {\mathfrak{l}}_{i\mathbb{C}})$ is a symmetric pair.
\newline{\rm{iii)}}\enspace
The simple root $\alpha_i$ occurs in the highest root
 with coefficient one.
\newline{\rm{iv)}}\enspace
 $(\mathfrak{g}_\mathbb{C}, \alpha_i)$ is in the following list
 if we label the simple roots $ \alpha_1, \dots, \alpha_n$
 in the Dynkin diagram as in Table~\ref{tbl:7.3.1}.

\medskip
\begin{eq-text}
{Type}\ $A_n$ \qquad   $\alpha_1, \alpha_2, \dots, \alpha_n$ 
\label{eqn:7.3.1}
\end{eq-text}
\begin{eq-text}
{Type}\ $B_n$ \qquad   $\alpha_1$ 
\label{eqn:7.3.2}
\end{eq-text}
\begin{eq-text}
{Type}\ $C_n$ \qquad   $\alpha_n$
\label{eqn:7.3.3}
\end{eq-text}
\begin{eq-text}
{Type}\ $D_n$ \qquad   $\alpha_1, \alpha_{n-1}, \alpha_n$
\label{eqn:7.3.4}
\end{eq-text}
\begin{eq-text}
{Type}\ $E_6$ \qquad   $\alpha_1, \alpha_{6}$
\label{eqn:7.3.5}
\end{eq-text}
\begin{eq-text}
{Type}\ $E_7$ \qquad   $\alpha_7$
\label{eqn:7.3.6}
\end{eq-text}

\medskip
\noindent
For types $G_2$, $F_4$, $E_8$,
 there are no maximal parabolic subalgebras with abelian nilradicals.
\end{lem}

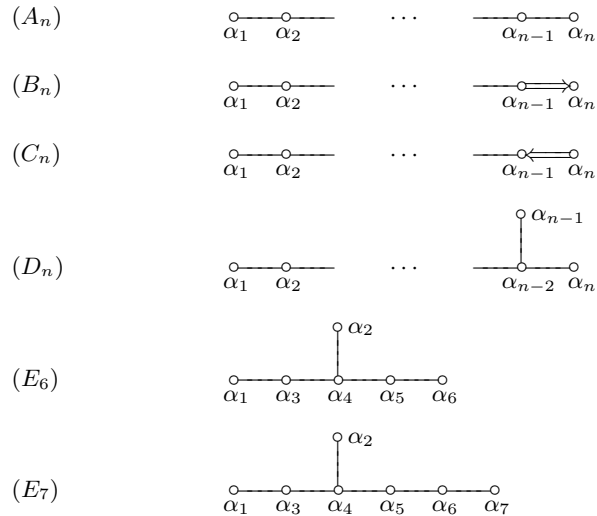
\addtocounter{table}{1}
\begin{table}[H]
\caption{}
\label{tbl:7.3.1}
\begin{align*}
&(A_n)
&&\circdown{\alpha_1} \tsume\yokobo\tsume
 \circdown{\alpha_2} \tsume\yokobo 
 \qquad\cdots\qquad \yokobo\tsume
 \circdown{\alpha_{n-1}} \tsume\yokobo\tsume
 \circdown{\ \alpha_n}
\\[\medskipamount]
&(B_n)
&&\circdown{\alpha_1} \tsume\yokobo\tsume
 \circdown{\alpha_2} \tsume\yokobo 
 \qquad\cdots\qquad \yokobo\tsume
 \circdown{\alpha_{n-1}} \tsume\Longlongrightarrow\tsume
 \circdown{\ \alpha_n}
\\[\medskipamount]
&(C_n)
&&\circdown{\alpha_1} \tsume\yokobo\tsume
 \circdown{\alpha_2} \tsume\yokobo 
 \qquad\cdots\qquad \yokobo\tsume
 \circdown{\alpha_{n-1}} \tsume\Longlongleftarrow\tsume
 \circdown{\ \alpha_n}
\\[\medskipamount]
&(D_n)
&&\circdown{\alpha_1} \tsume\yokobo\tsume
 \circdown{\alpha_2} \tsume\yokobo 
 \qquad\cdots\qquad \yokobo\tsume
 \rlap{\kern0.07em\tatebo}
 \rlap{\kern0.13em\raisebox{2.2em}{$\circ\, \alpha_{n-1}$}}
 \circdown{\alpha_{n-2}} \tsume\yokobo\tsume
 \circdown{\ \alpha_n}
\\[\medskipamount]
&(E_6)
&&\circdown{\alpha_1} \tsume\yokobo\tsume
 \circdown{\alpha_3} \tsume\yokobo\tsume
 \rlap{\kern0.07em\tatebo}
 \rlap{\kern0.13em\raisebox{2.2em}{$\circ\, \alpha_2$}}
 \circdown{\alpha_4} \tsume\yokobo\tsume
 \circdown{\alpha_5} \tsume\yokobo\tsume
 \circdown{\alpha_6}
\\[\medskipamount]
&(E_7)
&&\circdown{\alpha_1} \tsume\yokobo\tsume
 \circdown{\alpha_3} \tsume\yokobo\tsume
 \rlap{\kern0.07em\tatebo}
 \rlap{\kern0.13em\raisebox{2.2em}{$\circ\, \alpha_2$}}
 \circdown{\alpha_4} \tsume\yokobo\tsume
 \circdown{\alpha_5} \tsume\yokobo\tsume
 \circdown{\alpha_6}\tsume\yokobo\tsume
 \circdown{\alpha_7}
\end{align*}
\end{table}

\begin{proof}
See \cite{xrrs}\ for the equivalence
 (i) $\Leftrightarrow$ (iii) $\Leftrightarrow$ (iv).
The implication (iv) $\Rightarrow$ (ii) is straightforward.
For the convenience of the reader,
we present a table of the symmetric pairs 
$(\mathfrak{g}_{\mathbb{C}}, \mathfrak{l}_{i\mathbb{C}})$
corresponding to the index $i$ in (iv).

\medbreak

\begin{center}
\renewcommand{\arraystretch}{1.5}
\begin{tabular}{cccc}
   Type 
   & $\mathfrak{g}_{\mathbb{C}}$
   & $\mathfrak{l}_{i\mathbb{C}}$
   & $i$
\\
   $A_n$
   & $\mathfrak{sl}(n+1,\mathbb{C})$
   & \quad$\mathfrak{sl}(i,\mathbb{C})
      + \mathfrak{sl}(n+1-i,\mathbb{C}) + \mathbb{C}$ \qquad
   & $i=1,2,\dots,n$
\\
   $B_n$
   & $\mathfrak{so}(2n+1,\mathbb{C})$
   & $\mathfrak{so}(2n-1,\mathbb{C})+\mathbb{C}$
   & $ i=1 $
\\
   $C_n$
   & $\mathfrak{sp}(n,\mathbb{C})$
   & $\mathfrak{gl}(n,\mathbb{C})$
   & $ i=n $
\\
  $D_n$
   & $\mathfrak{so}(2n,\mathbb{C})$
   & $\mathfrak{so}(2n-2,\mathbb{C})
      + \mathbb{C} $
   & $ i=1 $
\\
   & $ \mathfrak{so}(2n,\mathbb{C})$
   & $\mathfrak{gl}(n,\mathbb{C}) $
   & $ i=n-1,n $
\\
   $E_6 $
   & $\mathfrak{e}_6$
   & $\mathfrak{so}(10,\mathbb{C}) + \mathbb{C}$
   & $ i=1,6 $
\\
   $E_7$
   & $\mathfrak{e}_7$
   & $\mathfrak{e}_6+\mathbb{C}$
   & $ i=1 $
\end{tabular}
\end{center}
If $(\mathfrak{g}_\mathbb{C}, {\mathfrak{l}}_{i\mathbb{C}})$ is a symmetric pair,
 then $[{\mathfrak{n}}^{-}_{i\mathbb{C}}, {\mathfrak{n}}^{-}_{i\mathbb{C}}]
 \subset {\mathfrak{n}}^{-}_{i\mathbb{C}}  \cap {\mathfrak{l}}_{i\mathbb{C}} = \{0\}$,
 whence (ii) $\Rightarrow$ (i).
\qed
\end{proof}

\addtocounter{fact}{1}
\begin{defn}
\label{def:7.3.2}
We say the representation $\F{\mathfrak{g}_\mathbb{C}}{k \omega_i}$ 
$(k=0, 1,2,\dots)$
 is of\/ \textit{pan type}, or a pan representation if
 $(\mathfrak{g}_\mathbb{C}, \alpha_i)$ satisfies one of
 (therefore, all of) the equivalent  conditions of
 Lemma~\ref{lem:7.3.1}. %
Here, {\bf pan} stands for a {\bf p}arabolic subalgebra
 with {\bf a}belian {\bf n}ilradical.
\end{defn}

\subsection{Examples of pan representations}
\label{subsec:7.4}
\begin{exa}
\label{ex:pangl}
Let $\mathfrak{g}_{\mathbb{C}} = \mathfrak{gl}(n,\mathbb{C})$
and $\lambda = (\lambda_1,\dots,\lambda_n) \in \mathbb{Z}^n$
with $\lambda_1 \ge \lambda_2 \ge \cdots \ge \lambda_n$.
(This $\mathfrak{g}_{\mathbb{C}}$ is not a simple Lie algebra,
but the above concept is defined similarly.) \ 
Then, $\pi_{\lambda}$ is of pan type if and only if
$$
\lambda_1 = \cdots = \lambda_i \ge \lambda_{i+1} = \cdots = \lambda_n
$$
for some $i$ $(1 \le i \le n-1)$.
Then, 
$(\mathfrak{l})_{i\mathbb{C}} \simeq
  \mathfrak{gl}(i,\mathbb{C}) + \mathfrak{gl}(n-i,\mathbb{C})$.

In particular,
the $k$th symmetric tensor representations
$S^k (\mathbb{C}^n)$ $(k \in \mathbb{N})$
and the $k$th exterior representations 
$\Lambda^k (\mathbb{C}^n)$ $(0 \le k \le n)$
are examples of pan representations since their highest weights are given
by $(k, 0, \dots, 0)$ and
$(\underbrace{1,\dots,1}_{k}, \, \underbrace{0,\dots,0}_{n-k} )$,
respectively.
\end{exa}

S. Okada \cite{xokada} studied branching laws for a specific class of
irreducible finite dimensional representations of classical Lie
algebras, 
which he referred to as ``rectangular-shaped representations''. 
The notion of ``pan representations'' is equivalent to that of
rectangular-shaped representations for type 
$(A_n)$, $(B_n)$, and $(C_n)$.
For type $(D_n)$,
$\pi_{k\omega_{n-1}}, \pi_{k\omega_n}$
$(k \in \mathbb{N})$
are  rectangular-shaped representations, 
while $\pi_{k\omega_1}$ $(k \in \mathbb{N})$
are not.

\subsection{Reduction to rank condition}
\label{subsec:7.5}
Suppose $(\mathfrak{g}_{\mathbb{C}}, \alpha_i)$ satisfies the
equivalent conditions in  Lemma~\ref{lem:7.3.1}. %
Let $\theta$ 
 be the complex involutive automorphism
 of the Lie algebra $\mathfrak{g}_\mathbb{C}$
that defines the symmetric pair 
$(\mathfrak{g}_{\mathbb{C}},\mathfrak{l}_{i\mathbb{C}})$.
We use the same letter $\theta$  to denote the corresponding
 holomorphic involution of a simply connected $G_\mathbb{C}$.
We take a maximal compact subgroup $G_U$ of $G_\mathbb{C}$
 such that $\theta G_U = G_U$.
Then $K := G_U^\theta = G_U \cap L_{i \mathbb{C}}$ 
becomes a maximal compact
 subgroup of ${L}_{i\mathbb{C}}$.

Let $\tau$ be
another complex involutive automorphism
of $\mathfrak{g}_{\mathbb{C}}$,
and
$(\mathfrak{g}_{\mathbb{C}},\mathfrak{h}_{\mathbb{C}})$
 the symmetric pair defined by $\tau$.
We also use the same letter $\tau$ to denote its lift to $G_{\mathbb{C}}$. 
We recall from Subsection~\ref{subsec:proofthmA} the `twisted' involution
$\tau^g$ for $g \in G_{\mathbb{C}}$ is given by 
$$
   \tau^g(x) = g \tau(g^{-1} x g) g^{-1} \quad (x \in G_\mathbb{C}) \,
   .
$$
\begin{lemma}
\label{lem:7.5}
Let $(\theta, \tau)$ be as above.
\newline
{\rm 1)}\enspace
There exist an involutive automorphism $\sigma$ of $G_U$
 and $g \in G_\mathbb{C}$
 satisfying the following three conditions
(by an abuse of notation,
we write $\tau$ for $\tau^{g}$):
\begin{eq-text}
 $\tau \mathfrak{g}_U = \mathfrak{g}_U$,
 $\sigma \theta = \theta \sigma$, $\sigma \tau = \tau \sigma$.
\label{eqn:7.5.1}
\end{eq-text}
\begin{eq-text}
The induced action of $\sigma$ on $G_U/K$ is anti-holomorphic.
\label{eqn:7.5.2}
\end{eq-text}
\begin{eq-text}
$(\mathfrak{g}_U)^{\sigma, -\tau, -\theta}$
 contains a maximal abelian subspace in  $(\mathfrak{g}_U)^{-\tau, -\theta}$.
\label{eqn:7.5.3}
\end{eq-text}
{\rm 2)}\enspace
For any $x \in G_U/K$,
there exists $h \in (G_U^\tau)_0$
such that $\sigma(x) = h \cdot x$.
In particular,
each $(G_U^\tau)_0$-orbit on $G_U/K$ is preserved by $\sigma$.
\end{lemma}

\begin{proof}
1) 
See \cite[Lemma~4.1]{visiblesymm} for the proof.
\newline
2) 
The second statement follows from the first statement and
  a similar argument of Lemma~\ref{lem:3.2}.
\end{proof}

\subsection{Proof of Theorem~\ref{thm:E}}
\label{subsec:7.7}
We are now ready to complete the proof of Theorem~\ref{thm:E} 
in Section~\ref{sec:1}.

Let $\pi = \F{\mathfrak{g}_\mathbb{C}}{k \omega_i}$ 
be a representation of pan type.
As in Subsection~\ref{subsec:7.2},
we consider the holomorphic line bundle 
$\mathcal{L}_{k\omega_i} \to G_{\mathbb{C}}/P_{i\mathbb{C}}^-$
and realize
  $\pi$ on the space of holomorphic sections
 $\mathcal{O}(\mathcal{L}_{k \omega_i})$. %
We fix a $G_U$-invariant inner product on 
$\mathcal{O}(\mathcal{L}_{k\omega_i})$.
With notation as in Subsection~\ref{subsec:7.5},
we have a diffeomorphism
$$
G_U/K \simeq G_{\mathbb{C}}/P_{i\mathbb{C}}^- \, ,
$$
through which the holomorphic line bundle 
$\mathcal{L}_{k\omega_i} \to G_{\mathbb{C}}/P_{i\mathbb{C}}^-$
is naturally identified with the $G_U$-equivariant holomorphic line
bundle 
$\mathcal{L} \to D$,
where we set
$\mathcal{L} := G_U \times_K \mathbb{C}_{k\omega_i}$
and
$D := G_U/K$
(a compact Hermitian symmetric space).

Now, applying Lemma~\ref{lem:7.5},
we take $\sigma$ and set
$H := (G_U^\tau)_0$.
We note that the complexification of the Lie algebra of $H$ 
is equal to $\mathfrak{h}_{\mathbb{C}}$ up to a conjugation by
$G_{\mathbb{C}}$. 
By Lemma~\ref{lem:7.5},
the condition \eqref{eqn:2.2.3} in Theorem~\ref{thm:2.2} is
satisfied. 
Furthermore, we see the condition \eqref{eqn:2.2.2} holds by a similar
argument of Lemma~\ref{lem:9.6}.
Therefore, the restriction $\pi|_{(G_U)^\tau_0}$ is multiplicity-free
by Theorem~\ref{thm:2.2}.
Hence,
 Theorem~\ref{thm:E} holds by Weyl's unitary trick. 
\qed

\subsection{Proof of Theorem~\ref{thm:F}}
\label{subsec:7.8}
Suppose $\pi_1$ and $\pi_2$ are representations of pan type.
We realize $\pi_1$ and $\pi_2$ on the space of
 holomorphic sections of holomorphic line bundles over
compact symmetric spaces
 $G_U/K_1$ and $G_U/K_2$,
 respectively.
We write $\theta_i$ for the corresponding involutive
 automorphisms of $G_U$ that define $K_i$ $(i=1,2)$.
In light of Lemma~\ref{lem:7.3.1}~(iv),  %
 we can assume that $\theta_1 \theta_2 = \theta_2 \theta_1$.
Then,
 applying Lemma~\ref{lem:7.5} %
to $(\theta_1, \theta_2)$
 we find an involution $\sigma' \in \Aut(G_U)$
 satisfying the following three conditions:
\begin{eq-text}
 $\sigma' \theta_i = \theta_i \sigma'$ $(i =1, 2)$.
\label{eqn:7.8.1}
\end{eq-text}
\begin{eq-text}
The induced action of $\sigma'$ on $G_U/K_i$ $(i=1,2)$ is anti-holomorphic.
\label{eqn:7.8.2}
\end{eq-text}
\begin{eq-text}
 $(\mathfrak{g}_U)^{\sigma', -\theta_1, -\theta_2}$
 contains a maximal abelian subspace of $(\mathfrak{g}_U)^{-\theta_1, -\theta_2}$.
\label{eqn:7.8.3}
\end{eq-text}
We remark that the condition \eqref{eqn:7.8.2} %
for $i=2$ is not included in Lemma~\ref{lem:7.5}, %
 but follows automatically by our choice of $\sigma$.

We define three involutive automorphisms $\tau$, $\theta$ and $\sigma$ on
 $G_U \times G_U$ by 
$\tau(g_1, g_2) := (g_2, g_1)$,
 $\theta:=(\theta_1, \theta_2)$
 and $\sigma := (\sigma', \sigma')$, respectively.
Then $(G_U \times G_U)^\tau = \diag(G_U)$.
By using the identification
$$
   (\mathfrak{g}_U \oplus \mathfrak{g}_U)^{-\tau}
   = \set{(X, -X)}{X \in \mathfrak{g}_U} 
   \overset{\sim}{\to} \mathfrak{g}_U \, ,
  \quad (X, -X) \mapsto X \, ,
$$
we have isomorphisms
\begin{align*}
   (\mathfrak{g}_U \oplus \mathfrak{g}_U)^{-\tau, -\theta}
   &\simeq (\mathfrak{g}_U)^{-\theta_1, -\theta_2} \, ,
\\
   (\mathfrak{g}_U \oplus \mathfrak{g}_U)^{\sigma, -\tau, -\theta}
   &\simeq (\mathfrak{g}_U)^{\sigma', -\theta_1, -\theta_2} \, .
\end{align*}
Thus, the condition \eqref{eqn:7.8.3} %
implies that 
  $(\mathfrak{g}_U \oplus \mathfrak{g}_U)^{\sigma, -\tau, -\theta}$
 contains a maximal abelian subspace of
  $(\mathfrak{g}_U \oplus \mathfrak{g}_U)^{-\tau, -\theta}$.
Then, by Lemma~\ref{lem:7.5} %
and by a similar argument of Lemma~\ref{lem:3.2} again,
 for any $(x, y) \in G_U/K_1 \times G_U/K_2$ there exists
 a $g \in G_U$ such that
 $\sigma'(x) = g\cdot x$ and $\sigma'(y) = g \cdot y$ simultaneously.
Now, Theorem~\ref{thm:F} %
follows readily from Theorem~\ref{thm:2.2}.
\qed

\subsection{List of multiplicity-free restrictions}
\label{subsec:7.9}
For the convenience of the reader,
 we present the list of
 the triple $(\mathfrak{g}_\mathbb{C}, \mathfrak{h}_\mathbb{C}, i)$ 
for which we can conclude from Theorem~\ref{thm:E} that
 the irreducible finite dimensional representation
 $\hwm{\mathfrak{g}_\mathbb{C}}{k \omega_i}$ of 
 a simple Lie algebra $\mathfrak{g}_\mathbb{C}$
 is multiplicity-free when restricted to $\mathfrak{h}_\mathbb{C}$
 for any $k \in \mathbb{N}$ by Theorem~\ref{thm:E}.
\begin{table}[H]
\caption{}
\label{tbl:7.9.1}
$$
\vbox{
\offinterlineskip
\def\tablerule{\noalign{\hrule}}
\halign{\strut#&\vrule#&
            \;\;\hfil#\hfil\hfil\,&\vrule#&
            \;\;\hfil#\hfil\hfil\,&\vrule#&
              \hfil#\hfil\hfil\,&\vrule#&
            \;\;\hfil#\hfil\hfil\,&\vrule#\cr\tablerule
 && ${\mathfrak{g}_\mathbb{C}}$ && $\mathfrak{h}_\mathbb{C}$
 &&&& $i$  
 &\cr\tablerule
    && $\mathfrak{sl}(n+1,\mathbb{C})$
    && $\mathfrak{sl}(p,\mathbb{C}) + \mathfrak{sl}(n+1-p, \mathbb{C}) + \mathbb{C}$ 
    &&&& $1, 2, \dots, n$ 
&\cr\tablerule
&& $\mathfrak{sl}(n+1,\mathbb{C})$ && $\mathfrak{so}(n+1,\mathbb{C})$
 &&&& $1, 2, \dots, n$ &\cr\tablerule
&& $\mathfrak{sl}(2m,\mathbb{C})$ && $\mathfrak{sp}(m,\mathbb{C})$
 &&&& $1, 2, \dots, 2m-1$ &\cr\tablerule
&& $\mathfrak{so}(2n+1,\mathbb{C})$ && $\mathfrak{so}(p,\mathbb{C})+\mathfrak{so}(2n+1-p, \mathbb{C})$
 &&&& $1$ &\cr\tablerule
&& $\mathfrak {sp}(n,\mathbb{C})$ && $\mathfrak{sp}(p,\mathbb{C}) + \mathfrak{sp}(n-p, \mathbb{C})$
 &&&& $n$ &\cr\tablerule
&& $\mathfrak {sp}(n,\mathbb{C})$ && $\mathfrak{gl}(n,\mathbb{C})$
 &&&& $n$ &\cr\tablerule
&& $\mathfrak {so}(2n,\mathbb{C})$ && $\mathfrak{so}(p,\mathbb{C}) + \mathfrak{so}(2n-p,\mathbb{C})$
 &&&& $1, n-1, n$ &\cr\tablerule
&& $\mathfrak {so}(2n,\mathbb{C})$ && $\mathfrak{gl}(n,\mathbb{C})$
 &&&& $1, n-1, n$ &\cr\tablerule
&& $\mathfrak {e}_6$ && $\mathfrak{so}(10, \mathbb{C})+\mathfrak{so}(2, \mathbb{C})$
 &&&& $1, 6$ &\cr\tablerule
&& $\mathfrak {e}_6$ && $\mathfrak{sl}(6, \mathbb{C})+\mathfrak{sl}(2, \mathbb{C})$
 &&&& $1, 6$ &\cr\tablerule
&& $\mathfrak {e}_6$ && $\mathfrak{f}_4$
 &&&& $1, 6$ &\cr\tablerule
&& $\mathfrak {e}_6$ && $\mathfrak{sp}(4, \mathbb{C})$
 &&&& $1, 6$ &\cr\tablerule
&& $\mathfrak {e}_7$ && $\mathfrak{e}_6+\mathfrak{so}(2, \mathbb{C})$
 &&&& $7$ &\cr\tablerule
&& $\mathfrak {e}_7$ && $\mathfrak{so}(12, \mathbb{C})+\mathfrak{sl}(2, \mathbb{C})$
 &&&& $7$ &\cr\tablerule
&& $\mathfrak {e}_7$ && $\mathfrak{sl}(8, \mathbb{C})$
 &&&& $7$ &\cr\tablerule
\noalign{\smallskip} 
\cr}}
$$
\end{table}

Some of the above cases were previously known to be multiplicity-free
 by case-by-case argument, in particular,
for the case 
$\operatorname{rank} \mathfrak{g}_{\mathbb{C}} = 
 \operatorname{rank} \mathfrak{h}_{\mathbb{C}}$. 
Among them, 
 the corresponding explicit branching
  laws have been studied 
   by S. Okada
  \cite{xokada} and H. Alikawa \cite{xalikawa}.

There are some few representations $\pi$ 
 that are not of pan
 type, 
 but are multiplicity-free
  when restricted to symmetric subgroups $H$.
Our method still works to capture
 such cases, but we do not go into details here 
(see \cite{xkleiden, xkgencar, xksovisible}).
\section{Generalization of the Hua--Kostant--Schmid Formula}
\label{sec:8}

This section discusses an explicit irreducible
decomposition formula of the restriction $\pi|_H$ where the triple
$(\pi,G,H)$ 
satisfies the following two conditions:

1)\enspace
 $\pi$ is a holomorphic discrete series representation
of scalar type (Definition~\ref{def:1.4}).

2)\enspace 
$(G,H)$ is a symmetric pair defined by an involution $\tau$ of
 holomorphic type
(Definition~\ref{def:holo-anti}).

We know a priori from Theorem~\ref{thm:B} (1) that
the branching law is discrete and multiplicity-free.
The main result of this section is Theorem~\ref{thm:gHKS},
which enriches this abstract property with an explicit multiplicity-free 
formula.
The formula for the special case
$H = K$ corresponds to the Hua--Kostant--Schmid formula
(\cite{xhua, xjohnson, xschmidherm}).
We also present explicit formulas for the irreducible decomposition of
the tensor product representation (Theorem \ref{thm:tensordeco}) and
of the restriction 
$U(p,q) \downarrow U(p-1,q)$
(Theorem \ref{thm:upqupq}).

Let us give a few comments on our proof of Theorem~\ref{thm:gHKS}.
Algebraically, our key machinery is Lemma~\ref{lem:HnG} which assures that the irreducible
$G$-decomposition is determined only by its $K$-structure.
Geometrically, a well-known method of taking normal derivatives
 (e.g.\ S. Martens \cite{xmartens}, Jakobsen--Vergne \cite{xjv}) 
 gives a general algorithm to obtain branching laws
for highest weight modules.
This algorithm yields explicit formulae
 by using the observation
 that
 the fiber of the normal bundle for $G^\tau/K^\tau \subset G/K$
 is the tangent space of
 another Hermitian symmetric space $G^{\tau \theta}/K^\tau$.
The key ingredient of the geometry here
 is the following nice properties
  of the two symmetric pairs
  $(G, G^\tau)$ and $(G, G^{\tau\theta})$:

a) $K \cap G^\tau = K \cap G^{\tau\theta}$,

b) $\mathfrak{p} = (\mathfrak{p} \cap \mathfrak{g}^\tau)
     \oplus (\mathfrak{p} \cap \mathfrak{g}^{\tau\theta})$.

Unless otherwise mentioned, 
 we shall assume $H$ is connected,
 that is, 
$H=G_0^{\tau}$ throughout this section.  

\subsection{Notation for highest weight modules}
\label{subsec:8.1}
We set up the notation and give a parametrization of irreducible
 highest weight modules 
for both finite and infinite dimensional cases.

First, we consider finite dimensional representations. 
 Let us take a Cartan subalgebra $\mathfrak{t}$ of 
a reductive Lie algebra
$\mathfrak{k}$
 and fix a positive system $\rt^+(\mathfrak{k}, \mathfrak{t})$.
We denote by $\hwm{\mathfrak{k}}{\mu}$ the irreducible finite dimensional
 representation of $\mathfrak{k}$ with highest weight $\mu$,
 if $\mu$ is a dominant integral weight.
A $\mathfrak{k}$-module $\hwm{\mathfrak{k}}{\mu}$ will be
 written also as $\hwm{K}{\mu}$
 if the action lifts to $K$.

Next, let
 $G$ be a connected reductive Lie group,
 $\theta$ a Cartan involution, 
$K = \{ g \in G: \theta g = g \}$, 
 ${\mathfrak {g}}={\mathfrak {k}}+ {\mathfrak {p}}$
 the corresponding Cartan decomposition
and $\mathfrak{g}_{\mathbb{C}} = \mathfrak{k}_{\mathbb{C}} 
      + \mathfrak{p}_{\mathbb{C}}$
its complexification.
We assume that
 there exists a central element $Z$ of ${\mathfrak {k}}$
 such that 
\begin{equation}
\label{eqn:gkpp}
     {\mathfrak {g}}_{\mathbb{C}}= {\mathfrak {k}}_{\mathbb{C}}+{\mathfrak {p}}_+ + {\mathfrak {p}}_-
\end{equation}
is the eigenspace decomposition
 of $\frac 1{\sqrt{-1}} \operatorname{ad} (Z)$ with eigenvalues 0, 1, and $-1$, 
respectively.
This assumption is satisfied if and only if
 $G$ is locally isomorphic to a direct product
of connected compact Lie groups and non-compact Lie groups of
Hermitian type 
(if $G$ is compact, we can simply take $Z=0$).

We set
\begin{equation}
\label{eqn:sqrtZ}
\widetilde{Z} := \frac{1}{\sqrt{-1}} Z \, .
\end{equation}

As in Definition~\ref{def:1.4},
we say an irreducible $(\mathfrak{g}_{\mathbb{C}},K)$-module 
$V$ is a \textit{highest weight module} if
$$
V^{\mathfrak{p}_+} = \{ v \in V:
   Yv = 0 \quad\mbox{for all $Y \in \mathfrak{p}_+$} \}
$$
is non-zero.
Then, $V^{\mathfrak{p}_+}$ is irreducible as a $K$-module,
and the $(\mathfrak{g}_{\mathbb{C}},K)$-module $V$ is determined
uniquely by the $K$-structure on $V^{\mathfrak{p}_+}$.
If $\mu$ is the highest weight of $V^{\mathfrak{p_+}}$,
we write $V$ as $\hwm{\mathfrak{g}}{\mu}$.
That is, the irreducible $(\mathfrak{g}_{\mathbb{C}},K)$-module
$\hwm{\mathfrak{g}}{\mu}$ is characterized by the $K$-isomorphism:
\begin{equation}
\label{eqn:pigpik}
(\hwm{\mathfrak{g}}{\mu})^{\mathfrak{p}_+}
 \simeq \hwm{\mathfrak{k}}{\mu} \, .
\end{equation}
An irreducible unitary highest weight representation $\pi$ of $G$ will
be denoted by $\hwm{G}{\mu}$ if the underlying 
$(\mathfrak{g}_{\mathbb{C}},K)$-module of $\pi$ is isomorphic to
$\hwm{\mathfrak{g}}{\mu}$.
Let $\Lambda_G$ be the totality of $\mu$ such that 
$\hwm{\mathfrak{g}}{\mu}$ lifts to an irreducible unitary representation of $G$.
For simply connected $G$,
irreducible unitary highest weight representations were classified,
that is,
the set
$\Lambda_G$ $(\subset \sqrt{-1} \mathfrak{t}^*)$
was explicitly found in \cite{xhew} and \cite{xjak}
(see also \cite{xej}).
In particular,
we recall from \cite{xhew} that
$$
\lambda(\widetilde{Z}) \in \mathbb{R}
\quad\text{for any $\lambda \in \Lambda_G$}
$$
and
\begin{equation}
\label{eqn:ubdlmd}
c_G := \sup_{\lambda\in\Lambda_G} \lambda(\widetilde{Z})
< \infty
\end{equation}
if $G$ is semisimple.

The highest weight module
$\hwm{\mathfrak{g}}{\mu}$ is the unique quotient of the generalized
Verma module
\begin{equation}
\label{eqn:Verma}
N^{\mathfrak{g}} (\mu) := U(\mathfrak{g}_{\mathbb{C}})
   \otimes_{U(\mathfrak{k}_{\mathbb{C}} + \mathfrak{p}_+)}
   \hwm{\mathfrak{k}}{\mu} \, ,
\end{equation}
where $\hwm{\mathfrak{k}}{\mu}$ is regarded as a module of the maximal
parabolic subalgebra $\mathfrak{k}_{\mathbb{C}} + \mathfrak{p}_+$
by making $\mathfrak{p}_+$ act trivially.
Furthermore,
 $\hwm{\mathfrak{g}}{\mu}$ has a
 $Z(\mathfrak{g}_\mathbb{C})$-infinitesimal character
 $\mu + \rho_\mathfrak{g} \in \mathfrak{t}^*_{\mathbb{C}}$ 
via the Harish-Chandra isomorphism
$$\Hom_{\text{$\mathbb{C}$-algebra}} 
 (Z(\mathfrak{g}_{\mathbb{C}}),\mathbb{C})
 \simeq \mathfrak{t}_{\mathbb{C}}^* / W \, ,
$$
where $Z(\mathfrak{g}_{\mathbb{C}})$ is the center of the
 enveloping algebra $U(\mathfrak{g}_{\mathbb{C}})$,
$W$ is the Weyl group of the root system
$\Delta(\mathfrak{g},\mathfrak{t})$,
 and
$\rho_{\mathfrak{g}}$ is 
half the sum of positive roots
 $\rt^+(\mathfrak{g}, \mathfrak{t}) 
 := \rt^+(\mathfrak{k}, \mathfrak{t}) \cup \rt(\mathfrak{p}_+, 
     \mathfrak{t})$.

\subsection{Strongly orthogonal roots}
\label{subsec:8.2}
Let $G$ be  a non-compact simple Lie group of Hermitian type,
and $\tau$ an involution 
 of holomorphic type which commutes with the Cartan involution 
$\theta$.

We take a Cartan subalgebra $\mathfrak{t}^\tau$ of
 the reductive Lie algebra
$$
\mathfrak{k}^\tau := \set{X \in \mathfrak{k}}{\tau X = X}
$$
 and extend it to a Cartan subalgebra $\mathfrak{t}$ of $\mathfrak{k}$.
We note that $\mathfrak{t}^\tau = \mathfrak{k}^\tau \cap \mathfrak{t}$.
The pair $(\mathfrak{k}, \mathfrak{k}^\tau)$
forms a reductive symmetric pair,
and $\mathfrak{t}$ plays an analogous role to the fundamental Cartan
subalgebra with respect to this symmetric pair.
Thus, using
 the same argument as 
 in \cite{xvalg}, 
we see that
if $\alpha \in \Delta(\mathfrak{k},\mathfrak{t})$ satisfies
$\alpha|_{\mathfrak{t}^\tau} = 0$
then $\alpha=0$.
Thus,
we can take positive systems $\rt^+(\mathfrak{k}, \mathfrak{t})$
 and $\rt^+(\mathfrak{k}^\tau, \mathfrak{t}^\tau)$
in a compatible way such that
\begin{equation}
\alpha|_{\mathfrak{t}^\tau}
\in \rt^+(\mathfrak{k}^\tau, \mathfrak{t}^\tau)
\quad\text{if $\alpha \in \rt^+(\mathfrak{k}, \mathfrak{t})$} \, .
\label{eqn:8.2.1}
\end{equation}

Since $\tau$ is of holomorphic type,
we have $\tau Z = Z$,
and therefore $\tau \mathfrak{p}_+ = \mathfrak{p}_+$.
Hence, we have a direct sum decomposition 
 $\mathfrak{p}_+ = \mathfrak{p}_+^\tau \oplus \mathfrak{p}_+^{-\tau}$,
 where we set 
$$
   \mathfrak{p}_+^{\pm\tau} := \set{X \in \mathfrak{p}_+}{\tau X = \pm X}
   \, .
$$

Let us consider the reductive subalgebra $\mathfrak{g}^{\tau\theta}$.
Its Cartan decomposition is given by
$$
\mathfrak{g}^{\tau\theta} =
   (\mathfrak{g}^{\tau\theta} \cap \mathfrak{g}^\theta)
     + (\mathfrak{g}^{\tau\theta} \cap \mathfrak{g}^{-\theta})
  = \mathfrak{k}^\tau + \mathfrak{p}^{-\tau} \, ,
$$
and its complexification is given by
\begin{equation}
\label{eqn:8.2.2}
\mathfrak{g}_{\mathbb{C}}^{\tau\theta}
  = \mathfrak{k}_{\mathbb{C}}^\tau \oplus
    \mathfrak{p}_+^{-\tau} \oplus \mathfrak{p}_-^{-\tau} \, . 
\end{equation}
The Cartan subalgebra $\mathfrak{t}^\tau$ of $\mathfrak{k}^\tau$ is
also a Cartan subalgebra of $\mathfrak{g}^{\tau\theta}$.

Let 
$
\rt(\mathfrak{p}_+^{-\tau}, \mathfrak{t}^\tau)
$
 be the set of weights of
 $\mathfrak{p}_+^{-\tau}$ with respect to $\mathfrak{t}^\tau$.
The roots $\alpha$ and $\beta$ are said to be {\it strongly orthogonal}
 if neither $\alpha + \beta$ nor $\alpha - \beta$ is a root.
We take a maximal set of strongly orthogonal roots
 $\{ \nu_1, \nu_2, \dots, \nu_l \}$
 in $\rt(\mathfrak{p}_+^{-\tau}, \mathfrak{t}^\tau)$
 such that
\newline\indent
i)\
 $\nu_1$ is the lowest root among the elements in 
$
\rt(\mathfrak{p}_+^{-\tau}, \mathfrak{t}^\tau),
$
\newline\indent
ii)\
 $\nu_{j+1}$ is the lowest root among the elements in 
$
\rt(\mathfrak{p}_+^{-\tau}, \mathfrak{t}^\tau)
$
that are strongly orthogonal to $\nu_1, \dots, \nu_j$.

A special case applied to $\tau = \theta$ shows
$\mathfrak{k}^\tau = \mathfrak{k}$,
$ \mathfrak{t}^\tau = \mathfrak{t}$,
$ \mathfrak{p}^{-\tau} = \mathfrak{p}$,
and
$\Delta(\mathfrak{p}_+^{-\tau},\mathfrak{t}^\tau) 
 = \Delta(\mathfrak{p}_+,\mathfrak{t})$.
In this case, 
we shall use the notation
$\{ {\bar{\nu}}_1, {\bar{\nu}}_2,\dots,{\bar{\nu}}_{\bar{l}} \}$ 
for a maximal set of strongly orthogonal roots
in $\Delta (\mathfrak{p}_+, \mathfrak{t})$ such that
\begin{enumerate}
    \renewcommand{\theenumi}{\roman{enumi}}
    \renewcommand{\labelenumi}{\theenumi)}
\item  %
${\bar{\nu}}_1$ is the lowest root among
$\Delta(\mathfrak{p}_+,\mathfrak{t})$,
\item  %
${\bar{\nu}}_{j+1}$ is the lowest root among the elements in
$\Delta(\mathfrak{p}_+,\mathfrak{t})$
that are strongly orthogonal to
${\bar{\nu}}_1,\dots,{\bar{\nu}}_j$
$(1 \le j \le {\bar{l}})$.
\end{enumerate}
Then, $\bar{l} = \rrank \mathfrak{g}$
by
 \cite{xkwolf}.
Likewise, in light of \eqref{eqn:8.2.2} for
  the Hermitian symmetric space
$G^{\tau\theta} / G^{\tau\theta} \cap K = 
   G^{\tau\theta} / G^{\tau,\theta}$,
 we have $l = \rrank \mathfrak{g}^{\tau \theta}$.
In general, $l \le \bar{l}$.

\subsection{Branching laws for semisimple symmetric pairs}
\label{subsec:8.5}
It follows from \eqref{eqn:pigpik} that
the highest weight module
 $\hwm{\mathfrak{g}}{\mu}$ is
 of scalar type, namely,
$(\hwm{\mathfrak{g}}{\mu})^{\mathfrak{p}_+}$ is one dimensional, 
if and only if 
\begin{equation}
 \langle \mu, \alpha \rangle = 0 
 \qquad
\text{ for any }
  \alpha \in \rt(\mathfrak{k}, \mathfrak{t}) \, .
\label{eqn:8.5.1}
\end{equation}
Furthermore,
the representation
$\hwm{G}{\mu}$ is 
 a (relative) holomorphic discrete series representation
 of $G$ if and only if
\begin{equation}
 \langle \mu + \rho_\mathfrak{g}, \alpha \rangle < 0 
 \qquad
\text{ for any }
  \alpha \in \rt(\mathfrak{p}_+, \mathfrak{t}) \, .
\label{eqn:8.5.2}
\end{equation}

We are now ready to state
 the branching law of holomorphic discrete series
representations
$\hwm{G}{\mu}$
of scalar type with respect to semisimple symmetric pairs $(G,H)$: 

\begin{theorem}
\label{thm:gHKS}
Let $G$ be a non-compact simple Lie group of Hermitian type.
Assume that $\mu \in \sqrt{-1}\, \mathfrak{t}^*$ satisfies 
\eqref{eqn:8.5.1} %
and \eqref{eqn:8.5.2}. %
Let $\tau$ be an involutive automorphism of $G$ of holomorphic type, 
$H=G_0^\tau$
(the identity component of $G^\tau$),
and $\{\nu_1,\dots,\nu_l \}$
be the set of strongly orthogonal roots in
$\Delta(\mathfrak{p}_+^{-\tau}, \mathfrak{t}^\tau)$
as 
in Subsection~\ref{subsec:8.2}.
Then, $\hwm{G}{\mu}$ decomposes discretely into a multiplicity-free
sum of irreducible 
 $H$-modules: 
\begin{equation}
  \hwm{G}{\mu} |_H
 \simeq
\sideset{}{^\oplus}\sum_{
 \Sb a_1 \ge \dots \ge a_l \ge 0\\ a_1, \dots, a_l \in \mathbb{N} \endSb}
       \hwm{H}{\mu|_{\mathfrak{t}^\tau}- \sum_{j=1}^l a_j \nu_j}
\quad
\text{(discrete Hilbert sum)}.
\label{eqn:8.5.3}
\end{equation}
\end{theorem}

The formula for the case  $H=K$ (that is, $\tau=\theta$)
was previously found by
L.-K. Hua (implicit in the classical case), B. Kostant (unpublished) and 
 W. Schmid \cite{xschmidherm}
(see also Johnson \cite{xjohnson} for an algebraic proof).
In this case, 
each summand in the right side is finite dimensional.

For $\tau \ne \theta$,
some special cases have been also studied by H. Jakobsen,
M. Vergne, J. Xie, W. Bertram and J. Hilgert
\cite{xbehi, xjak, xjv, xxie}.
Further, quantitative results by means of reproducing kernels are
obtained in \cite{xsaid}.
The formula \eqref{eqn:8.5.3} in the above
 generality was first given by the
author \cite{xkmfjp}.

We shall give a proof of Theorem~\ref{thm:gHKS} in
Subsection~\ref{subsec:8.3}.

\subsection{Irreducible decomposition of tensor products}
\label{subsec:8.6}

As we saw in Example~\ref{ex:gpmfd},
the pair $(G \times G, \diag(G))$ forms a symmetric pair. 
Correspondingly, the
 tensor product representation can be regarded as a special (and
easy) case of restrictions of representations with respect to 
symmetric pairs.
This subsection provides a decomposition formula
of the tensor product of two holomorphic discrete series
representations of scalar type.
This is regarded as a counterpart of 
Theorem~\ref{thm:gHKS}
for tensor
product representations. 

We recall from Subsection~\ref{subsec:8.2} that
$\{\bar{\nu}_1,\dots,\bar{\nu}_{\bar{l}}\}$
is a maximal set of strongly orthogonal roots in
$\Delta(\mathfrak{p}_+,\mathfrak{t})$
and $\bar{l}=\rrank\mathfrak{g}$.

\begin{theorem}
\label{thm:tensordeco}
Let $G$ be a non-compact simple Lie group of Hermitian type. 
Assume that $\mu_1, \mu_2 \in \sqrt{-1}\, \mathfrak{t}^*$
satisfy the conditions \eqref{eqn:8.5.1} and \eqref{eqn:8.5.2}.
Then, the tensor product representation
$\hwm{G}{\mu_1} \widehat{\otimes} \hwm{G}{\mu_2}$
decomposes discretely into a multiplicity-free sum of
irreducible $G$-modules: 
$$
\hwm{G}{\mu_1} \widehat{\otimes} \hwm{G}{\mu_2}
\simeq \sum_{\substack{a_1 \ge\dots\ge a_{\bar{l}} \ge 0 \\
                       a_1,\dots,a_{\bar{l}} \,\in\, \mathbb{N}}}
       \hwm{G}{\mu_1 + \mu_2 - \sum_{j=1}^{\bar{l}} a_j {\bar{\nu}}_j}.
$$
\end{theorem}
The proof of Theorem~\ref{thm:tensordeco} will be given in
Subsection~\ref{subsec:pf tensordeco}.

\subsection{Eigenvalues of the central element $Z$}
\label{subsec:Spec}

Our  proof of Theorems~\ref{thm:gHKS} and
\ref{thm:tensordeco} depends on the algebraic lemma 
 that the $K$-type formula
determines the irreducible decomposition of the whole group 
(see Lemma~\ref{lem:HnG}).
This is a very strong assertion, 
which
 fails in general for non-highest weight modules.
This subsection collects some nice properties peculiar to  
highest weight modules that will be used in
the proof of
Lemma~\ref{lem:HnG}.

For a $K$-module $V$, 
we define a subset of $\mathbb C$ by  
$$
     \operatorname{Spec}_{\widetilde{Z}}(V):=\{\text{eigenvalues of
     $\widetilde Z$ on $V$}\} \, , 
$$
where we set 
$$
     \widetilde Z:=\frac {1}{\sqrt{-1}} Z \, .
$$  
For instance, 
$\operatorname{Spec}_{\widetilde{Z}}(V)$ 
is a singleton if $V$ is an irreducible $K$-module.
We also note that
$\operatorname{Spec}_{\widetilde Z}
  ({\mathfrak {g}}_{\mathbb{C}})=\{0, \pm 1\}$
by 
 \eqref{eqn:gkpp}.

\begin{lemma}
\label{lem:spec}
Suppose $V$ is an irreducible \gk-module.  
Then,
\newline
{\rm{1)}}\enspace
$\operatorname{Spec}_{\widetilde{Z}}(V)\subset a_0 + {\mathbb{Z}}$
 for some $a_0 \in {\mathbb{C}}$.  
\newline
{\rm{2)}}\enspace 
If $\sup \operatorname{Re} \operatorname{Spec}_{\widetilde{Z}}(V)< \infty$, 
 then $V$ is a highest weight module.  
\newline
{\rm{3)}}\enspace 
If $V$ is a highest weight module $\hwm {\mathfrak g} {\lambda}$,
 then $\operatorname{Spec}_{\widetilde{Z}}(V) 
\subset -\mathbb{N} + \lambda(\widetilde Z)$
 and $\sup\operatorname{Re} \operatorname{Spec}_{\widetilde Z}(V)
= \operatorname{Re}\lambda(\widetilde Z)$.    
\newline
{\rm{4)}}\enspace
If $V$ is a unitary highest weight module,
 then $\operatorname{Spec}_{\widetilde Z}(V) \subset (-\infty, c_G]$, 
 where $c_G$ is a constant depending on $G$.      
\newline
{\rm{5)}}\enspace
If both $V$ and $F$ are highest weight modules of finite length, 
then any irreducible subquotient $W$
 of $V \otimes F$ is also a highest weight module.  
\end{lemma}
\begin{proof}
1)\enspace
For $a \in {\mathbb{C}}$, 
we write the eigenspace of $\widetilde Z$
 as 
$
     V_{a}:=\set{v \in V}{\widetilde Z v = a v}
$. 
Then, 
it follows from the Leibniz rule that
$$
     {\mathfrak {p}}_+ V_{a} \subset V_{a+1} \, , \quad
     {\mathfrak {k}}_{\mathbb{C}} V_{a} \subset V_{a} \, ,
     \quad\text{and}\quad 
     {\mathfrak {p}}_- V_{a} \subset V_{a-1} \, .   
$$
An iteration of this argument shows that 
$$
   \operatorname{Spec}_{\widetilde{Z}}(U({\mathfrak {g}}_{\mathbb{C}})V_{a})
    \subset a + {\mathbb{Z}} \, .    
$$ 
Now
 we take $a_0$ such that $V_{a_0}\ne \{0\}$.  
Since $V$ is irreducible, we have
$V=U(\mathfrak g_{\mathbb C})V_{a_0}$, 
and therefore
$\operatorname{Spec}_{\widetilde Z}(V) \subset a_0 + \mathbb Z$.   
\newline
2)\enspace
Suppose $\sup \operatorname{Re}\operatorname{Spec}_{\widetilde Z}(V) < \infty$.  
Since  $\operatorname{Re}\operatorname{Spec}_{\widetilde Z}(V)$ is 
  discrete by (1),
there exists $a \in \operatorname{Spec}_{\widetilde{Z}}(V)$ 
such that $\operatorname{Re}a$ attains its maximum.
Then 
$$
     {\mathfrak {p}}_+ V_{a} \subset  V_{a+1}=\{0\} \, .  
$$
Thus, 
 $V_{a} \subset V^{{\mathfrak{p}}_+}$.  
Hence, 
$V$ is a highest weight module.   
\newline
3)\enspace
The highest weight module $\hwm{\mathfrak {g}}{\lambda}$
 is isomorphic to the unique irreducible 
quotient  of the generalized Verma module
$N^{\mathfrak {g}}(\lambda)
 =U(\mathfrak g_{\mathbb C})
   \otimes _{U(\mathfrak k_{\mathbb{C}} + {\mathfrak p}_+)}\hwm{\mathfrak k}{\lambda}$.  
By the Poincar\'{e}--Birkhoff--Witt theorem,
$N^{\mathfrak{g}}(\lambda)$ is isomorphic to
$
S(\mathfrak{p}_-) \otimes \pi_\lambda^{\mathfrak{k}}
$
as a $\mathfrak{k}$-module.
Thus, any $\mathfrak k$-type $\hwm{\mathfrak k}{\mu}$
 occurring in $\hwm{\mathfrak g}{\lambda}$ is of the form
$$
     \mu = \lambda+ \sum_{\alpha \in \Delta({\mathfrak p}_-, {\mathfrak t})}m_{\alpha} \alpha
$$ 
for some $m_{\alpha} \in {\mathbb{N}}$.    
As $\alpha(\widetilde Z)=-1$ for any $\alpha \in \Delta(\mathfrak p_-, \mathfrak t)$, 
 we have
\begin{equation}
\label{eqn:eigenKG}
  \mu(\widetilde Z)=\lambda(\widetilde Z)-
  \sum_{\alpha\in\Delta(\mathfrak{p}_-,\mathfrak{t})} m_{\alpha} \, .  
\end{equation}
In particular, 
we have the following equivalence:
\begin{equation}
\label{eqn:Rekg}
\operatorname{Re}\mu(\widetilde Z)=\operatorname{Re}\lambda(\widetilde Z)
\ \Longleftrightarrow\ \mu=\lambda \, ,
\end{equation}
and we also  have 
\begin{equation}
\label{eqn:Ngl}
  \operatorname{Spec}_{\widetilde Z}(\hwm{\mathfrak g}{\lambda})
  \subset
  \set{\lambda(\widetilde Z)-
   \sum_{\alpha\in\Delta(\mathfrak{p}_-,\mathfrak{t})} 
   m_{\alpha}}{m_{\alpha} \in {\mathbb{N}}}
  =-{\mathbb N} + \lambda(\widetilde Z) \, .  
\end{equation}
Furthermore, 
since the $\mathfrak{k}$-type 
$\hwm{\mathfrak k} {\lambda}$ occurs in $\hwm{\mathfrak g}{\lambda}$, 
we have $\lambda(\widetilde{Z}) \in
 \operatorname{Spec}_{\widetilde Z}(\hwm{\mathfrak g}{\lambda})$.  
Here,  
$
     \sup \operatorname{Re}\operatorname{Spec}_{\widetilde Z}(\hwm{\mathfrak g}{\lambda})
      = \operatorname{Re} \lambda(\widetilde Z)
$.
\newline
4)\enspace
This statement follows from \eqref{eqn:ubdlmd} and from (3).
\newline
5)\enspace
For two subsets $A$ and $B$ in $\mathbb C$, 
we write 
 $A+B:=\set{a+b \in \mathbb C}{a \in A, b \in B}$.  
Then, 
$
     \operatorname{Spec}_{\widetilde Z}(V \otimes F) 
     \subset
     \operatorname{Spec}_{\widetilde Z}(V)
     + 
     \operatorname{Spec}_{\widetilde Z}(F)\, .  
$
Therefore, 
\begin{align*}
     \sup \operatorname{Re} \operatorname{Spec}_{\widetilde Z}(W)
&     \le
    \sup \operatorname{Re}\operatorname{Spec}_{\widetilde Z} (V)    
\\
&    \le
    \sup \operatorname{Re} \operatorname{Spec}_{\widetilde Z}(V)
    + 
    \sup \operatorname{Re} \operatorname{Spec}_{\widetilde Z}(F)
    < \infty.
\end{align*}  
Hence, 
 $W$ is also a highest weight module by (2).
\qed
\end{proof}

\subsection{Bottom layer map}
\label{subsec:ktog}

The following lemma finds an irreducible summand
(`bottom layer') from the $K$-type structure.

\begin{lemma}
\label{lem:ktog}
Let $V$ be a \gk-module.
We assume that $V$ decomposes into an algebraic direct sum of 
(possibly, infinitely many)
irreducible highest weight modules.  
We set 
$$
     \operatorname{Supp}_{\mathfrak k}(V)
   :=\set{\mu \in \sqrt{-1}\mathfrak t^*}
           {\Hom_{\mathfrak k}(\hwm{\mathfrak k}{\mu}, V)\ne \{0\}} \, .  
$$
If the evaluation map
$$
    \operatorname{Supp}_{\mathfrak k}(V) \to \mathbb{R}\, , \quad
    \mu \mapsto \operatorname{Re} \mu(\widetilde Z)
$$
attains its maximum at $\mu_0$, 
then 
$$
     \Hom_{(\mathfrak g_{\mathbb C}, K)}(\hwm{\mathfrak g}{\mu_0},
     V)\ne\{0\} \, .  
$$
\end{lemma}
\begin{proof}
Take a non-zero map 
 $q \in \Hom_{\mathfrak k}(\hwm{\mathfrak k}{\mu_0}, V)$.  
As $V$ is an algebraic direct sum
 of irreducible highest weight modules, 
there exists a projection 
 $p:V \to \hwm{\mathfrak g}{\lambda}$
 for some $\lambda$
 such that $p \circ q \ne 0$.  
This means that  $\hwm{\mathfrak k}{\mu_0}$ occurs
 in $\hwm{\mathfrak g}{\lambda}$, 
and therefore we have 
$$
     \operatorname{Re}\mu_0(\widetilde Z)
     \le \sup \operatorname{Re}\operatorname{Spec}_{\widetilde Z}(\hwm{\mathfrak g}{\lambda})
                             =\operatorname{Re}\lambda(\widetilde Z)
     \, .
$$
Here, the last equality is by Lemma~\ref{lem:spec}~(3).

Conversely, 
the maximality of $\mu_0$
 implies that 
$\operatorname{Re}\mu_0(\widetilde Z)\ge \operatorname{Re}\lambda(\widetilde Z)$.  
Hence,
$\operatorname{Re} \mu_0 (\widetilde{Z})
 = \operatorname{Re} \lambda(\widetilde{Z})$,
and we have then
 $\mu_0 = \lambda$
by \eqref{eqn:Rekg}.
Since $\hwm{\mathfrak{g}}{\lambda}$ is an irreducible summand of $V$,
we have  
 $\Hom_{(\mathfrak g_{\mathbb C}, K)}(\hwm{\mathfrak g}{\mu_0}, V) \ne 
  \{0\}$.    
\qed
\end{proof}
\subsection{Determination of the $\mathfrak{g}_{\mathbb{C}}$-structure 
by $K$-types}
\label{subsec:8.11}
In general,
the $K$-type formula is not sufficient to determine the irreducible
decomposition of a unitary representation even in the discretely
decomposable case.
However, this is the case if any irreducible summand is a highest
weight module. 
Here is the statement that we shall use as a main machinery of the proof of
Theorems~\ref{thm:gHKS} and \ref{thm:tensordeco}.
\begin{lemma}
\label{lem:HnG}
Suppose $(\pi, \mathcal H)$ is a $K$-admissible unitary representation of $G$,
 which splits discretely into a Hilbert direct sum
 of irreducible unitary highest weight representations of $G$.  
Let $\mathcal{H}_K$ be the space of $K$-finite vectors of
 $\mathcal{H}$. 
Assume that there exists a function 
$n_\pi: \mathfrak{t}_{\mathbb{C}}^* \to \mathbb{N}$
such that $\mathcal{H}_K$ is isomorphic to the following
direct sum as $\mathfrak{k}$-modules:
\begin{equation}
\label{eqn:Hnalg}
    \mathcal H_K \simeq \bigoplus_{\lambda} n_{\pi}(\lambda) \hwm{\mathfrak g}{\lambda} 
     \quad
     \text{(algebraic direct sum)}.  
\end{equation}
Then,
$n_\pi(\lambda)\ne0$
only if 
$\lambda \in \Lambda_G$, that is,
$\hwm{\mathfrak{g}}{\lambda}$ lifts to an irreducible unitary
representation $\hwm{G}{\lambda}$ of $G$.
Furthermore,
the identity \eqref{eqn:Hnalg} holds as a 
$(\mathfrak{g}_{\mathbb{C}},K)$-module isomorphism, 
and 
the unitary representation $\pi$ has the following decomposition into
irreducible unitary representations of $G$:
\begin{equation}
\label{eqn:HnG}
     \pi \simeq {\sum_{\lambda}}^{\oplus} n_{\pi}(\lambda) \hwm G{\lambda}
     \quad
     \text{(discrete Hilbert sum)}.
\end{equation}
\end{lemma}
\begin{proof}
We write an abstract irreducible decomposition of $\mathcal H$
 as 
\begin{equation*}
    \mathcal H \simeq {\sum_{\lambda\in\Lambda_G}}^{\oplus} m_{\lambda} \hwm
    G{\lambda} \quad\text{(discrete Hilbert sum)}.
\end{equation*}
Since $\mathcal{H}$ is $K$-admissible,
the multiplicity $m_{\lambda}< \infty$ 
for all $\lambda$,
and
 we have an isomorphism of \gk-modules
with the same multiplicity $m_\lambda$
 (see \cite[Theorem 2.7]{xkaspm}):
\begin{equation}
\label{eqn:Hmk}
     \mathcal H_K \simeq \bigoplus_{\lambda\in\Lambda_G} 
m_{\lambda} \hwm{\mathfrak g}{\lambda} 
     \quad
     (\text{algebraic direct sum}).
\end{equation}
Let us show $n_\pi(\lambda) = m_\lambda$ for all $\lambda$.
For this,
we begin with an observation that
$$
\operatorname{Spec}_{\widetilde{Z}} (\mathcal{H}_K)
= \bigcup_{\substack{\lambda\text{ such that}\\
                     m_\lambda\ne0}}
   \operatorname{Spec}_{\widetilde{Z}} (\hwm{\mathfrak{g}}{\lambda})
$$
is a subset in $\mathbb{R}$ and has an upper bound.
This follows from
Lemma~\ref{lem:spec} (4) applied to each irreducible summand in
\eqref{eqn:Hmk}.

First, we consider the case where
there exists $a \in \mathbb{R}$ such that
\begin{equation}
\label{eqn:congr}
\lambda(\widetilde{Z})\equiv a \bmod \mathbb{Z}
\quad\text{for any $\lambda$ satisfying $n_\pi(\lambda)\ne0$.}
\end{equation}
Then, the set
\begin{equation}
\label{eqn:wtH}
\{ \lambda(\widetilde{Z}): \lambda\in\mathfrak{t}^*_{\mathbb{C}}, 
   n_\pi(\lambda)\ne 0 \}
\end{equation}
is contained in 
$\operatorname{Spec}_{\widetilde{Z}}(\mathcal{H}_K)$
by \eqref{eqn:Hnalg},
and is  discrete by \eqref{eqn:congr}.
Hence,
it is a discrete subset of $\mathbb{R}$ with an upper bound.
Thus, 
we can find $\mu_0 \in \mathfrak{t}_{\mathbb{C}}^*$
such that $n_\pi(\mu_0) \ne 0$
and that
$\mu_0(\widetilde{Z})$ attains its maximum in \eqref{eqn:wtH}.
In turn, 
the evaluation map $\operatorname{Supp}_{\mathfrak k}(\mathcal H_K) \to \mathbb R$, 
$\mu \mapsto \mu(\widetilde Z)$
  attains its maximum at $\mu_0 \in \operatorname{Supp}_{\mathfrak
k}(\mathcal H_K)$
by \eqref{eqn:Hnalg} and Lemma~\ref{lem:spec} (3).
Therefore, 
$\Hom_{(\mathfrak g_{\mathbb C}, K)}(\hwm{\mathfrak g}{\mu_0}, \mathcal H_K)\ne \{0\}$
 by Lemma \ref{lem:ktog}.
Thus, we have shown $m_{\mu_0} \ne 0$,
that is,
$\hwm{G}{\mu_0}$ occurs as a subrepresentation in $\mathcal{H}$.

Next, we consider the unitary representation $\pi'$ on 
$$
\mathcal{H}' := 
\sideset{}{^\oplus} \sum_{\lambda\ne\mu_0}
   m_\lambda \hwm{G}{\lambda}
   \oplus (m_{\mu_0} -1) \hwm{G}{\mu_0} \, ,
$$
the orthogonal complement of a subrepresentation
$\hwm{G}{\mu_0}$ in $\mathcal{H}$.
Then, 
the $K$-type formula \eqref{eqn:Hnalg} for
$(\pi',\mathcal{H}')$ holds if we set
$$
n_{\pi'} (\lambda)
  := \begin{cases}
        n_\pi (\lambda) - 1   & (\lambda = \mu_0)\, , \\
        n_\pi (\lambda)       & (\lambda \ne \mu_0)\, .
     \end{cases} 
$$
Hence,
by the downward induction on
$\sup\operatorname{Spec}_{\widetilde{Z}}
 (\mathcal{H}_K)$,
we have
$n_\pi(\lambda) = m_\lambda$ for all $\lambda$.

For the general case, 
 let $A$ be the set of complete representatives
 of $\{\lambda(\widetilde{Z}) \in \mathbb{C} \mod \mathbb Z:
 n_\pi({\lambda}) \ne 0\}$.  
For each $a \in A$,
we define a subrepresentation $\mathcal{H}_a$ of $\mathcal{H}$ by
$$
     \mathcal H_{a} 
  :=\sideset{}{^\oplus}
\sum_{\lambda(\widetilde{Z}) \equiv a \operatorname{mod}\mathbb Z}
 m_{\lambda}\hwm G{\lambda}
\quad\text{(discrete Hilbert sum)}.  
$$
Then, 
we have an isomorphism of unitary representations of $G$:
$$
   \mathcal H \simeq {\sum_{a \in A}}^{\oplus} \mathcal
   H_{a} \, .
$$
Since
$\operatorname{Spec}_{\widetilde{Z}}(\hwm{\mathfrak{g}}{\lambda}) 
 \subset a+\mathbb{Z}
$
if and only if
$\lambda(\widetilde{Z})\equiv a \bmod \mathbb{Z}$
by Lemma~\ref{lem:spec} (3),
we get from \eqref{eqn:Hnalg} the following $K$-isomorphism
\begin{equation}
\label{eqn:Hnalga}
(\mathcal{H}_a)_K \simeq
   \bigoplus_{\lambda(\widetilde{Z})\equiv a\bmod\mathbb{Z}}
   n_\pi(\lambda) \hwm{\mathfrak{g}}{\lambda}
\end{equation}
for each $a \in A$.
Therefore, 
our proof for the first step assures 
$n_\pi(\lambda)=m_\lambda$ for any $\lambda$
such that $\lambda\equiv a \bmod \mathbb{Z}$.
Since $a \in A$ is arbitrary,
we obtain Lemma in the general case.
\qed
\end{proof}

\subsection{Proof of Theorem~\ref{thm:gHKS}}
\label{subsec:8.3}

In this section, we give
 a proof of Theorem~\ref{thm:gHKS}.
This is done by showing a more general formula in Lemma~\ref{lem:8.3}
 without 
the scalar type assumption \eqref{eqn:8.5.1}. 
Then, Theorem~\ref{thm:gHKS} follows readily from
 Lemma~\ref{lem:8.3} because the  
 assumption \eqref{eqn:8.5.1} makes 
$\dim \hwm{\mathfrak{k}}{\mu} = 1$
and
$\mathbb{S}_{(a_1,\dots,a_l)} (\mu) =
 \{ \mu - \sum_{j=1}^l a_j \nu_j \}$
(see \eqref{eqn:Samu} for notation).

For a discussion below,
it is convenient to use the concept of a multiset.
Intuitively, a multiset is a set counted with multiplicities;
for example,
$\{a,a,a,b,c,c\}$.
More precisely,
a multiset $\mathbb{S}$ consists of a set $S$ and a function
$m: S \to \{ 0,1,2,\dots,\infty \}$.
If $\mathbb{S}' = \{S,m'\}$
is another multiset on $S$ such that 
$m'(x) \le m(x)$ for all $x \in S$,
we say $\mathbb{S}'$ is a {\it{submultiset}} of $\mathbb{S}$ and write
$\mathbb{S}' \subset \mathbb{S}$.

Suppose we are in the setting of Subsection~\ref{subsec:8.2}
and recall $\tau$ is an involution of holomorphic type. 
For a $\rt^+(\mathfrak{k}, \mathfrak{t})$-dominant weight $\mu$,
 we introduce a multiset 
$\mathbb{S}(\mu)$ 
 consisting of $\rt^+(\mathfrak{k}^\tau, \mathfrak{t}^\tau)$-dominant 
weights: 
$$
\mathbb{S}(\mu):=
\bigcup_{\Sb a_1 \ge \dots \ge a_l \ge 0\\ a_1, \dots, a_l \in \mathbb{N} \endSb}
\mathbb{S}_{(a_1,\dots,a_l)} (\mu) \, ,
$$
where we define the multiset
$\mathbb{S}_{(a_1,\dots,a_l)} (\mu)$
by
\begin{equation}
\label{eqn:Samu}
\parbox{25em}{$\{$highest weight of irreducible $\mathfrak{k}^\tau$-modules
occurring in\newline
\hspace*{1em} $\hwm{\mathfrak{k}^\tau}{-\sum_{j=1}^l a_j \nu_j} \otimes
\hwm{\mathfrak{k}}{\mu}|_{\mathfrak{k}^\tau} $
counted with multiplicities$\}$.}
\end{equation}
Because 
the central element $\widetilde{Z} = \frac{1}{\sqrt{-1}} Z$ of 
$\mathfrak{k}_{\mathbb{C}}$ 
acts on the irreducible representation 
$\hwm{\mathfrak{k}}{\mu}$ by the scalar $\mu(\widetilde{Z})$
and because $\nu_j (\widetilde{Z}) = 1$ for all $j$ $(1 \le j \le l)$, 
any element $\nu$ in $\mathbb{S}_{(a_1,\dots,a_l)}(\mu)$
satisfies
$
\nu(\widetilde{Z}) = - \sum_{j=1}^l a_j + \mu(\widetilde{Z})
$.
Therefore, the multiplicity of each element of the multiset
$\mathbb{S}(\mu)$ is finite.

\begin{lemma}
\label{lem:8.3}
Let $\tau$ be an involution of $G$
 of holomorphic type,
 and $H = G_0^\tau$.
If $\hwm{G}{\mu}$ is 
a (relative) holomorphic discrete series representation of $G$,
 then it decomposes discretely into irreducible
 $H$-modules as:
\begin{equation*}
  \hwm{G}{\mu} |_H
 \simeq
\sideset{}{^\oplus}\sum_{\nu \in \mathbb{S}(\mu)}
       \hwm{H}{\nu}
\quad
\text{(discrete Hilbert sum).}
\end{equation*}
\end{lemma}
{
\renewcommand{\proofname}{Proof of Lemma~\ref{lem:8.3}} %
\begin{proof}
It follows from Fact~\ref{fact:3.4.1}~(1)
that $\hwm{G}{\mu}$ is $(H\cap K)$-admissible,
and splits discretely into a Hilbert direct sum of irreducible unitary
representations of $H$.

Applying Lemma~\ref{lem:HnG} to $H = G_0^\tau$,
we see that Lemma~\ref{lem:8.3} is deduced from the following
$\mathfrak{k}^\tau$-isomorphism: 
\begin{equation}
\label{eqn:kgeneral}
  \hwm{\mathfrak{g}}{\mu}
 \simeq
\bigoplus_{\nu \in \mathbb{S}(\mu)}\hwm{\mathfrak{g}^\tau}{\nu}
\quad
\text{(algebraic direct sum).}
\end{equation}
The rest of the proof is devoted to showing \eqref{eqn:kgeneral}.

Since $\hwm{G}{\mu}$ is a holomorphic discrete series,
$\hwm{\mathfrak{g}}{\mu}$ is isomorphic to the generalized Verma
module $N^{\mathfrak{g}}(\mu)
 = U(\mathfrak{g}_\mathbb{C}) 
\otimes_{U(\mathfrak{k}_\mathbb{C} + \mathfrak{p}_+)} \hwm{\mathfrak{k}}{\mu}
$ as a $\mathfrak{g}$-module, 
which in turn is isomorphic to the $\mathfrak{k}$-module 
$
S(\mathfrak{p}_-) \otimes \hwm{\mathfrak{k}}{\mu}
$.

According to the decomposition
$\mathfrak{p}_- = \mathfrak{p}_-^\tau \oplus \mathfrak{p}_-^{-\tau}$
as $\mathfrak{k}^\tau$-modules, 
we have then the following $\mathfrak{k}^\tau$-isomorphism:
\begin{equation}
\label{eqn:pipp}
\hwm{\mathfrak{g}}{\mu}
 \simeq S(\mathfrak{p}_-) \otimes \hwm{\mathfrak{k}}{\mu}
 \simeq S(\mathfrak{p}_-^\tau) \otimes
         S(\mathfrak{p}_-^{-\tau}) \otimes
         \hwm{\mathfrak{k}}{\mu}  \, .
\end{equation}
Now, we consider
the Hermitian symmetric space
$G^{\tau\theta} / G^{\tau,\theta}$,
for which the complex structure is given by the decomposition
$\mathfrak{g}_{\mathbb{C}}^{\tau\theta} = 
 \mathfrak{k}_{\mathbb{C}}^\tau \oplus
 \mathfrak{p}_+^{-\tau} \oplus
 \mathfrak{p}_-^{-\tau}
$
(see  \eqref{eqn:8.2.2}).
Then, the  Hua--Kostant--Schmid formula 
(\cite[ Behauptung c]{xschmidherm}) applied
 to $G^{\tau\theta}/G^{\tau,\theta}$
 decomposes the symmetric algebra 
 $S(\mathfrak{p}_-^{-\tau})$ 
 into irreducible 
$\mathfrak{k}^\tau$-modules:
\begin{equation}
\label{eqn:HKStau}
 S(\mathfrak{p}_-^{-\tau})
 \simeq
\bigoplus
  \Sb a_1 \ge \dots \ge a_l \ge 0\\ a_1, \dots, a_l \in \mathbb{N} \endSb
       \hwm{\mathfrak{k}^\tau}{- \sum_{j=1}^l a_j \nu_j} \, .
\end{equation}
It follows from the definition of $\mathbb{S}(\mu)$ that we have the
following irreducible decomposition as $\mathfrak{k}^\tau$-modules:
$$
S(\mathfrak{p}_-^{-\tau}) \otimes \hwm{\mathfrak{k}}{\mu}
 \simeq \bigoplus_{\nu\in\mathbb{S}(\mu)}
      \hwm{\mathfrak{k}^\tau}{\nu}.
$$
Combining this with \eqref{eqn:pipp},
we get a $\mathfrak{k}^\tau$-isomorphism
$$
\hwm{\mathfrak{g}}{\mu}
 \simeq \bigoplus_{\nu\in\mathbb{S}(\mu)}
         S(\mathfrak{p}_-^\tau) \otimes
         \hwm{\mathfrak{k}^\tau}{\nu} \, .
$$
Next, we consider the Verma module
$N^{\mathfrak{g}^\tau}(\nu) =
 U(\mathfrak{g}^\tau_{\mathbb{C}})
 \otimes_{U(\mathfrak{k}^\tau_{\mathbb{C}}+\mathfrak{p}^\tau_+)}
 \hwm{\mathfrak{k}^\tau}{\nu}$ of the subalgebra $\mathfrak{g}^\tau$.
Then,
$\hwm{\mathfrak{g}^\tau}{\nu}$
is the unique irreducible quotient of
$N^{\mathfrak{g}^\tau}(\nu)$.
We shall show later that
$N^{\mathfrak{g}^\tau}(\nu)$
is irreducible as a $\mathfrak{g}^\tau$-module,
but at this stage we denote by
$\hwm{\mathfrak{g}^\tau}{\nu},
 \hwm{\mathfrak{g}^\tau}{\nu'},
 \hwm{\mathfrak{g}^\tau}{\nu''},
 \ldots$
the totality of irreducible subquotient modules of 
$N^{\mathfrak{g}^\tau}(\nu)$. 
(There are at most finitely many subquotients,
and all of them are highest weight modules.) 
Then, as $\mathfrak{k}^\tau$-modules,
we have the following isomorphisms:
\begin{align*}
S(\mathfrak{p}_-^\tau) \otimes \hwm{\mathfrak{k}^\tau}{\nu}
& \simeq N^{\mathfrak{g}^\tau}(\nu)
\\
& \simeq \hwm{\mathfrak{g}^\tau}{\nu} \oplus
         \hwm{\mathfrak{g}^\tau}{\nu'} \oplus
         \hwm{\mathfrak{g}^\tau}{\nu''} \oplus \cdots \, .
\end{align*}
Therefore, we get a $\mathfrak{k}^\tau$-isomorphism:
$$
\hwm{\mathfrak{g}}{\mu}
 \simeq \bigoplus_{\nu\in\mathbb{S}(\mu)}
     (\hwm{\mathfrak{g}^\tau}{\mu} \oplus
      \hwm{\mathfrak{g}^\tau}{\nu'} \oplus
      \hwm{\mathfrak{g}^\tau}{\nu''} \oplus \cdots) \, .
$$
Accordingly, 
the restriction $\hwm{G}{\mu}|_H$ splits discretely into irreducible
unitary representations of $H$ by Lemma~\ref{lem:HnG}:
$$
\hwm{G}{\mu} |_H
   \simeq \; \sideset{}{^\oplus}\sum_{\nu \in \mathbb{S}(\mu)}
     (\hwm{H}{\nu} \oplus
      \hwm{H}{\nu'} \oplus
      \hwm{H}{\nu''} \oplus \cdots ) \, .
$$
Since $\hwm{G}{\mu}$ is a (relative) holomorphic discrete series
representation of $G$,
all irreducible summands in the right-hand side must be (relative)
holomorphic discrete series representations of $H$ by
Fact~\ref{fact:3.4.1} (1).
Therefore,
$N^{\mathfrak{g}^\tau}(\nu)$ is irreducible,
and the other subquotients 
$\hwm{\mathfrak{g}^\tau}{\nu'}, \hwm{\mathfrak{g}^\tau}{\nu''},
  \ldots$ 
do not appear.
Hence, 
the $\mathfrak{k}^\tau$-structures of the both sides of
\eqref{eqn:kgeneral} are the same.
Thus, Lemma~\ref{lem:8.3} is proved.
\qed
\end{proof}
}

\subsection{Proof of Theorem~\ref{thm:tensordeco}}
\label{subsec:pf tensordeco}
For two irreducible representations
$\hwm{\mathfrak{k}}{\mu_1}$ and
$\hwm{\mathfrak{k}}{\mu_2}$,
we define a multiset
$\mathbb{S}(\mu_1,\mu_2)$
consisting of $\Delta^+(\mathfrak{k},\mathfrak{t})$-dominant weights
by 
$$
\mathbb{S}(\mu_1,\mu_2)
  := \bigcup_{\substack{a_1\ge\dots\ge a_{\bar{l}}\ge0 \\
                        a_1,\dots,a_{\bar{l}}\in\mathbb{N} }}
     \mathbb{S}_{(a_1,\dots,a_{\bar{l}})}
      (\mu_1,\mu_2) \, ,
$$
where
$\mathbb{S}_{(a_1,\dots,a_{\bar{l}})} (\mu_1,\mu_2)$
is the multiset consisting of highest weights of irreducible
$\mathfrak{k}$-modules occurring in
$\hwm{\mathfrak{k}}{-\sum_{j=1}^{\bar{l}} a_j\bar{\nu}_j}
 \otimes \hwm{\mathfrak{k}}{\mu_1} 
 \otimes \hwm{\mathfrak{k}}{\mu_2}
$
counted with multiplicities.

Theorem~\ref{thm:tensordeco} 
 is derived from the following more general formula:

\begin{lemma}
\label{lem:tensordeco}
The tensor product of two (relative)
holomorphic discrete series representations 
$\hwm{G}{\mu_1}$ and $\hwm{G}{\mu_2}$ 
decomposes as follows:
$$
\hwm{G}{\mu_1} \widehat{\otimes} \, \hwm{G}{\mu_2}
\simeq \;
  \sideset{}{^\oplus}\sum_{\nu \in \mathbb{S}(\mu_1,\mu_2)}
         \hwm{G}{\nu} \, .
$$
\end{lemma}

\begin{proof}
We define two injective maps by: 
\begin{align*}
&\diag: \mathfrak{p}_+ \to \mathfrak{p}_+ \oplus \mathfrak{p}_+ \, ,
\quad
X \mapsto (X,X) \, ,
\\
&\diag': \mathfrak{p}_+ \to \mathfrak{p}_+ \oplus \mathfrak{p}_+ \, ,
\quad
X \mapsto (X,-X) \, .
\end{align*}
It then follows that we have $\mathfrak{k}$-isomorphisms:
\begin{align*}
S(\mathfrak{p}_-) \otimes S(\mathfrak{p}_-)
& \simeq S(\mathfrak{p}_- \oplus \mathfrak{p}_-)
\\
& \simeq S(\diag(\mathfrak{p}_-)) \otimes
         S(\diag'(\mathfrak{p}_-) )
\\
& \simeq \bigoplus_{\substack{a_1 \ge\cdots\ge a_{\bar{l}}\ge0\\
                              a_1,\dots,a_{\bar{l}}\in\mathbb{N}}}
          S(\diag(\mathfrak{p}_-))
         \otimes
         \hwm{\mathfrak{k}}{-\sum^{\bar{l}}_{j=1}a_j\bar{\nu}_j} \, . 
\end{align*}
This brings us the following $\mathfrak{k}$-isomorphisms:
\begin{align*}
\hwm{\mathfrak{g}}{\mu_1} \otimes \hwm{\mathfrak{g}}{\mu_2}
& \simeq S(\mathfrak{p}_-) \otimes \hwm{\mathfrak{k}}{\mu_1}
     \otimes S(\mathfrak{p}_-) \otimes \hwm{\mathfrak{k}}{\mu_2}
\\
& \simeq \bigoplus_{\nu\in\mathbb{S}(\mu_1,\mu_2)}
   S(\diag(\mathfrak{p}_-)) \otimes \hwm{\mathfrak{k}}{\nu} 
\\
& \simeq \bigoplus_{\nu\in\mathbb{S}(\mu_1,\mu_2)}
     N^{\mathfrak{g}}_\nu \, .
\end{align*}
The rest of the proof goes similarly to that of Lemma~\ref{lem:8.3}.
\qed
\end{proof}

\subsection{Restriction $U(p,q) \downarrow U(p-1,q)$
and $SO(n,2) \downarrow SO(n-1,2)$}
\label{subsec:8.4}

Suppose $(G,H)$ is a reductive symmetric pair whose complexification 
 $(\mathfrak{g}_\mathbb{C}, \mathfrak{h}_\mathbb{C})$
 is one of the following types:

\smallskip
$(\mathfrak{sl}(n, \mathbb{C}), \mathfrak{gl}(n-1, \mathbb{C}))$ 
(or
 $(\mathfrak{gl}(n,\mathbb{C}),\mathfrak{gl}(1,\mathbb{C})
     +\mathfrak{gl}(n-1,\mathbb{C}))$),

\smallskip
 $(\mathfrak{so}(n, \mathbb{C}), \mathfrak{so}(n-1, \mathbb{C}))$.

\smallskip
\noindent
As is classically known (see \cite{xvk}),
for compact $(G,H)$
such as $(U(n),U(1) \times U(n-1))$ or $(SO(n),SO(n-1))$, 
any irreducible finite dimensional representation $\pi$ of $G$
is multiplicity-free when restricted to $H$.
For non-compact $(G,H)$
such as $(U(p,q), U(1)\times U(p-1,q))$ or
$(SO(n,2), SO(n-1,2))$,
an analogous theorem still holds for highest weight representations $\pi$:
\begin{theorem}
\label{thm:minus1}
If\/
 $(\mathfrak{g}, \mathfrak{h}) = ({\mathfrak{u}}(p,q), 
  \mathfrak{u}(1)+{\mathfrak{u}}(p-1,q))$
 or $(\mathfrak{so}(n,2), \mathfrak{so}(n-1,2))$,
then 
 any irreducible unitary highest weight representation
 of $G$  decomposes discretely into a multiplicity-free sum of irreducible 
 unitary highest weight representations of $H$. 
\end{theorem}

In contrast to Theorem~\ref{thm:A},
the distinguishing feature of Theorem~\ref{thm:minus1}
 is that $\pi$ is not necessarily of scalar type
but an arbitrary unitary highest weight module.
The price to pay is that the pair $(G, H)$ is very special.
We do not give the proof here that uses the vector bundle version of 
Theorem~\ref{thm:2.2} (see \cite{mfbdle}).
Instead, we give an explicit decomposition formula for holomorphic discrete
series $\pi$.
The proof of Theorem~\ref{thm:minus1} for the 
case $(G,H) = (SO_0(n,2), SO_0(n-1,2))$ can be also found in
 Jakobsen and Vergne
 \cite[Corollary~3.1]{xjv}.

\subsection{Branching law for $U(p,q)\downarrow U(p-1,q)$}
This subsection gives an explicit branching law of a holomorphic discrete 
series representation $\hwm{G}{\mu}$ of
$G=U(p,q)$ when restricted to
$H=U(1)\times U(p-1,q)$.
Owing to \eqref{eqn:8.5.2},
such $\hwm{G}{\mu}$ is parametrized by
$\mu=(\mu_1,\dots,\mu_{p+q}) \in \mathbb{Z}^{p+q}$
satisfying
$$
\mu_1\ge\cdots\ge\mu_p, \mu_{p+1}\ge\cdots\ge\mu_{p+q},
   \mu_{p+q}\ge\mu_1 +p+q \, .
$$
Here is the formula:
\begin{theorem}[Branching law $U(p,q)\downarrow U(p-1,q)$]
\label{thm:upqupq}
Retain the above setting.
Then, 
 the branching law of $\hwm{G}{\mu}$ of the restriction to the
subgroup
$H$
is multiplicity-free for any $\mu$; it is given as follows:
\begin{equation}
\label{holoupq}
\hwm{G}{\mu}|_H
\simeq {}
    \sideset{}{^\oplus}\sum_{a=0}^\infty 
    \sideset{}{^\oplus}\sum_{\substack{
           {\mu_1\ge\lambda_2\ge\mu_2\ge\cdots\ge\lambda_p\ge\mu_p} \\
           {\lambda_{p+1}\ge\mu_{p+1}\ge\cdots\ge\lambda_{p+q}\ge\mu_{p+q}}\\
           {\sum_{i=1}^q(\lambda_{p+i}-\mu_{p+i})=a}}}
      \mathbb{C}_{\sum_{i=1}^p \mu_i - \sum_{i=1}^p \lambda_i - a}
 \boxtimes
\hwm{U(p-1,q)}{(\lambda_2,\dots,\lambda_p,\lambda_{p+1},\dots,\lambda_{p+q})}
\, . 
\end{equation}
\end{theorem}

\begin{proof}
For 
$
(G,H) \equiv (G,G^\tau)
= (U(p,q), U(1) \times U(p-1,q))
$,
we have
\begin{alignat*}{2}
&G^{\tau\theta} 
&&\simeq U(1,q) \times U(p-1),
\\
&H \cap K
\ (= K^\tau = K^{\tau\theta})
&&\simeq U(1) \times U(p-1) \times U(q),
\end{alignat*}
$\mathfrak{t}^\tau = \mathfrak{t}$,
and
$$
\Delta^+ (\mathfrak{p}_+^{-\tau}, \mathfrak{t}^\tau)
= \{ e_1 - e_{p+i} : 1 \le i \le q \}
$$
by using the standard basis of
$\Delta(\mathfrak{g}, \mathfrak{t})
= \{ \pm (e_i - e_j) :
1 \le i < j \le p+q \}$.
Thus,
$l = \rrank G^{\tau\theta} = 1$
and $\nu_1 = e_1 - e_{p+1}$.
Hence, the $\mathfrak{k}^\tau$-type formula
\eqref{eqn:HKStau}
amounts to 
\begin{align}
\label{eqn:Supq}
S (\mathfrak{p}_-^{-\tau})
&\simeq  \bigoplus_{a=0}^\infty
 \hwm{H\cap K}{-a(e_1 - e_{p+1})}
\nonumber
\\
&\simeq  \bigoplus_{a=0}^\infty
 \mathbb{C}_{-a} \boxtimes \mathbf{1} \boxtimes
 \hwm{U(q)}{(a,0,\dots,0)}
\end{align}
as $H \cap K \simeq U(1) \times U(p-1) \times U(q)$ modules.
Here, $\mathbf{1}$ denotes the trivial one dimensional representation
of $U(p-1)$.

On the other hand, 
we recall a classical branching formula
$U(p) \downarrow U(p-1)$:
$$
\hwm{U(p)}{(\mu_1,\dots,\mu_p)} |_{U(1)\times U(p-1)}
\simeq \bigoplus_{\mu_1 \ge \lambda_2 \ge \mu_2 \ge
                  \cdots \ge \lambda_p \ge
                  \mu_p} 
\mathbb{C}_{\sum_{i=1}^p \mu_i - \sum_{i=2}^p \lambda_i}
\otimes \hwm{U(p-1)}{(\lambda_2,\dots,\lambda_p)} \, ,
$$
whereas the classical Pieri rule says
$$
\hwm{U(q)}{(a,0,\dots,0)} \otimes \hwm{U(q)}{(\mu_{p+1},\dots,\mu_{p+q})}
\simeq
\bigoplus_{\substack{\lambda_{p+1}\ge\mu_{p+1}\ge\cdots
                      \ge\lambda_{p+q}\ge\mu_{p+q}\\
                    {\sum_{i=1}^q (\lambda_{p+i}-\mu_{p+i})=a}}}
 \hwm{U(q)}{(\lambda_{p+1},\dots,\lambda_{p+q})} \, .
$$
These two formulae together with \eqref{eqn:Supq} yield
the following $\mathfrak{k}^\tau$-isomorphisms:
\begin{align*}
&S(\mathfrak{p}_-^{-\tau}) \otimes
 \hwm{\mathfrak{k}}{\mu} |_{\mathfrak{k}^\tau}
\\
&{}\simeq
 \bigoplus_{a=0}^\infty
   (( \mathbb{C}_{-a} \boxtimes \mathbf{1}) \otimes
      \hwm{U(p)}{(\mu_1,\dots,\mu_p)}
   |_{U(1)\times U(p-1)} )
 \boxtimes (\hwm{U(q)}{(a,0,\dots,0)} \otimes
            \hwm{U(q)}{(\mu_{p+1},\dots,\mu_{p+q})})
\\
&{}\simeq
 \bigoplus_{a=0}^\infty \ 
 \bigoplus_{\substack{
           {\mu_1\ge\lambda_2\ge\mu_2\ge\cdots\ge\lambda_p\ge\mu_p} \\
           {\lambda_{p+1}\ge\mu_{p+1}\ge\cdots\ge\lambda_{p+q}\ge\mu_{p+q}}\\
           {\sum_{i=1}^q(\lambda_{p+i}-\mu_{p+i})=a}}}
 \mathbb{C}_{\sum_{i=1}^p\mu_i-\sum_{i=2}^p\lambda_i-a}
 \boxtimes
 \hwm{U(p-1)}{(\lambda_2,\dots,\lambda_p)}
 \boxtimes
 \hwm{U(q)}{(\lambda_{p+1},\dots,\lambda_{p+q})} \, .
\end{align*}
In view of the $\mathfrak{k}^\tau$-isomorphisms
$$
\hwm{\mathfrak{g}}{\mu}\simeq S(\mathfrak{p}^\tau_-) \otimes
   S(\mathfrak{p}^{-\tau}_-) \otimes \hwm{\mathfrak{k}}{\mu}|_{\mathfrak{k}^\tau}
$$
and $N^{\mathfrak{g}^\tau} (\nu) \simeq S(\mathfrak{p}^\tau_-) \otimes
      \hwm{\mathfrak{k}^\tau}{\nu}$,
we have now shown that the $\mathfrak{k}^\tau$-structure of
$\hwm{\mathfrak{g}}{\mu}$ coincides with that of
$$
\bigoplus_{a=0}^\infty \ 
 \bigoplus_{\substack{
           {\mu_1\ge\lambda_2\ge\mu_2\ge\cdots\ge\lambda_p\ge\mu_p} \\
           {\lambda_{p+1}\ge\mu_{p+1}\ge\cdots\ge\lambda_{p+q}\ge\mu_{p+q}}\\
           {\sum_{i=1}^q(\lambda_{p+i}-\mu_{p+i})=a}}}
N^{\mathfrak{g}^\tau}
   (\sum_{i=1}^p \mu_i - \sum_{i=2}^p \lambda_i - a,
     \lambda_2, \dots, \lambda_{p+q}) \, .
$$
As in the last part of the proof of Theorem~\ref{thm:gHKS},
we see that any generalized Verma module occurring in the right-hand side 
is irreducible
(and is isomorphic to the underlying
$(\mathfrak{g}^\tau_{\mathbb{C}}, H\cap K)$-module of a holomorphic discrete 
series of $H$).
Therefore, Theorem follows from Lemma~\ref{lem:HnG}.
\qed
\end{proof}

\section{Appendix: Associated Bundles on Hermitian Symmetric Spaces}
\label{sec:9}

In this Appendix,
 we explain standard operations on homogeneous vector bundles.
The results are well-known and elementary,
but we recall them briefly
 for the convenience of the reader.
The main goal is Lemma~\ref{lem:9.6} %
which is used to verify the condition \eqref{eqn:2.2.2} in 
Theorem~\ref{thm:2.2}.

\subsection{Homogeneous vector bundles}
\label{subsec:9.1}
Let $M$ be a real manifold,
 and $V$ a (finite dimensional) vector space over $\mathbb{C}$.
Suppose that we are given an open covering $\{ U_\alpha \}$ of $M$
 and transition functions
$$
  g_{\alpha \beta} : U_\alpha \cap U_\beta \to GL_\mathbb{C}(V)
$$
 satisfying the following compatibility conditions:
$$
    g_{\alpha \beta} \; g_{\beta \gamma} \; g_{\gamma \alpha}
    \equiv \operatorname{id}
    \quad \text{on} \
    U_\alpha \cap U_\beta \cap U_\gamma \; ;
\qquad
    g_{\alpha \alpha} 
    \equiv \operatorname{id}
    \quad \text{on} \
    U_\alpha \, .
$$
A complex vector bundle $\mathcal{V}$ over $M$ with typical fiber $V$
 is constructed as 
 the equivalence class of
$
 \coprod_\alpha (U_\alpha \times  V),
$
where
$
 (x,v) \in U_\beta \times V \text{ and } (y,w) \in U_\alpha \times V
$
 are defined to be equivalent if $y = x$ and $w = g_{\alpha \beta}(x) v$.
Then,
 the space of sections $\Gamma(M, \mathcal{V})$ is identified with the collection
\begin{equation}
   \set{(f_\alpha)}{f_\alpha\in C^\infty(U_\alpha, V), 
     \ f_\alpha(x) = g_{\alpha \beta}(x) f_\beta(x), 
    \text{ for } x \in U_\alpha \cap U_\beta} \, .
\label{eqn:9.1.1}
\end{equation}
If $M$ is a complex manifold and if
 every  $g_{\alpha \beta}$ is holomorphic
 (or anti-holomorphic),
 then $\mathcal{V} \to M$ is a holomorphic 
(or anti-holomorphic, respectively)
 vector bundle.
Next, let $G$ be a Lie group,
 $K$ a closed subgroup of $G$,
 and $M:= G/K$ the homogeneous manifold.
Then, we can take an open covering $\{ U_\alpha \}$ of $M$
 such that for each $\alpha$ there is a local section
 $\varphi_\alpha : U_\alpha \to G$ 
 of the principal bundle $G \to G/K$.
Given a representation $\chi : K \to GL_\mathbb{C}(V)$,
 we define the homogeneous vector bundle,
  $\mathcal{V} := G \times_K (\chi, V)$.
Then $\mathcal{V}$ is associated with the transition functions:
$$
   g_{\alpha \beta}: U_\alpha \cap U_\beta \to GL_{\mathbb{C}}(V),
  \quad
   g_{\alpha \beta}(x) := \chi(\varphi_\alpha(x)^{-1}
  \varphi_\beta(x)) \, .
$$
The section space $\Gamma(M, \mathcal{V})$ is identified with
 the following subspace of $C^\infty(G,V)$:
\begin{equation}
   \set{f \in C^\infty(G, V)}{f(g k) = \chi^{-1}(k) f(g), \ \text{for }
 g \in G, k \in K} \, .
\label{eqn:9.2.1}
\end{equation}

\subsection{Pull-back of vector bundles}
\label{subsec:9.3}
Let $G'$ be a Lie group,
 $K'$ a closed subgroup of $G'$,
 and $M' := G'/K'$ the homogeneous manifold.
Suppose that $\sigma : G' \to G$ is a Lie group homomorphism
 such that $\sigma(K') \subset K$. 
We use the same letter $\sigma$ to denote
 by the induced map $M' \to M$,
 $g' K' \mapsto \sigma(g') K$.
Then the pull-back of the vector bundle $\mathcal{V} \to M$,
 denoted by $\sigma^* \mathcal{V} \to M'$,
 is associated to the representation
$$
  \chi \circ \sigma : K' \to GL_{\mathbb{C}}(V) \, .
$$
Then we have a commuting diagram of the pull-back of sections 
(see \eqref{eqn:9.2.1}):  %
\begin{alignat*}{7}
  &\sigma^* 
  &&: {}
  &&\;\;\Gamma(M, \mathcal{V}) 
  &&\to \;\;
  &&\Gamma(M', \sigma^* \mathcal{V})\, ,
  \quad
  &&(f_\alpha)_\alpha 
  &&\mapsto (f_\alpha \circ \sigma)_\alpha \, ,
\\
   &
   && 
   &&{}\qquad \cap
   &&
   && \qquad\cap
   &&
   && 
\\
   &\sigma^* 
   &&: {}
   &&C^\infty(G, V) 
   &&\to \;\;
   &&C^\infty(G', V)\, ,
   \quad
   &&\;\;f 
   &&\mapsto 
    \;\;f \circ \sigma \, .
\end{alignat*}
\subsection{Push-forward of vector bundles}
\label{subsec:9.4}
Suppose that $V$ and $W$ are complex vector spaces
 and that $\xi: V \to W$ is an anti-linear bijective map.
Then,
 we have an anti-holomorphic group isomorphism
$$
   GL_\mathbb{C}(V) \to GL_\mathbb{C}(W)\, ,
   \quad
    g \mapsto g^\xi := \xi \circ g \circ \xi^{-1} \, .
$$
Let $\mathcal{V} \to M$ be a complex vector bundle
 with transition functions
 $g_{\alpha \beta}: U_\alpha \cap U_\beta \to GL_{\mathbb{C}}(V)$.
Then,
 one constructs a complex vector bundle
 $\xi_* \mathcal{V} \to M$
 with the transition functions
 $g_{\alpha \beta}^\xi : U_\alpha \cap U_\beta \to GL_\mathbb{C}(W)$.
We have a natural homomorphism
$$
  \xi_* : \Gamma(M, \mathcal{V}) \to \Gamma(M, \xi_* \mathcal{V})\, ,
   \quad
   (f_\alpha) \mapsto (\xi \circ f_\alpha) \, ,
$$
 which is well-defined because 
 the compatibility condition in \eqref{eqn:9.1.1} %
is satisfied as follows:
If $x \in U_\alpha \cap U_\beta$ then 
$$
  g_{\alpha \beta}^\xi(x) (\xi\circ f_\beta)(x)
 =
  (\xi \circ g_{\alpha \beta}(x)\circ \xi^{-1})  (\xi\circ f_\beta)(x)
 =
  \xi \circ g_{\alpha \beta}(x) f_\beta(x)
 =
  \xi \circ  f_\alpha(x) \, .
$$

If $\mathcal{V}$ is the homogeneous vector bundle $G \times_K (\chi, V)$
 associated to a representation $\chi : K \to GL_{\mathbb{C}}(V)$,
 then
 $\xi_*\mathcal{V}$ is isomorphic to 
the homogeneous vector bundle
$G \times_K (\chi^\xi, W)$
 associated to the representation 
$$
 \chi^\xi : K \to GL_{\mathbb{C}}(W)\, ,
 \quad
 k \mapsto \chi^\xi(k) := \xi \circ \chi(k) \circ \xi^{-1} \, .
$$

\subsection{A sufficient condition for the isomorphism
$\xi_*\sigma^*\mathcal{V}\simeq\mathcal{V}$}
\label{subsec:9.5}
We are particularly interested
 in the case where
 $G' = G$, $K' = K$, 
 $V = W = \mathbb{C}$
   and $\xi(z) := \bar z$
 (the complex conjugate of $z$)
 in the setting of 
 Subsections~\ref{subsec:9.3} and \ref{subsec:9.4}.

By the identification of
 $GL_\mathbb{C}(\mathbb{C})$ with $\mathbb{C}^\times$,
 we have $g^\xi = \overline{g}$
 for $g \in GL_\mathbb{C}(V)\simeq \mathbb{C}^{\times}$.
Then,
 $\chi^\xi$ coincides with the conjugate representation 
$$
  \overline{\chi} : K \to GL_\mathbb{C}(W) \simeq \mathbb{C}^\times\, ,
  \quad
  k \mapsto \overline{\chi(k)}
$$
 for $\chi \in \Hom(K, \mathbb{C}^\times)$.
Thus,
 we have an isomorphism of $G$-equivariant holomorphic line bundles: 
\begin{equation}
 \xi_* \sigma^* \mathcal{V} \simeq G 
   \times_K (\overline{\chi\circ \sigma}, \mathbb{C})
\label{eqn:9.5.1}
\end{equation}
with the following correspondence of sections:
$$
  \xi_* \circ \sigma^* :
  \Gamma(M, \mathcal{V}) \to \Gamma(M, \xi_* \sigma^* \mathcal{V})\, ,
  \quad
   (f_\alpha) \mapsto (\overline{f_\alpha \circ \sigma}) \, .
$$

We now apply the formula \eqref{eqn:9.5.1} 
to the setting where
 $M=G/K$ is an irreducible Hermitian symmetric space.

\begin{lemma}
\label{lem:9.6}
Let $\chi : K \to \mathbb{C}^\times$ be a unitary character.
We denote by $\mathcal{V}$ the homogeneous line bundle $G \times_K (\chi, \mathbb{C})$.
Suppose $\sigma$ is an involutive automorphism of $G$
of anti-holomorphic type (see Definition~\ref{def:holo-anti}).
Then we have an isomorphism of $G$-equivariant holomorphic line bundles:
$
   \xi_* \sigma^* \mathcal{V} \simeq \mathcal{V}. 
$
\end{lemma}

\begin{proof}
In view of \eqref{eqn:9.5.1}, %
 it suffices to show
$
   \overline{\chi \circ \sigma} = \chi.
$
As   %
 the character $\chi$ of $K$ is unitary,
 we have $\overline{\chi(k)} = \chi(k^{-1})$
 for any $k \in K$.
Let $Z$ be a generator of the center 
$\mathfrak{c}(\mathfrak{k})$ of $\mathfrak{k}$.
Since $\sigma$ is of anti-holomorphic type,
we have $\sigma Z = -Z$,
and then
$$
  \overline{\chi \circ \sigma(\exp t Z)}
 =
  \overline{\chi (\exp (- t Z))}
 =
   \chi (\exp t Z)
 \qquad
 (t \in \mathbb{R})\, .
$$
On the other hand,
 if $k \in [K, K]$,
 then $\overline{\chi \circ \sigma(k)} = 1 = \chi(k)$
because $[K,K]$ is a connected semisimple Lie group.
As $K= \exp \mathfrak{c}(\mathfrak{k}) \cdot [K,K]$,
 we have shown
$
   \overline{\chi \circ \sigma} = \chi.
$
Hence Lemma.
\qed
\end{proof}

\begin{acknowledgement}

I owe much to J. Faraut
 for enlightening discussions,
 in particular,
 for clarifying the idea of his joint work
  \cite{xft}  with E. Thomas at an early stage 
  on the occasion of
 the CIMPA school on Invariant Theory
 in Tunisia 
 organized by P. Torasso in 1996,
 when I formalized 
 Theorem~\ref{thm:A}
 and gave its proof in the classical case.

Some of the results of this article 
 were presented in various places
 including the conferences
 \lq\lq Representation Theory at Saga\rq\rq\ 
 in 1997
 (see  \cite{xkmfjp})
  organized by K. Mimachi, at
 \lq\lq Summer Solstice Days at Paris University VI and 
VII\rq\rq\ in 1999,
 at Copenhagen 
organized by B. \O rsted and H. Schlichtkrull
in 1999,
 and at Oberwolfach
 organized by
  A. Huckleberry and K.-H. Neeb and J. Wolf 
  in 2000 and 2004, at
the MSRI program ``Integral Geometry'' organized by S. Gindikin,
L. Barchini, and R. Zierau in 2001,
at the workshop
``Representation theory of Lie groups, harmonic analysis on
 homogeneous spaces and quantization''
at the Lorentz center in Leiden in  2002
(see \cite{xkleiden})
 organized by G. van Dijk and
V. F. Molchanov, 
 and
 at  the workshop ``Representation theory
  and automorphic forms''
in Korea in 2005 coorganized by J. Yang
 and W. Schmid with the author.
I express my deep gratitude to the organizers
 of these conferences
and to the participants for 
helpful and stimulating 
comments on various occasions.
I also thank an anonymous referee who read the manuscript very
 carefully. 
Special thanks are due to Ms. Suenaga
for her help in preparing for
the final manuscript.

\end{acknowledgement}

\printindex
\end{document}